\documentclass[12pt,twoside,reqno,leqno]{amsart}
\usepackage{url}
\usepackage{amsmath}
\usepackage{amssymb}
\usepackage{amsfonts}
\usepackage{amsthm}
\usepackage[mathscr]{eucal}
\usepackage{amscd}
\usepackage{amsmath}
\usepackage{latexsym}
\usepackage[dvips]{graphics}
\usepackage{graphics,epstopdf}
\usepackage[pdftex]{graphicx}
\theoremstyle{plain}

\newtheorem{prop}{Proposition}
\theoremstyle{definition}

 \theoremstyle{remark}
\newtheorem{rem}{Remark}
\newtheorem{example}{Example}

\allowdisplaybreaks 
\def\DJ{\leavevmode\setbox0=\hbox{D}\kern0pt\rlap
{\kern.04em\raise.188\ht0\hbox{-}}D}
\begin{document}
\title[Unified Multi-Tupled Fixed Point Theorems] {Unified Multi-Tupled Fixed Point Theorems involving
Mixed Monotone Property in Ordered Metric Spaces}
\author[Alam, Imdad and Ali]{Aftab Alam, Mohammad Imdad and Javid Ali$^{\ast}$}
\thanks{$^\ast$Correspondence:\;{\rm javid@amu.ac.in}\bigskip}\maketitle
\begin{center}
{\footnotesize Department of Mathematics\\Aligarh Muslim
University\\ Aligarh-202002, India.\\
Email addresses: aafu.amu@gmail.com, mhimdad@gmail.com,
javid@amu.ac.in\\}
\end{center}
{\footnotesize{\noindent {\bf Abstract.}} In the present article, we
introduce a unified notion of multi-tupled fixed points and utilize
the same to prove some existence and uniqueness unified multi-tupled
fixed point theorems for Boyd-Wong type nonlinear contractions satisfying generalized mixed monotone property in ordered
metric spaces. Our results unify several classical and well-known
$n$-tupled (including coupled, tripled and quadrupled ones) fixed
point results existing in the literature.


}
\section{Introduction and Preliminaries}
\label{SC:Introduction} In our entire presentation, we use throughout the following symbols
and notations.
\begin{enumerate}
\item [{(1)}] $\mathbb{N}_0$ stands for the set of nonnegative
integers ($i.e.$ $\mathbb{N}_0=\mathbb{N}\cup \{0\}$).
\item [{(2)}] $m,l \in \mathbb{N}_0$.
\item [{(3)}] $n$ stands for a fixed natural number greater than
 1.
\item [{(4)}] $I_n$ denotes the set $\{1,2,...,n\}$ and we use $i,j,k\in I_n.$
\item [{(5)}] $\iota_n$ denotes a fixed nontrivial partition $\{A,B\}$ of $I_n$ ($i.e.$ $\iota_n:=\{A,B\}$, where $A$ and $B$ are nonempty subsets of $I_n$ such that $A\cup B=I_n$ and
$A\cap B=\emptyset$).
\item [{(6)}] As usual, for a nonempty set $X$, $X^n$ denotes the cartesian
product of $n$ identical copies of $X$, $i.e.$, $X^n:=X\times X
\times\stackrel{(n)}{...}\times X$. We often call $X^n$ as the
$n$-dimensional product set induced by $X$.
\item [{(7)}] A sequence in $X$ is denoted by $\{x^{(m)}\}$ and a
sequence in $X^n$ is denoted by $\{{\rm U}^{(m)}\}$ where
U$^{(m)}=(x^{(m)}_1,x^{(m)}_2,...,x^{(m)}_n)$ such that for each
$i\in I_n$, $\{x^{(m)}_i\}$ is a sequence in $X$.
\end{enumerate}
Naturally, given a mapping $f:X\rightarrow X$, an element $x\in X$ satisfying $f(x)=x$ is called a fixed point of $f$. Presic \cite{P1,P2} extends
the notion of fixed points for the mapping $F:X^n\rightarrow X$ as follows:\\

\noindent\textbf{Definition 1} \cite{P1,P2}. Let $X$ be a nonempty
set and $F:X^{n} \rightarrow X$ a mapping. An element $x\in X$  is
called a fixed point of $F$ if
$$F(x,x,...,x)=x.$$
In 1975, particularly for $n=2$, Opoitsev \cite{C01,C02} initiated a
weaker notion of fixed point, which satisfies $F(x,y)=x$ and
$F(y,x)=y$ instead of $F(x,x)=x$ and hence for $y=x$, this reduces
to Definition 1. Using this notion, Opoitsev and Khurodze \cite{C03}
proved some results for nonlinear operators on ordered Banach
spaces. Unknowingly, in 1987, Guo and Lakshmikantham \cite{C04}
reconsidered this concept for mixed monotone operators defined on a
real Banach space equipped with a partial ordering by a cone besides
calling this notion as coupled fixed point.\\

\noindent\textbf{Definition 2} \cite{C01}-\cite{C04}. Let $X$ be a
nonempty set and $F:X^{2} \rightarrow X$ a mapping. An element
$(x,y)\in X$ is called a coupled fixed point of $F$ if
$$F(x,y)=x\;{\rm and}\;F(y,x)=y.$$

Inspired by  the results of Guo and Lakshmikantham \cite{C04}, several authors ($e.g.$
\cite{C05}-\cite{C14}) studied and developed the theory of coupled
fixed points for mixed monotone operators in the context of ordered Banach spaces.\\

Recall that a set $X$ together with a partial order $\preceq$ (often
denoted by $(X,\preceq)$) is called an ordered set. In this context,
$\succeq$ denotes the dual order of $\preceq$ ($i.e.\;x\succeq y$
means $y\preceq x$). Two elements $x$ and $y$ in an ordered set
$(X,\preceq)$ are said to be comparable if either $x\preceq y$ or
$y\preceq x$ and denote it as $x\prec\succ y$. In respect of a pair
of self-mappings $f$ and $g$ defined on an ordered set
$(X,\preceq)$, we say that $f$ is $g$-increasing (resp.
$g$-decreasing) if for any $x,y\in X$; $g(x)\preceq g(y)$ implies
$f(x)\preceq f(y) (\rm{resp.}\; f(x)\succeq f(y))$. As per standard
practice, $f$ is called $g$-monotone if $f$ is either $g$-increasing
or $g$-decreasing. Notice that with $g=I$ (the identity mapping),
the notions of $g$-increasing, $g$-decreasing and $g$-monotone
mappings transform into increasing, decreasing and
monotone mappings respectively.\\

In 2006, Bhasker and Lakshmikantham \cite{C1} extended the idea of
monotonicity for the mapping $F:X^{2} \rightarrow X$ by introducing
the notion of mixed monotone property in ordered metric spaces and
obtained some coupled fixed point theorems for linear contractions
satisfying mixed monotone property with application in existence and
uniqueness of a solution of periodic boundary value problems.
Although, some variants of such results were earlier reported in 2001 by Zhang \cite{C13}.\\

\noindent\textbf{Definition 3} \cite{C1}. Let $(X,\preceq)$ be an
ordered set and $F: X^{2} \rightarrow X$ a mapping. We say that $F$
has mixed monotone property if $F$ is increasing in its first
argument and is decreasing in its second argument, $i.e.$, for any
$x,y\in X,$
$$x_{1}, x_{2}\in X, x_{1}\preceq x_{2}\Rightarrow F(x_{1},y)\preceq F(x_{2},y),$$
$$y_{1}, y_{2}\in X, y_{1}\preceq y_{2} \Rightarrow F(x,y_{1})\succeq
F(x,y_{2}).$$ Later, Lakshmikantham and \'{C}iri\'{c} \cite{C2}
generalized the notions of coupled fixed point and mixed monotone
property for a pair of mappings, which
runs as follows:\\

\noindent\textbf{Definition 4} \cite{C2}. Let $X$ be a nonempty set
and $F:X^{2} \rightarrow X$ and $g: X\rightarrow X$ two mappings. An
element $(x,y)\in X$ is called a coupled coincidence point of $F$
and $g$ if
$$F(x,y)=g(x)\;{\rm and}\;F(y,x)=g(y).$$

\noindent\textbf{Definition 5} \cite{C2}. Let $(X,\preceq)$ be an
ordered set and $F: X^{2} \rightarrow X$ and $g: X\rightarrow X$ two
mappings. We say that $F$ has mixed $g$-monotone property if $F$ is
$g$-increasing in its first argument and is $g$-decreasing in its
second argument, $i.e.$, for any $x,y\in X,$
$$x_{1}, x_{2}\in X, g(x_{1})\preceq g(x_{2})\Rightarrow F(x_{1},y)\preceq F(x_{2},y),$$
$$y_{1}, y_{2}\in X, g(y_{1})\preceq g(y_{2}) \Rightarrow F(x,y_{1})\succeq
F(x,y_{2}).$$ Notice that under the restriction $g=I,$ the identity
mapping on $X,$ Definitions 4 and 5 reduce to Definitions 2 and 3 respectively.\\

As a continuation of this trends, various authors extended the
notion of coupled fixed (coincidence) point and mixed monotone
($g$-monotone) property for the mapping $F:X^n\rightarrow X,\;n\geq
3$ in different ways. A natural extension of mixed monotone property
introduced by Berinde and Borcut \cite{T1} (for $n=3$), Karapinar
and Luong \cite{Q1} (for $n=4$), Imdad $et\;al.$ \cite{n1} (for even
$n$) and Gordji and Ramezani \cite{NX1} and Ert$\ddot{\rm
u}$rk and Karakaya \cite{NX4,NX4-} (for general $n$) runs as follows:\\

\noindent\textbf{Definition 6} (\cite{T1}-\cite{NX4}). Let
$(X,\preceq)$ be an ordered set and $F: X^{n} \rightarrow X$ a
mapping. We say that $F$ has alternating mixed monotone property if
$F$ is increasing in its odd position argument and is decreasing in
its even position argument, $i.e.$, for any $x_1,x_2,...,x_n\in X,$
$$\underline{x}_1,\overline{x}_1\in X, \underline{x}_1\preceq \overline{x}_1\Rightarrow F(\underline{x}_1,x_{2},...,x_{n})\preceq
F(\overline{x}_1,x_{2},...,x_{n})$$
$$\underline{x}_2,\overline{x}_2\in X, \underline{x}_2\preceq
\overline{x}_2\Rightarrow F(x_{1},\underline{x}_2,...,x_{n})\succeq
F(x_{1},\overline{x}_2,...,x_{n})$$
\indent\hspace{6cm}$\vdots$\\
\indent\hspace{2.3cm}$\underline{x}_n,\overline{x}_n\in X,
\underline{x}_n\preceq \overline{x}_n \Rightarrow
{\begin{cases}F(x_{1},x_{2},...,\underline{x}_n)\preceq
F(x_{1},x_{2},...,\overline{x}_n)\;{\rm if\;} n\;{\rm is\; odd,}\cr
\hspace{0.0in}F(x_{1},x_{2},...,\underline{x}_n)\succeq
F(x_{1},x_{2},...,\overline{x}_n)\;{\rm if\;}n\;{\rm is\;
even.}\cr\end{cases}}$\\
Although, the authors in \cite{T1}-\cite{NX4} used the word `mixed
monotone property', but we use `alternating mixed monotone property'
to differ another extension of
mixed monotone property (see Definition 7).\\

Another extension of Definition 3 is $p$-monotone property
introduced by Berzig and Samet \cite{HD1} as follows:\\

\noindent\textbf{Definition 7} \cite{HD1}. Let $(X,\preceq)$ be an
ordered set, $F: X^{n} \rightarrow X$ a mapping and $1\leq p<n$. We
say that $F$ has $p$-mixed monotone property if $F$ is increasing
for the range of components from 1 to $p$ and is decreasing for the
range of components from $p+1$ to $n$, $i.e.$, for any
$x_1,x_2,...,x_n\in X,$
$$\underline{x}_1,\overline{x}_1\in X,\;\underline{x}_1\preceq
\overline{x}_1\Rightarrow
F(\underline{x}_1,x_{2},...,x_{p},...,x_{n})\preceq
F(\overline{x}_1,x_{2},...,x_{p},...,x_{n})$$
$$\underline{x}_2,\overline{x}_2\in X,\;\underline{x}_2\preceq\overline{x}_2\Rightarrow F(x_{1},\underline{x}_2,...,x_{p},...,x_{n})\preceq
F(x_{1},\overline{x}_2,...,x_{p},...,x_{n})$$
\indent\hspace{6cm}$\vdots$
$$\underline{x}_p,\overline{x}_p\in X,\; \underline{x}_p\preceq
\overline{x}_p\Rightarrow
F(x_{1},x_{2},...,\underline{x}_p,...,x_{n})\preceq
F(x_{1},x_{2},...,\overline{x}_p,...,x_{n})$$
\indent\hspace{1cm}$\underline{x}_{p+1},\overline{x}_{p+1}\in X,\;
\underline{x}_{p+1}\preceq \overline{x}_{p+1}\Rightarrow
F(x_{1},...,x_{p},\underline{x}_{p+1},...,x_{n})\succeq
F(x_{1},...,x_{p},\overline{x}_{p+1},...,x_{n})$\\
\indent\hspace{1cm}$\underline{x}_{p+2},\overline{x}_{p+2}\in X,\;
\underline{x}_{p+2}\preceq \overline{x}_{p+2}\Rightarrow
F(x_{1},...,x_{p+1},\underline{x}_{p+2},...,x_{n})\succeq
F(x_{1},...,x_{p+1},\overline{x}_{p+2},...,x_{n})$\\
\indent\hspace{6cm}$\vdots$\\
\indent\hspace{1cm}$\underline{x}_n,\overline{x}_n\in X,\;\underline{x}_n\preceq \overline{x}_n\Rightarrow F(x_{1},x_{2},...,x_{p},...,\underline{x}_n)\succeq F(x_{1},x_{2},...,x_{p},...,\overline{x}_n)$\\

In 2012, Rold$\acute{\rm a}$n $et\;al.$ \cite{MD1,MD2} introduced a
generalized notion of mixed
monotone property. Although, the authors of \cite{MD2}, termed the same as `mixed
monotone property (w.r.t. $\{{\rm A,B}\}$)'. For the sack of brevity, we prefer to call the same as
`$\iota_n$-mixed
monotone property'.\\

\noindent\textbf{Definition 8} (see \cite{MD1,MD2}). Let
$(X,\preceq)$ be an ordered set and $F:X^{n} \rightarrow X$ a
mapping. We say that $F$ has $\iota_n$-mixed monotone property if
$F$ is increasing in arguments
of A and is decreasing in arguments of B, $i.e.$, for any $x_1,x_2,...,x_n\in X$, \\
\indent$\underline{x}_i,\overline{x}_i\in X$,
$\underline{x}_i\preceq \overline{x}_i$ \\
\indent\hspace{3mm}$\Rightarrow
F(x_{1},x_{2},...,x_{i-1},\underline{x}_i,x_{i-1},...,x_{n})\preceq
F(x_{1},x_{2},...,x_{i-1},\overline{x}_i,x_{i-1},...,x_{n})$ for each $i\in$ A,\\
\indent$\underline{x}_i,\overline{x}_i\in X$,
$\underline{x}_i\preceq \overline{x}_i$\\
\indent\hspace{3mm}$\Rightarrow
F(x_{1},x_{2},...,x_{i-1},\underline{x}_i,x_{i-1},...,x_{n})\succeq
F(x_{1},x_{2},...,x_{i-1},\overline{x}_i,x_{i-1},...,x_{n})$ for each $i\in$ B.\\

In particular, on setting $\iota_n:=\{{\rm A,B}\}$ such that
A$=\{2s-1:s\in\{1,2,...,[\frac{n+1}{2}]\}\}\; i.e.$ the set of all
odd numbers in $I_n$ and B$=\{2s:s\in\{1,2,...,[\frac{n}{2}]\}\}\;
i.e.$ the set of all even numbers in $I_n$, Definition 8 reduces to
the definition of alternating mixed monotone property, while on
setting $\iota_n:=\{{\rm A,B}\}$ such that A$=\{1,2,...,p\}$ and
B$=\{p,p+1,...,n\}$, where $1\leq p<n$, Definition 8 reduces to
the definition of $p$-mixed monotone property.\\

\noindent\textbf{Definition 9} (see \cite{MD1}). Let $(X,\preceq)$
be an ordered set and $F:X^{n} \rightarrow X$ and $g:X\rightarrow X$
two mappings. We say that $F$ has $\iota_n$-mixed $g$-monotone
property if $F$ is $g$-increasing in arguments
of A and is $g$-decreasing in arguments of B, $i.e.$, for any $x_1,x_2,...,x_n\in X$, \\
\indent$\underline{x}_i,\overline{x}_i\in X$,
$g(\underline{x}_i)\preceq g(\overline{x}_i)$ \\
\indent\hspace{3mm}$\Rightarrow
F(x_{1},x_{2},...,x_{i-1},\underline{x}_i,x_{i-1},...,x_{n})\preceq
F(x_{1},x_{2},...,x_{i-1},\overline{x}_i,x_{i-1},...,x_{n})$ for each $i\in$ A,\\
\indent$\underline{x}_i,\overline{x}_i\in X$,
$g(\underline{x}_i)\preceq g(\overline{x}_i)$\\
\indent\hspace{3mm}$\Rightarrow
F(x_{1},x_{2},...,x_{i-1},\underline{x}_i,x_{i-1},...,x_{n})\succeq
F(x_{1},x_{2},...,x_{i-1},\overline{x}_i,x_{i-1},...,x_{n})$ for each $i\in$ B.\\
Notice that under the restriction $g=I,$ the identity
mapping on $X,$ Definition 9 reduces to Definition 8.\\

In the same continuation
Paknazar $et\;al.$ \cite{NX3} introduced the concept of new
$g$-monotone property for the mapping $F:X^{n} \rightarrow X$, which merely
depends on the first argument of $F$. Thereafter, Karapinar
$et\;al.$ \cite{NX3+} noticed that multi-tupled coincidence theorems
involving new $g$-monotone property (proved by Paknazar $et\;al.$
\cite{NX3}) can be reduced to
corresponding (unidimensional) coincidence theorems.\\

In an attempt to extend the notion of coupled fixed point from $X^2$
to $X^3$ and $X^4$ various authors introduced the concepts of
tripled and quadrupled fixed points respectively. Here it can be
pointed out that these notions were defined in different ways by their
respective authors so as to make their notions compatible under the corresponding mixed monotone
property. The following definitions of tripled and quadrupled fixed points
are available in literature.\\

\noindent\textbf{Definition 10}. Let $X$ be a nonempty set and
$F:X^{3} \rightarrow X$ a mapping. An element $(x_1,x_2,x_3)\in X^3$
is called a tripled/triplet fixed point of $F$ if
\begin{enumerate}
\item [{$\bullet$}](Berinde and Borcut \cite{T1}) $F(x_{1},x_{2},x_3)=x_{1},\;F(x_{2},x_{1},x_2)=x_{2},\;F(x_{3},x_{2},x_1)=x_{3}.$
\item [{$\bullet$}](Wu and Liu \cite{TQ}) $F(x_{1},x_{2},x_3)=x_{1},\;F(x_{2},x_{3},x_2)=x_{2},\;F(x_{3},x_{2},x_1)=x_{3}.$
\item [{$\bullet$}](Berzig and Samet \cite{HD1})
$F(x_{1},x_{2},x_3)=x_{1},\;F(x_{2},x_{1},x_3)=x_{2},\;F(x_{3},x_{3},x_2)=x_{3}.$\\
\end{enumerate}

\noindent\textbf{Definition 11}. Let $X$ be a nonempty set and
$F:X^{4} \rightarrow X$ a mapping. An element $(x_1,x_2,x_3,x_4)\in
X^4$ is called a quadrupled/quartet fixed point of $F$ if
\begin{enumerate}
\item [{$\bullet$}](Karapinar and Luong \cite{Q1}) $F(x_{1},x_{2},x_3,x_4)=x_{1},\;F(x_{2},x_{3},x_4,x_1)=x_{2},$\\
\indent\hspace{4.6cm}
$F(x_{3},x_{4},x_1,x_2)=x_{3},\;F(x_{4},x_{1},x_2,x_3)=x_{4}.$
\item [{$\bullet$}](Wu and Liu \cite{TQ}) $F(x_{1},x_{4},x_3,x_2)=x_{1},\;F(x_{2},x_{1},x_4,x_3)=x_{2},$\\
\indent\hspace{3cm}$F(x_{3},x_{2},x_1,x_4)=x_{3},\;F(x_{4},x_{3},x_2,x_1)=x_{4}.$
\item [{$\bullet$}](Berzig and Samet \cite{HD1}) $F(x_{1},x_{2},x_3,x_4)=x_{1},\;F(x_{1},x_{2},x_4,x_3)=x_{2},$\\
\indent\hspace{4cm}$F(x_{3},x_{4},x_2,x_1)=x_{3},\;F(x_{3},x_{4},x_1,x_2)=x_{4}.$
\end{enumerate}

In the same continuation, the notion of coupled fixed point is
extended for the mapping $F:X^n\rightarrow X$ by various authors in
different ways (similar to tripled and quadrupled ones). Also this
notion is available under different names as adopted by various
authors in their respective papers such as:
\begin{enumerate}
\item [{$\bullet$}] {\em $n$-tupled fixed point} (see Imdad $et\;al.$ \cite{n1})
\item [{$\bullet$}] {\em $n$-tuple fixed point} (see Karapinar and Rold$\acute{\rm a}$n \cite{NX4+}, Al-Mezel $et\;al.$ \cite{MD3}, Rad $et\;al.$ \cite{NX5})
\item [{$\bullet$}] {\em $n$-tuplet fixed point} (see Ert$\ddot{\rm u}$rk and Karakaya \cite{NX4,NX4-})
\item [{$\bullet$}] {\em $n$-fixed point} (see Gordji and Ramezani \cite{NX1}, Paknazar $et\;al.$ \cite{NX3})
\item [{$\bullet$}] {\em Fixed point of $n$-order} (see Samet and Vetro \cite{NFI}, Berzig and Samet \cite{HD1})
\item [{$\bullet$}] {\em Multidimensional fixed point} (see Rold$\acute{\rm a}$n $et\;al.$ \cite{MD1}, Dalal $et\;al.$ \cite{R2})
\item [{$\bullet$}] {\em Multiplied fixed point} (see Olaoluwa and Olaleru \cite{NX6})
\item [{$\bullet$}] {\em Multivariate coupled fixed point} (see Lee and Kim \cite{NX7}).
\end{enumerate}
Here, it is worth mentioning that we prefer to use `$n$-tupled
fixed point' due to its natural analogy with earlier used terms namely: coupled ($2$-tupled), tripled ($3$-tupled) and quadrupled ($4$-tupled).\\

After the appearance of multi-tupled fixed points, some authors paid
attention to unify the different types of multi-tupled fixed points.
A first attempt of this kind was given by Berzig and Samet
\cite{HD1}, wherein authors defined a one-to-one correspondence
between alternating mixed monotone property and $p$-mixed monotone
property and utilized the same to define a unified notion of $n$-tupled
fixed point by using 2$n$ mappings from $I_n$ to $I_n$. Later,
Rold$\acute{\rm a}$n $et\;al.$ \cite{MD1} extended the notion of
$n$-tupled fixed point of Berzig and Samet \cite{HD1} so as to
make $\iota_n$-mixed monotone property working and introduced the
notion of $\Upsilon$-fixed point based on $n$ mappings from $I_n$ to
$I_n$. To do this, Rold$\acute{\rm a}$n $et\;al.$ \cite{MD1}
considered the following family
$$\Omega_{{\rm A,B}}:=\{\sigma:I_n\rightarrow I_n:\sigma(A)\subseteq A\;{\rm and}\;\sigma(B)\subseteq
B\}$$ and
$$\Omega^{\prime}_{{\rm A,B}}:=\{\sigma:I_n\rightarrow I_n:\sigma(A)\subseteq B\;{\rm and}\;\sigma(B)\subseteq A\}.$$
Let $\sigma_1,\sigma_2,...,\sigma_n$ be $n$ mappings from $I_n$ into
itself and let $\Upsilon$ be $n$-tuple
$(\sigma_1,\sigma_2,...,\sigma_n)$.\\

\noindent\textbf{Definition 12} \cite{MD1,MD2}. Let $X$ be a
nonempty set and $F:X^{n} \rightarrow X$ a mapping. An element
$(x_1,x_2,...,x_n)\in X^n$ is called a $\Upsilon$-fixed point of $F$
if
$$F(x_{\sigma_i(1)},x_{\sigma_i(2)},...,x_{\sigma_i(n)})=x_i\;\;\forall~i\in I_n.$$

\begin{rem} (\cite{NX4+,MD3}).
In order to ensure the existence of $\Upsilon$-coincidence/fixed
points, it is very important to assume that the $\iota_n$-mixed
$g$-monotone property is compatible with the permutation of the
variables, $i.e.$, the mappings of
$\Upsilon=(\sigma_1,\sigma_2,...,\sigma_n)$ should verify:
$$\sigma_i\in \Omega_{{\rm A,B}}\;{\rm if}\;i\in{\rm A\;\;and}\;\;\sigma_i\in \Omega^{\prime}_{{\rm A,B}}\;{\rm if}\;i\in{\rm B}.$$
\end{rem}

In this paper, we observe that the $n$-mappings involved in
$\Upsilon$-fixed point are not independent to each other. We can
represent these mappings in the form of only one mapping, which is
infact a binary operation on $I_n$. Using this fact, we refine and
modify the notion of $\Upsilon$-fixed point and introduce the notion
of $\ast$-fixed point, where $\ast$ is a binary operation on $I_n$.
Moreover, we represent the binary operation $\ast$ in the form of a
matrix. Due to this, the notion of $\ast$-fixed point becomes
relatively more natural and effective as compared to $\Upsilon$-fixed
point. Furthermore, we present some $\ast$-coincidence theorems for
a pair of mappings $F:X^{n} \rightarrow X$ and $g:X\rightarrow X$
under Boyd-Wong type nonlinear contractions satisfying
$\iota_n$-mixed $g$-monotone property in ordered metric spaces. Our
results unify several multi-tupled fixed/coincidence point results
of the existing literature.

\section{Ordered Metric Spaces and Control Functions}
\label{SC: Ordered Metric Spaces and Control Functions}

In this section, we summarize some order-theoretic metrical notions
and possible relations between some existing control functions
besides indicating a recent coincidence theorem for nonlinear
contractions in ordered metric spaces. Here it can be pointed out that major part of the present work is essentially contained in \cite{PGF13,PGF17,PGF18}.\\

\noindent\textbf{Definition 13} \cite{PGF2}. A triplet
$(X,d,\preceq)$ is called an ordered metric space if $(X,d)$ is a
metric space and $(X,\preceq)$ is an ordered set. Moreover, if
$(X,d)$ is a complete metric space, we say that $(X,d,\preceq)$ is
an ordered complete metric space.\\

\noindent\textbf{Definition 14} \cite{PGF18}. Let $(X,d,\preceq)$ be an ordered metric space. A nonempty subset $Y$ of $X$ is
called a subspace of $X$ if $Y$ itself is an ordered metric space equipped with the metric $d_{Y}$ and partial order $\preceq_{Y}$ defined by:\\
$$d_{Y}(x,y)=d(x,y)\;\forall~x,y\in Y$$ and
$$x \preceq_{Y} y\Leftrightarrow x\preceq y\;\forall~x,y\in Y.$$

As per standard practice, we can define the notions of increasing,
decreasing, monotone, bounded above and bounded below sequences
besides bounds (upper as well as lower) of a sequence in an ordered
set $(X,\preceq)$, which on the set of real numbers with natural
ordering coincide with their usual senses (see Definition 8
\cite{PGF13}). Let $(X,d,\preceq)$ be an ordered metric space and
$\{x^{(m)}\}$ a sequence in $X$. We adopt the following notations:
\begin{enumerate}
\item[{(i)}] if $\{x^{(m)}\}$ is increasing and $x^{(m)}\stackrel{d}{\longrightarrow} x$ then we denote it symbolically by $x^{(m)}\uparrow x,$
\item[{(ii)}] if $\{x^{(m)}\}$ is decreasing and $x^{(m)}\stackrel{d}{\longrightarrow} x$ then we denote it symbolically
by $x^{(m)}\downarrow x,$
\item[{(iii)}] if $\{x^{(m)}\}$ is monotone and $x^{(m)}\stackrel{d}{\longrightarrow} x$ then we denote it symbolically
by $x^{(m)}\uparrow\downarrow x.$\\
\end{enumerate}
\noindent\textbf{Definition 15} \cite{PGF17}. An ordered metric
space $(X,d,\preceq)$ is called O-complete if every monotone Cauchy
sequence in $X$ converges.
\begin{rem} \cite{PGF17} Every ordered complete metric space is O-complete.\end{rem}

\noindent\textbf{Definition 16} \cite{PGF18}. Let $(X,d,\preceq)$ be
an ordered metric space. A subset $E$ of $X$ is called O-closed if
for any sequence $\{x_n\}\subset E$,
$$x_n \uparrow\downarrow x\Rightarrow x\in E.$$

\begin{rem} \cite{PGF18} Every closed subset of an
ordered metric space is O-closed.\end{rem}

\begin{prop} \cite{PGF18} Let $(X,d,\preceq)$
be an {\rm O}-complete ordered metric space.  A subspace $E$ of $X$
is {\rm O}-closed iff $E$ is {\rm O}-complete.\end{prop}

\noindent\textbf{Definition 17} \cite{PGF17}. Let $(X,d,\preceq)$ be
an ordered metric space, $f:X\rightarrow X$ a mapping and $x\in X$.
Then $f$ is called O-continuous at $x$ if for any sequence
$\{x^{(m)}\}\subset X$,
$$x^{(m)} \uparrow\downarrow x\Rightarrow f(x^{(m)})\stackrel{d}{\longrightarrow} f(x).$$
Moreover, $f$ is called O-continuous if it is O-continuous at each point of $X$.\\
\begin{rem} \cite{PGF17} Every continuous mapping defined on an ordered metric
space is O-continuous.\end{rem}

{\noindent\textbf{Definition 18} \cite{PGF17}. Let $(X,d,\preceq)$
be an ordered metric space, $f$ and $g$ two self-mappings on $X$ and
$x\in X$. Then $f$ is called $(g,{\rm O})$-continuous at $x$ if for
any sequence $\{x^{(m)}\}\subset X$,
$$g(x^{(m)}) \uparrow\downarrow g(x)\Rightarrow f(x^{(m)})\stackrel{d}{\longrightarrow} f(x).$$
Moreover, $f$ is called $(g,{\rm O})$-continuous if it is $(g,{\rm
O})$-continuous at each point of $X$.\\

\noindent\textbf{Definition 19} \cite{PGF17}. Let $(X,d,\preceq)$ be
an ordered metric space and $f$ and $g$ two self-mappings on $X$. We
say that $f$ and $g$ are O-compatible if for any sequence
$\{x^{(m)}\}\subset X$ and for any $z\in X$,
$$g(x^{(m)})\uparrow\downarrow z\;{\rm and}\;f(x^{(m)})\uparrow\downarrow z\Rightarrow\lim\limits_{m\to \infty}d(gfx^{(m)},fgx^{(m)})=0.$$

Notice that the above notion is slightly weaker than the notion of
$O$-compatibility (of Luong and Thuan \cite{CP5}) as they
\cite{CP5} assumed that only the sequence $\{gx^{(m)}\}$ is monotone
but here both
$\{gx^{(m)}\}$ and $\{fx^{(m)}\}$ be assumed monotone.\\

The following notion is formulated by using certain properties on
ordered metric space (in order to avoid the necessity of continuity
requirement on underlying mapping)
utilized by earlier authors especially from \cite{PF2, PGF3, C1, C2} besides some other ones.\\

\noindent\textbf{Definition 20} \cite{PGF13}. Let $(X,d,\preceq)$
be an ordered metric space and $g$ a self-mapping on $X.$ We say
that
\begin{enumerate}
\item[{(i)}] $(X,d,\preceq)$ has {\it
g-ICU}\;(increasing-convergence-upper bound) property if $g$-image
of every increasing convergent sequence $\{x_n\}$ in $X$ is bounded
above by $g$-image of its limit (as an upper bound), $i.e.,$
$$x_n\uparrow x \Rightarrow g(x_n)\preceq g(x)\;\;\forall~ n\in \mathbb{N}_0
,$$
\item[{(ii)}] $(X,d,\preceq)$ has {\it g-DCL}\;(decreasing-convergence-lower bound) property if
$g$-image of every decreasing convergent sequence $\{x_n\}$ in $X$
is bounded below by $g$-image of its limit (as a lower bound),
$i.e.,$
$$x_n\downarrow x \Rightarrow g(x_n)\succeq g(x)\;\;\forall~ n\in \mathbb{N}_0\;{\rm and}$$
\item[{(iii)}] $(X,d,\preceq)$ has {\it g-MCB}\;(monotone-convergence-boundedness)
property if it has both {\it g-ICU} as well as {\it g-DCL} property.
\end{enumerate}
Notice that under the restriction $g=I,$ the identity mapping on $X,$ the notions of {\it g-ICU} property, {\it g-DCL} property and {\it g-MCB} property transform to
{\it ICU} property, {\it DCL} property and {\it MCB} property respectively.\\

The following family of control functions is essentially due to Boyd
and Wong \cite{B2}.
$$\Psi=\Big\{\varphi:[0,\infty)\to [0,\infty):\varphi(t)<t\;{\rm for~each}\; t>0\;{\rm and~}\varphi\;{\rm is~
right-upper~semicontinuous}\Big\}.$$ Mukherjea \cite{B3} introduced
the following family of control functions:
$$\Theta=\Big\{\varphi:[0,\infty)\to [0,\infty):\varphi(t)<t\;{\rm for~each}\; t>0\;{\rm and~}\varphi\;{\rm is~
right~ continuous}\Big\}.$$ The following family of control
functions found in literature is more natural.
$$\Im=\Big\{\varphi:[0,\infty)\to
[0,\infty):\varphi(t)<t\;{\rm for~each}\; t>0\;{\rm and~}
\varphi\;{\rm is~ continuous}\Big\}.$$ The following family of
control functions is due to Lakshmikantham and \'{C}iri\'{c} \cite{C2}.\\
$$\Phi=\Big\{\varphi:[0,\infty)\to [0,\infty):\varphi(t)<t\;{\rm for~each}\; t>0\;{\rm and~}\lim\limits_{r\to t^+}\varphi(r)<t \;{\rm
for~each~t>0}\Big\}.$$ The following family of control functions is
indicated in  Boyd and Wong \cite{B2} but was later used in
Jotic \cite{B4}.\\
$$\Omega=\Big\{\varphi:[0,\infty)\to [0,\infty):\varphi(t)<t\;{\rm for~each}\; t>0\;{\rm and~}\limsup\limits_{r\to t^+}\varphi(r)<t \;{\rm for~each~t>0}\Big\}.$$
Recently, Alam $et\;al.$ \cite{PGF13} studied the following
relation among above classes of control functions.\\

\begin{prop} \cite{PGF13}. The class $\Omega$
enlarges the classes $\Psi,~\Theta,~\Im$ and $\Phi$ under the
following inclusion relation:
$$\Im\subset\Theta\subset\Psi\subset\Omega~{\rm and}~\Im\subset\Theta\subset\Phi\subset\Omega.$$
\end{prop}

\noindent\textbf{Definition 21}. Let $X$ be a nonempty set and $f$
and $g$ two self-mappings on $X$. Then
 an element $x\in X$ is called a
coincidence point of $f$ and $g$ if
$$f(x)=g(x)=\overline{x},$$
for some $\overline{x}\in X$. Moreover, $\overline{x}$ is called a
point of coincidence of $f$ and $g$. Furthermore, if
$\overline{x}=x$, then $x$ is called a
common fixed point of $f$ and $g$.\\

The following coincidence theorems are crucial results to prove our main results.\\

\noindent\textbf{Lemma 1.} Let $(X,d,\preceq)$ be an ordered metric
space and $E$ an O-complete subspace of $X$. Let $f$ and $g$ be two
self-mappings on $X$. Suppose that the following conditions hold:
\begin{enumerate}
\item [{(i)}] $f(X)\subseteq g(X)\cap E$,
\item [{(ii)}] $f$ is $g$-increasing,
\item [{(iii)}] $f$ and $g$ are O-compatible,
\item [{(iv)}] $g$ is O-continuous,
\item [{(v)}] either $f$ is O-continuous or $(E,d,\preceq)$ has {\it g-MCB} property,
\item [{(vi)}] there exists $x_{0}\in X$ such that
$g(x_{0})\prec\succ f(x_{0})$,
\item [{(vii)}] there exists $\varphi\in \Omega$ such that
$$d(fx,fy)\leq\varphi(d(gx,gy))\;\;\forall~x,y\in X~{\rm with}~
g(x)\prec\succ g(y).$$
 \end{enumerate}
Then $f$ and $g$ have a coincidence point. Moreover, if the
following condition is also hold:\\
\indent\hspace{5mm}(viii) for each pair $x,y\in X$, $\exists~z\in X$ such that $g(x)\prec\succ g(z)$ and $g(y)\prec\succ g(z)$,\\
then $f$ and $g$ have a unique point of coincidence, which remains also a unique common fixed point.\\

\noindent\textbf{Lemma 2.} Let $(X,d,\preceq)$ be an ordered metric
space and $E$ an O-complete subspace of $X$. Let $f$ and $g$ be two
self-mappings on $X$. Suppose that the following conditions hold:
\begin{enumerate}
\item [{(i)}] $f(X)\subseteq E\subseteq g(X)$,
\item [{(ii)}] $f$ is $g$-increasing,
\item [{(iii)}] either $f$ is $(g,{\rm O})$-continuous or $f$ and $g$ are continuous or $(E,d,\preceq)$ has {\it MCB} property,
\item [{(iv)}] there exists $x_{0}\in X$ such that
$g(x_{0})\prec\succ f(x_{0})$,
\item [{(v)}] there exists $\varphi\in \Omega$ such that
$$d(fx,fy)\leq\varphi(d(gx,gy))\;\;\forall~x,y\in X~{\rm with}~
g(x)\prec\succ g(y).$$
 \end{enumerate}
Then $f$ and $g$ have a coincidence point. Moreover, if the
following condition is also hold:\\
\indent\hspace{5mm}(vi) for each pair $x,y\in X$, $\exists~z\in X$ such that $g(x)\prec\succ g(z)$ and $g(y)\prec\succ g(z)$,\\
then $f$ and $g$ have a unique point of coincidence.\\

We skip the proofs of above lemmas as they are proved in Alam $et\;al.$ \cite{PGF13,PGF17,PGF18}.

\section{Extended Notions Upto Product Sets}
\label{SC:Extended Notions Upto Product Sets}

With a view to extend the domain of the mapping
$f:X\longrightarrow X$ to $n$-dimensional product set $X^n$, we
introduce the variants of some existing notions namely:
fixed/coincidence points, commutativity, compatibility, continuity,
$g$-continuity etc. for the mapping $F:X^{n} \rightarrow X$. On the
lines of Herstein \cite{ALG}, a binary operation $\ast$ on a set $S$
is a mapping from $S\times S$ to $S$ and a permutation $\pi$ on a
set $S$ is a one-one mapping from a $S$ onto itself. Throughout this
paper, we adopt the following notations.
\begin{enumerate}
\item [{(1)}] In order to understand a binary operation $\ast$ on $I_n$, we denote the image of any element $(i,k)\in I_n\times I_n$ under $\ast$ by $i_k$ rather
than $\ast(i,k)$.
\item [{(2)}] A binary operation $\ast$ on $I_n$ can be identically
represented by an $n\times n$ matrix throughout its ordered image
such that the first and second components run over rows and columns
respectively, $i.e.$,
$$\ast=[m_{ik}]_{n\times n}\; {\rm where}\; m_{ik}=i_k \;{\rm for\;each}\; i,k\in I_n.$$
\item [{(3)}] A permutation $\pi$ on $I_n$ can be identically represented by an $n$-tuple throughout its ordered image,
$i.e.$,
$$\pi=(\pi(1),\pi(2),...,\pi(n)).$$
\item [{(4)}] $\mathfrak{B_{n}}$ denotes the family of all binary operations $\ast$ on
$I_n$, $i.e.,$
$$\mathfrak{B_{n}}=\{\ast:\ast:I_n\times I_n\rightarrow I_n\}.$$
\item [{(5)}] For any fixed $\iota_n$, $\mathcal{U}_{\iota_n}$ denotes the family of all binary operations $\ast$ on $I_n$ satisfying the following conditions:\\
\indent $(a)$ $\ast(A\times A)\subset A$\\
\indent $(b)$ $\ast(A\times B)\subset B$\\
\indent $(c)$ $\ast(B\times A)\subset B$\\
\indent $(d)$ $\ast(B\times B)\subset A.$
\end{enumerate}
\begin{rem} The following facts are straightforward:
\begin{enumerate}
\item [{(i)}] for each $i\in I_n$, $\{i_1,i_2,...,i_n\}\subseteq I_n.$
\item [{(ii)}] $\mathcal{U}_{\iota_n}\subset \mathfrak{B_{n}}.$
\end{enumerate}
\end{rem} \noindent\textbf{Definition 22}. Let $X$ be a nonempty
set, $\ast\in\mathfrak{B_{n}}$ and $F:X^{n} \rightarrow X$ a
mapping. An element $(x_{1},x_{2},...,x_{n})\in X^{n}$  is called an
$n$-tupled fixed point of $F$ w.r.t. $\ast$ (or, in short,
$\ast$-fixed point of $F$) if
$$F(x_{i_1},x_{i_2},...,x_{i_n})=x_i\;\;{\rm for\;each}\;i\in I_n.$$
Selection of $\ast$ for tripled fixed points of Berinde and Borcut
\cite{T1}, Wu and Liu \cite{TQ} and Berzig and Samet \cite{HD1} are
respectively:
$$\left[\begin{matrix}
1 &2 &3\\
2 &1 &2\\
3 &2 &1\\
\end{matrix}\right],\;
\left[\begin{matrix}
1 &2 &3\\
2 &3 &2\\
3 &2 &1\\
\end{matrix}\right]\;{\rm and}\;
\left[\begin{matrix}
1 &2 &3\\
2 &1 &3\\
3 &3 &2\\
\end{matrix}\right].$$
Selection of $\ast$ for quadrupled fixed points of Karapinar and
Luong \cite{Q1}, Wu and Liu \cite{TQ} and Berzig and Samet
\cite{HD1} are respectively:
$$\left[\begin{matrix}
1 &2 &3 &4\\
2 &3 &4 &1\\
3 &4 &1 &2\\
4 &1 &2 &3\\
\end{matrix}\right],\;
\left[\begin{matrix}
1 &4 &3 &2\\
2 &1 &4 &3\\
3 &2 &1 &4\\
4 &3 &2 &1\\
\end{matrix}\right]\;{\rm and}\;
\left[\begin{matrix}
1 &2 &3 &4\\
1 &2 &4 &3\\
3 &4 &2 &1\\
3 &4 &1 &2\\
\end{matrix}\right].$$
\begin{rem} To ensure the existence of $\ast$-fixed point for a mapping satisfying $\iota_n$-mixed monotone property defined on an ordered metric space, the class
$\mathfrak{B_{n}}$ must be restricted to the subclass
$\mathcal{U}_{\iota_n}$ ($i.e.$ necessarily
$\ast\in\mathcal{U}_{\iota_n}$) so that $\iota_n$-mixed monotone
property can work.
\end{rem}
\begin{prop} The notion of $\ast$-fixed point is equivalent to $\Upsilon$-fixed point.\end{prop}
\noindent{\bf Proof.} Let $(x_{1},x_{2},...,x_{n})\in X^{n}$  is a
$\Upsilon$-fixed point of the mapping $F:X^n\rightarrow X$, where
$\Upsilon=(\sigma_1,\sigma_2,...,\sigma_n)$. Define $\ast:I_n\times
I_n\rightarrow I_n$ by
$$i_k=\sigma_i(k)\;\forall~i,k\in I_n,$$
which implies that $(x_{1},x_{2},...,x_{n})\in X^{n}$  is a
$\ast$-fixed point of $F$.\\
Conversely, suppose that $(x_{1},x_{2},...,x_{n})\in X^{n}$  is an
$\ast$-fixed point of the mapping $F$. Let
$\sigma_1$,$\sigma_2$,...,$\sigma_n$ be the row $n$-tuples of the
matrix representation of $\ast$, $i.e.$,
$$\ast=\left[\begin{matrix}
\sigma_1\\
\sigma_2\\
\vdots\\
\sigma_n\\
\end{matrix}\right]$$
so that $\sigma_1$,$\sigma_2$,...,$\sigma_n$ forms $n$ mappings from
$I_n$ into itself and $\sigma_i(k)=i_k\;\forall~i,k\in I_n.$ Denote
$\Upsilon=(\sigma_1,\sigma_2,...,\sigma_n)$, which amounts to say
that $(x_{1},x_{2},...,x_{n})\in X^{n}$  is a $\Upsilon$-fixed point
of $F$.\\
Moreover, in order to hold $\iota_n$-mixed monotone property, the arguments in Remark 1 and Remark 6 are equivalent.\\

\noindent\textbf{Definition 23}. Let $X$ be a nonempty set,
$\ast\in\mathfrak{B_{n}}$ and $F:X^{n} \rightarrow X$ and $g:
X\rightarrow X$ two mappings. An element $(x_{1},x_{2},...,x_{n})\in
X^{n}$  is called an $n$-tupled coincidence point of $F$ and $g$
w.r.t. $\ast$ (or, in short, $\ast$-coincidence point of $F$ and
$g$) if
$$F(x_{i_1},x_{i_2},...,x_{i_n})=g(x_i)\;\;{\rm for\;each}\;i\in I_n.$$
In this case $(gx_1,gx_2,...,gx_n)$ is called point of
$\ast$-coincidence of $F$ and $g$.\\
Notice that if $g$ is an identity mapping on $X$ then Definition
23 reduces to Definition 22.\\

\noindent\textbf{Definition 24}. Let $X$ be a nonempty set,
$\ast\in\mathfrak{B_{n}}$ and $F:X^{n} \rightarrow X$ and $g:
X\rightarrow X$ two mappings. An element $(x_{1},x_{2},...,x_{n})\in
X^{n}$  is called a common $n$-tupled fixed point of $F$ and $g$
w.r.t. $\ast$ (or, in short, common $\ast$-fixed point of $F$ and
$g$) if
$$F(x_{i_1},x_{i_2},...,x_{i_n})=g(x_i)=x_i\;\;{\rm for\;each}\;i\in I_n.$$

In the following lines, we define four special types $n$-tupled
fixed points, which are somewhat natural.\\

\noindent\textbf{Definition 25}. Let $X$ be a nonempty set and
$F:X^{n} \rightarrow X$ a mapping. An element
$(x_{1},x_{2},...,x_{n})\in X^{n}$ is called a forward cyclic
$n$-tupled fixed point of $F$ if
$$F(x_{i},x_{i+1},...,x_n,x_1,...,x_{i-1})=x_i\;{\rm for\;each}\;i\in I_n$$
$i.e.$\\
\indent\hspace{5cm}$F(x_{1},x_{2},...,x_{n})=x_{1},$\\
\indent\hspace{5cm}$F(x_{2},x_{3},...,x_{n},x_{1})=x_{2},$\\
\indent\hspace{5.6cm}$\vdots$\\
\indent\hspace{5cm}$F(x_{n},x_{1},x_{2},...,x_{n-1})=x_{n}.$\\
This was initiated by Samet and Vetro \cite{NFI}. To obtain this we define $\ast$ as\\
\vspace{0.2cm}\\
\indent\hspace{3cm}$i_k={\begin{cases}i+k-1\;\;\;\;\;\;\;\;\;\;1\leq
k\leq{n-i+1}\cr
\hspace{0.0in}i+k-n-1\;\;\;\;{n-i+2}\leq k\leq{n}\cr\end{cases}}$\\
\vspace{0.2cm}\\
$i.e.$\indent\hspace{2.6cm}$\ast=\left[\begin{matrix}
1 &2 &\cdots &{n-1} &n\\
2 &3 &\cdots &{n} &{1}\\
\vdots &\vdots &&\vdots &\vdots \\
n &1 &\cdots &{n-2} &{n-1}\\
\end{matrix}\right]_{n\times n}$\\

\noindent\textbf{Definition 26}. Let $X$ be a nonempty set and
$F:X^{n} \rightarrow X$ a mapping. An element
$(x_{1},x_{2},...,x_{n})\in X^{n}$  is called a backward cyclic
$n$-tupled fixed point of $F$ if
$$F(x_{i},x_{i-1},...,x_{1},x_n,x_{n-1},...,x_{i+1})=x_i\;{\rm for\;each}\;i\in I_n$$
$i.e.$\\
\indent\hspace{5cm}$F(x_{1},x_{n},x_{n-1},...,x_{2})=x_{1},$\\
\indent\hspace{5cm}$F(x_{2},x_{1},x_{n},...,x_{3})=x_{2},$\\
\indent\hspace{5.6cm}$\vdots$\\
\indent\hspace{5cm}$F(x_{n},x_{n-1},x_{n-2},...,x_{1})=x_{n}.$\\
To obtain this we define $\ast$ as\\
\vspace{0.2cm}\\
\indent\hspace{3cm}$i_k={\begin{cases}i-k+1\;\;\;\;\;\;\;\;\;\;1\leq
k\leq{i}\cr
\hspace{0.0in}n+i-k+1\;\;\;\;{i+1}\leq k\leq{n-1}\cr\end{cases}}$\\
\vspace{0.2cm}\\
$i.e.$\indent\hspace{2.6cm}$\ast=\left[\begin{matrix}
1 &n &{n-1} &\cdots &2\\
2 &1 &{n} &\cdots &3\\
\vdots &\vdots &\vdots &&\vdots \\
n &{n-1} &{n-2} &\cdots &1\\
\end{matrix}\right]_{n\times n}$\\

\noindent\textbf{Definition 27}. Let $X$ be a nonempty set and
$F:X^{n} \rightarrow X$ a mapping. An element
$(x_{1},x_{2},...,x_{n})\in X^{n}$  is called a 1-skew cyclic
$n$-tupled fixed point of $F$ if
$$F(x_{i},x_{i-1},...,x_2,x_1,x_2,...,x_{n-i+1})=x_i\;{\rm for\;each}\;i\in I_n.$$
This was introduced by Gordji and Ramezani \cite{NX1}. To obtain this we define $\ast$ as\\
\vspace{0.2cm}\\
\indent\hspace{3cm}$i_k={\begin{cases}i-k+1\;\;\;\;\;\;\;\;\;\;1\leq
k\leq{i}\cr
\hspace{0.0in}k-i+1\;\;\;\;{i+1}\leq k\leq{n}\cr\end{cases}}$\\

\noindent\textbf{Definition 28}. Let $X$ be a nonempty set and
$F:X^{n} \rightarrow X$ a mapping. An element
$(x_{1},x_{2},...,x_{n})\in X^{n}$ is called a $n$-skew cyclic
$n$-tupled fixed point of $F$ if
$$F(x_{i},x_{i+1},...,x_{n-1},x_n,x_{n-1},...,x_{n-i+1})=x_i\;{\rm for\;each}\;i\in I_n.$$
To obtain this we define $\ast$ as\\
\indent\hspace{3cm}$i_k={\begin{cases}i+k-1\;\;\;\;\;\;\;\;\;\;1\leq
k\leq{n-i+1}\cr \hspace{0.0in}2n-i-k+1\;\;\;\;{n-i+2}\leq
k\leq{n}\cr\end{cases}}$

\begin{rem} In particular for $n=4$, forward cyclic and backward cyclic
$n$-tupled fixed points reduce to quadrupled fixed points of
Karapinar and Luong \cite{Q1} and Wu and Liu \cite{TQ} respectively.
Also, for $n=3$, 1-skew cyclic and $n$-skew cyclic
$n$-tupled fixed points reduce to tripled fixed points of Berinde
and Borcut \cite{T1} and Wu and Liu \cite{TQ} respectively.
\end{rem}

\noindent\textbf{Definition 29}. A binary operation $\ast$ on $I_n$
is called permuted if each row of matrix representation of $\ast$
forms a permutation on $I_n.$
\begin{example} On $I_3$, consider two binary operations\\
\indent\hspace{2cm}$\ast$=$\left[\begin{matrix}
1 &2 &3\\
2 &1 &3\\
3 &2 &1\\
\end{matrix}\right]$,\;\;
$\circ$=$\left[\begin{matrix}
1 &2 &3\\
2 &1 &3\\
3 &3 &2\\
\end{matrix}\right]$\\
$\ast$ is permuted as each of rows $(1,2,3),(2,1,3),(3,2,1)$ is a
permutation on $I_3$. While $\circ$ is not permuted as last row
$(3,3,2)$ is not permutation on $I_3$.
\end{example}
It is clear that binary operations defined for forward cyclic
and backward cyclic $n$-tupled fixed points are permuted
while for 1-skew cyclic and $n$-skew cyclic $n$-tupled
fixed points are not permuted.
\begin{prop} A permutation $\ast$ on $I_n$ is permuted iff for each $i\in I_n,$
$$\{i_1,i_2,...,i_n\}=I_n.$$
\end{prop}
\noindent\textbf{Definition 30}. Let $(X,d)$ be a metric space,
$F:X^{n} \rightarrow X$ a mapping and $(x_1,x_2,...,x_n)\in X^n$. We
say that $F$ is continuous at $(x_1,x_2,...,x_n)$ if for any
sequences $\{x_1^{(m)}\},\{x_2^{(m)}\},$ $..., \{x_n^{(m)}\}\subset
X$,
$$x_1^{(m)}\stackrel{d}{\longrightarrow} x_1,\;x_2^{(m)}\stackrel{d}{\longrightarrow} x_2,...,\;x_n^{(m)}\stackrel{d}{\longrightarrow} x_n$$
$$\Longrightarrow F(x_1^{(m)},x_2^{(m)},...,x_n^{(m)})\stackrel{d}{\longrightarrow}F(x_1,x_2,...,x_n).$$
Moreover, $F$ is called continuous if it is continuous at each point of $X^n$.\\

\noindent\textbf{Definition 31}. Let $(X,d)$ be a metric space and
$F:X^{n} \rightarrow X$ and $g: X\rightarrow X$ two mappings and
$(x_1,x_2,...,x_n)\in X^n$. We say that $F$ is $g$-continuous at
$(x_1,x_2,...,x_n)$ if for any sequences
$\{x_1^{(m)}\},\{x_2^{(m)}\},...,\{x_n^{(m)}\}\subset X$,
$$g(x_1^{(m)})\stackrel{d}{\longrightarrow} g(x_1),\;g(x_2^{(m)})\stackrel{d}{\longrightarrow} g(x_2),...,\;g(x_n^{(m)})\stackrel{d}{\longrightarrow} g(x_n)$$
$$\Longrightarrow F(x_1^{(m)},x_2^{(m)},...,x_n^{(m)})\stackrel{d}{\longrightarrow}F(x_1,x_2,...,x_n).$$
Moreover, $F$ is called $g$-continuous if it is $g$-continuous at each point of $X^n$.\\

Notice that setting $g=I$ (identity mapping on $X$), Definition 31
reduces to Definition 30.\\

\noindent\textbf{Definition 32}. Let $(X,d,\preceq)$ be an ordered
metric space, $F:X^{n} \rightarrow X$ a mapping and
$(x_1,x_2,...,x_n)\in X^n$. We say that $F$ is O-continuous at
$(x_1,x_2,...,x_n)\in X^n$ if for any sequences
$\{x_1^{(m)}\},\{x_2^{(m)}\},...,\{x_n^{(m)}\}\subset X$,
$$x_1^{(m)}\uparrow\downarrow x_1,\;x_2^{(m)}\uparrow\downarrow x_2,...,\;x_n^{(m)}\uparrow\downarrow x_n$$
$$\Longrightarrow F(x_1^{(m)},x_2^{(m)},...,x_n^{(m)})\stackrel{d}{\longrightarrow}F(x_1,x_2,...,x_n).$$
Moreover, $F$ is called O-continuous if it is O-continuous at each point of $X^n$.\\

\noindent\textbf{Definition 33}. Let $(X,d,\preceq)$ be an ordered
metric space, $F:X^{n} \rightarrow X$ and $g: X\rightarrow X$ two
mappings and $(x_1,x_2,...,x_n)\in X^n$. We say that $F$ is $(g,{\rm
O})$-continuous at $(x_1,x_2,...,x_n)$ if for any sequences
$\{x_1^{(m)}\},\{x_2^{(m)}\},...,\{x_n^{(m)}\}\subset X$,
$$g(x_1^{(m)})\uparrow\downarrow g(x_1),\;g(x_2^{(m)})\uparrow\downarrow g(x_2),...,\;g(x_n^{(m)})\uparrow\downarrow g(x_n)$$
$$\Longrightarrow F(x_1^{(m)},x_2^{(m)},...,x_n^{(m)})\stackrel{d}{\longrightarrow}F(x_1,x_2,...,x_n).$$
Moreover, $F$ is called $(g,{\rm
O})$-continuous if it is $(g,{\rm
O})$-continuous at each point of $X^n$.\\

Notice that setting $g=I$ (identity mapping on $X$), Definition 33
reduces to Definition 32.
\begin{rem} Let $(X,d,\preceq)$ be an ordered
metric space and $g: X\rightarrow X$ a mapping. If $F:X^{n}
\rightarrow X$ is a continuous (resp. $g$-continuous) mapping then
$F$ is also O-continuous (resp. $(g,{\rm O})$-continuous).
\end{rem}

\noindent\textbf{Definition 34}. Let $X$ be a nonempty set and
$F:X^{n} \rightarrow X$ and $g: X\rightarrow X$ two mappings. We say
that $F$ and $g$ are commuting if for all $x_1,x_2,...,x_n\in X,$
$$g(F(x_1,x_2,...,x_n))=F(gx_1,gx_2,...,gx_n).$$

\noindent\textbf{Definition 35}. Let $(X,d)$ be a metric space and
$F:X^{n} \rightarrow X$ and $g: X\rightarrow X$ two mappings. We say
that $F$ and $g$ are $\ast$-compatible if for any sequences
$\{x_1^{(m)}\},\{x_2^{(m)}\},...,\{x_n^{(m)}\}\subset X$ and for any
$z_1,z_2,...,z_n\in X$,
$$g(x_i^{(m)})\stackrel{d}{\longrightarrow} z_i\;{\rm and}\;F(x_{i_1}^{(m)},x_{i_2}^{(m)},...,x_{i_n}^{(m)})\stackrel{d}{\longrightarrow} z_i\;\;{\rm for~each}~i\in I_n$$
$$\Longrightarrow\lim\limits_{m\to \infty}d(gF(x_{i_1}^{(m)},x_{i_2}^{(m)},...,x_{i_n}^{(m)}),F(gx_{i_1}^{(m)},gx_{i_2}^{(m)},...,gx_{i_n}^{(m)}))=0\;\;{\rm for~each}~i\in I_n.$$

\noindent\textbf{Definition 36}. Let $(X,d,\preceq)$ be an ordered
metric space and $F:X^{n} \rightarrow X$ and $g: X\rightarrow X$ two
mappings. We say that $F$ and $g$ are $(\ast,{\rm O})$-compatible if
for any sequences
$\{x_1^{(m)}\},\{x_2^{(m)}\},...,\{x_n^{(m)}\}\subset X$ and for any
$z_1,z_2,...,z_n\in X$,
$$g(x_i^{(m)})\uparrow\downarrow z_i\;{\rm and}\;F(x_{i_1}^{(m)},x_{i_2}^{(m)},...,x_{i_n}^{(m)})\uparrow\downarrow z_i\;\;{\rm for~each}~i\in I_n$$
$$\Longrightarrow\lim\limits_{m\to \infty}d(gF(x_{i_1}^{(m)},x_{i_2}^{(m)},...,x_{i_n}^{(m)}),F(gx_{i_1}^{(m)},gx_{i_2}^{(m)},...,gx_{i_n}^{(m)}))=0\;\;{\rm for~each}~i\in I_n.$$

\noindent\textbf{Definition 37}. Let $X$ be a nonempty set and
$F:X^{n} \rightarrow X$ and $g: X\rightarrow X$ two mappings. We say
that $F$ and $g$ are weakly $\ast$-compatible if for any
$x_1,x_2,...,x_n\in X,$
$$g(x_i)=F(x_{i_1},x_{i_2},...,x_{i_n})\;\;{\rm for~each}~i\in I_n$$
$$\Longrightarrow g(F(x_{i_1},x_{i_2},...,x_{i_n}))=F(gx_{i_1},gx_{i_2},...,gx_{i_n})\;\;{\rm for~each}~i\in I_n.$$
\begin{rem} Evidently, in an ordered metric space, commutativity $\Rightarrow$
$\ast$-compatibility $\Rightarrow$ $(\ast,{\rm O})$-compatibility
$\Rightarrow$ weak $\ast$-compatibility. \end{rem}
\begin{prop} If $F:X^{n} \rightarrow X$ and $g: X\rightarrow X$ are
weakly $\ast$-compatible, then every point of $\ast$-coincidence of
$F$ and $g$ is also an $\ast$-coincidence point of $F$
 and $g$.\end{prop}
\noindent{\bf Proof.} Let
$(\overline{x}_1,\overline{x}_2,...,\overline{x}_n)\in X^n$ be a
point of $\ast$-coincidence of $F$ and $g$, then $\exists~
x_1,x_2,...,x_n\in X$ such that
$F(x_{i_1},x_{i_2},...,x_{i_n})=g(x_i)=\overline{x}_i$ for each
$i\in I_n$. Now, we have to show that
$(\overline{x}_1,\overline{x}_2,...,\overline{x}_n)$ is a
$\ast$-coincidence point of $F$ and $g.$ On using weak
$\ast$-compatibility of $F$ and $g$, for each $i\in I_n$, we have
\begin{eqnarray*}
g(\overline{x}_i)&=& g(F(x_{i_1},x_{i_2},...,x_{i_n}))\\
&=& F(gx_{i_1},gx_{i_2},...,gx_{i_n})\\
&=&F(\overline{x}_{i_1},\overline{x}_{i_2},...,\overline{x}_{i_n}),
\end{eqnarray*}
which implies that
$(\overline{x}_1,\overline{x}_2,...,\overline{x}_n)$ is an
$\ast$-coincidence point of $F$ and $g$.

\section{Auxiliary Results}
\label{SC:Auxiliary Results}
The classical technique involved in the proofs of the
multi-tupled fixed point results due to Bhasker and
Lakshmikantham \cite{C1}, Berinde and Borcut \cite{T1}, Karapinar
and Luong \cite{Q1}, Imdad $et\;al.$ \cite{n1}, Berzig and Samet
\cite{HD1}, Rold$\acute{\rm a}$n $et\;al.$ \cite{MD1} $etc.$ is
very long specially due to the involvement of $n$ coordinates of the elements and the
sequences in $X^n$. In 2011, Berinde \cite{C6}, generalized the
coupled fixed point results of Bhasker and Lakshmikantham \cite{C1}
by using the corresponding fixed point theorems on ordered metric
spaces. Recently, utilizing this technique several authors such as:
Jleli $et\;al.$ \cite{PGF10}, Samet $et\;al.$ \cite{R1}, Wu and Liu
\cite{PGF12}, Wu and Liu \cite{TQ}, Dalal $et\;al.$ \cite{R2},
Radenovi$\acute{\rm c}$ \cite{R3}, Al-Mezel $et\;al.$ \cite{MD3},
Rold$\acute{\rm a}$n $et\;al.$ \cite{MD2}, Rad $et\;al.$ \cite{NX5},
Sharma $et\;al.$ \cite{R4} $etc.$ proved some multi-tupled fixed
point results. The technique of reduction of multi-tupled fixed
point results from corresponding fixed point results is fascinating,
relatively simpler, shorter and more effective than classical
technique. Due to this fact, we also prove our results using later
technique. In this section, we discuss some basic results, which
provide the tools for reduction of the multi-tupled fixed point
results from the corresponding fixed point results. Before doing
so, we consider the following induced notations.
\begin{enumerate}
\item [{(1)}] For any U$=(x_{1},x_{2},...,x_{n})\in X^n$, for an $\ast\in \mathfrak{B_{n}}$ and for each $i\in I_n$, U$^{\ast}_i$ denotes the ordered element
 $(x_{i_1},x_{i_2},...,x_{i_n})$ of $X^n$.
\item [{(2)}] For each $\ast\in \mathfrak{B_{n}}$, a mapping $F:X^n\rightarrow X$ induce an associated mapping $F_\ast:X^n\rightarrow X^n$ defined by
$$F_\ast({\rm U})=(F{\rm U}^{\ast}_1,F{\rm U}^{\ast}_2,...,F{\rm U}^{\ast}_n)\;\;\forall~{\rm U}\in X^n.$$
\item [{(3)}] A mapping $g:X\rightarrow X$ induces an associated mapping $G:X^n\rightarrow
X^n$ defined by
$$G({\rm U})=(gx_1,gx_2,...,gx_n)\;\;\forall~{\rm U}=(x_{1},x_{2},...,x_{n})\in X^n. $$
\item [{(4)}] For a metric space $(X,d)$, $\Delta_n$ and $\nabla_n$ denote
two metrics on product set $X^n$ defined by:\\ for all
U=$(x_{1},x_{2},...,x_{n})$, V=$(y_{1},y_{2},...,y_{n})\in X^n,$
$$\Delta_n({\rm U,V})=\frac{1}{n}\sum\limits_{i=1}^{n}d(x_i,y_i)$$
$$\nabla_n({\rm U,V})=\max\limits_{i\in I_n}d(x_i,y_i).$$
\item [{(5)}] For any ordered set $(X,\preceq)$ and a fixed $\iota_n$, $\sqsubseteq_{\iota_n}$ denotes a partial order on $X^n$ defined
by:\\
for all U=$(x_{1},x_{2},...,x_{n})$, V=$(y_{1},y_{2},...,y_{n})\in
X^n,$
$${\rm U}\sqsubseteq_{\iota_n}{\rm V}\Leftrightarrow x_i\preceq y_i \;{\rm {for\; each} }\; i\in A\; {\rm {and} }\;x_i\succeq y_i\; {\rm {for\; each} }\; i\in B.$$
\end{enumerate}
\begin{rem} The following facts are straightforward:
\begin{enumerate}
\item [{(i)}] $F_\ast(X^n)\subseteq (FX^n)^n.$
\item [{(ii)}] $G(X^n)=(gX)^n.$
\item [{(iii)}] $(GU)^\ast_i=G(U^\ast_i)\;\forall~{\rm U}\in X^n.$
\item [{(iv)}] $\frac{1}{n}\nabla_n\leq\Delta_n\leq\nabla_n$ ($i.e.$ both the
metrics $\Delta_n$ and $\nabla_n$ are equivalent).
\end{enumerate}
\end{rem}
\noindent{\bf Lemma 3.} Let $X$ be a nonempty set, $E\subseteq X$,
$F:X^{n} \rightarrow X$ and $g: X\rightarrow X$ two mappings and
$\ast\in\mathfrak{B_{n}}$.
\begin{enumerate}
\item [{(i)}] If $F(X^n)\subseteq g(X)\cap E$ then $F_\ast(X^n)\subseteq (FX^n)^n\subseteq
G(X^n)\cap E^n$.
\item [{(ii)}] If $F(X^n)\subseteq E\subseteq g(X)$ then $F_\ast(X^n)\subseteq
(FX^n)^n\subseteq E^n \subseteq G(X^n)$.
\item [{(iii)}] An element $(x_{1},x_{2},...,x_{n})\in
X^n$ is $\ast$-coincidence point of $F$ and $g$ iff
$(x_{1},x_{2},...,x_{n})$ is a coincidence point of $F_{\ast}$ and $G$.
\item [{(iv)}] An element $(\overline{x}_{1},\overline{x}_{2},...,\overline{x}_{n})\in
X^n$ is point of $\ast$-coincidence of $F$ and $g$ iff
$(\overline{x}_{1},\overline{x}_{2},...,\overline{x}_{n})$ is a
point of coincidence of $F_{\ast}$ and $G$.
\item [{(v)}] An element $(x_{1},x_{2},...,x_{n})\in X^n$ is common $\ast$-fixed point of $F$ and $g$ iff
$(x_{1},x_{2},$$...,$$x_{n})$ is a common fixed point of $F_{\ast}$ and $G$.
\end{enumerate}
\noindent{\bf Proof.} The proof of the lemma is straightforward and hence it
is left to the reader.\\

\noindent{\bf Lemma 4.} Let $(X,\preceq)$ be an ordered set, $g:
X\rightarrow X$ a mapping and $\ast\in\mathcal{U}_{\iota_n}$. If
$G({\rm U})\sqsubseteq_{\iota_n}G({\rm V})$ for some U,V$\in X^n$
then
\begin{enumerate}
\item [{(i)}] $G({\rm U}^\ast_i)\sqsubseteq_{\iota_n}G({\rm V}^\ast_i)$\;\;for each $i\in A$,
\item [{(ii)}] $G({\rm U}^\ast_i)\sqsupseteq_{\iota_n}G({\rm V}^\ast_i)$\;\;for each $i\in B$.
\end{enumerate}
\noindent{\bf Proof.} Let U=$(x_{1},x_{2},...,x_{n})$ and
V=$(y_{1},y_{2},...,y_{n})$, then we have
$$(gx_{1},gx_{2},...,gx_{n})\sqsubseteq_{\iota_n}(gy_{1},gy_{2},...,gy_{n}),$$
which implies that
$$g(x_i)\preceq g(y_i) \;{\rm {for\; each} }\; i\in A\; {\rm {and} }\;g(x_i)\succeq g(y_i)\; {\rm {for\; each} }\; i\in B.\eqno(1)$$
Now, we consider the following
cases:\\
{\bf Case I:} $i\in A$. Then by the definition of
$\mathcal{U}_{\iota_n}$, we have
$$i_k\in A\;{\rm {for\; each} }\; k\in A\; {\rm {and} }\;i_k\in B\; {\rm {for\; each} }\; k\in B.\eqno(2)$$
Using (1) and (2), we obtain
$$g(x_{i_k})\preceq g(y_{i_k})\;{\rm {for\; each} }\; k\in A\; {\rm {and} }\;g(x_{i_k})\succeq g(y_{i_k})\; {\rm {for\; each} }\; k\in B,$$
which implies that
$$(gx_{i_1},gx_{i_2},...,gx_{i_n})\sqsubseteq_{\iota_n}(gy_{i_1},gy_{i_2},...,gy_{i_n}),$$
$i.e.$ $$G({\rm U}^\ast_i)\sqsubseteq_{\iota_n}G({\rm
V}^\ast_i)\;\;{\rm for~
each~} i\in A.$$ Hence, (i) is proved.\\
{\bf Case II:} $i\in B$. Then by the definition of
$\mathcal{U}_{\iota_n}$, we have
$$i_k\in B\;{\rm {for\; each} }\; k\in A\; {\rm {and} }\;i_k\in A\; {\rm {for\; each} }\; k\in B.\eqno(3)$$
Using (1) and (3), we obtain
$$g(x_{i_k})\succeq g(y_{i_k})\;{\rm {for\; each} }\; k\in A\; {\rm {and} }\;g(x_{i_k})\preceq g(y_{i_k})\; {\rm {for\; each} }\; k\in B,$$
which implies that
$$(gx_{i_1},gx_{i_2},...,gx_{i_n})\sqsupseteq_{\iota_n}(gy_{i_1},gy_{i_2},...,gy_{i_n}),$$
$i.e.$ $$G({\rm U}^\ast_i)\sqsupseteq_{\iota_n}G({\rm
V}^\ast_i)\;\;{\rm for~
each~} i\in B.$$ Hence, (ii) is proved.\\

\noindent{\bf Lemma 5.} Let $(X,\preceq)$ be an ordered set,
$F:X^{2} \rightarrow X$ and $g: X\rightarrow X$ two mappings and
$\ast\in\mathcal{U}_{\iota_n}$. If $F$ has $\iota_n$-mixed
$g$-monotone property then $F_{\ast}$ is $G$-increasing in ordered set $(X^n,\sqsubseteq_{\iota_n})$.\\
\noindent{\bf Proof.} Take U=$(x_{1},x_{2},...,x_{n})$,
V=$(y_{1},y_{2},...,y_{n})\in X^n$ with $G({\rm
U})\sqsubseteq_{\iota_n}G({\rm V})$. Consider the following
cases:\\
{\bf Case I:} $i\in A$. Owing to Lemma 3, we obtain
$$G({\rm U}^\ast_i)\sqsubseteq_{\iota_n}G({\rm V}^\ast_i),$$ which implies
that
$$g(x_{i_k})\preceq g(y_{i_k})\;{\rm {for\; each} }\; k\in A\; {\rm {and} }\;g(x_{i_k})\succeq g(y_{i_k})\; {\rm {for\; each} }\; k\in B.\eqno(4)$$
On using (4) and $\iota_n$-mixed $g$-monotone property of $F$, we
obtain
\begin{eqnarray*}
F({\rm U}^\ast_i)&=& F(x_{i_1},x_{i_2},...,x_{i_n})\\
&\preceq& F(y_{i_1},x_{i_2},...,x_{i_n})\\
&\preceq& F(y_{i_1},y_{i_2},...,x_{i_n})\\
&\preceq& \cdots\\
&\preceq& F(y_{i_1},y_{i_2},...,y_{i_n})\\
&=& F({\rm V}^\ast_i)
\end{eqnarray*}
so that
$$F({\rm U}^\ast_i)\preceq F({\rm V}^\ast_i)\;{\rm for~each}~i\in A.\eqno(5)$$
{\bf Case II:} $i\in B$. Owing to Lemma 3, we obtain
$$G({\rm U}^\ast_i)\sqsupseteq_{\iota_n}G({\rm V}^\ast_i),$$ which implies
that
$$g(x_{i_k})\succeq g(y_{i_k})\;{\rm {for\; each} }\; k\in A\; {\rm {and} }\;g(x_{i_k})\preceq g(y_{i_k})\; {\rm {for\; each} }\; k\in B.\eqno(6)$$
On using (6) and $\iota_n$-mixed $g$-monotone property of $F$, we
obtain
\begin{eqnarray*}
F({\rm U}^\ast_i)&=& F(x_{i_1},x_{i_2},...,x_{i_n})\\
&\succeq& F(y_{i_1},x_{i_2},...,x_{i_n})\\
&\succeq& F(y_{i_1},y_{i_2},...,x_{i_n})\\
&\succeq& \cdots\\
&\succeq& F(y_{i_1},y_{i_2},...,y_{i_n})\\
&=& F({\rm V}^\ast_i)
\end{eqnarray*}
so that
$$F({\rm U}^\ast_i)\succeq F({\rm V}^\ast_i)\;{\rm for~each}~i\in B.\eqno(7)$$
From (5) and (7), we get
\begin{eqnarray*}
F_\ast({\rm U})&=&(F{\rm U}^{\ast}_1,F{\rm U}^{\ast}_2,...,F{\rm U}^{\ast}_n)\\
&\sqsubseteq_{\iota_n}& (F{\rm V}^{\ast}_1,F{\rm V}^{\ast}_2,...,F{\rm V}^{\ast}_n)\\
&=& F_\ast({\rm V}).
\end{eqnarray*}
Hence, $F_\ast$ is $G$-increasing.\\

\noindent{\bf Lemma 6.}  Let $(X,d)$ be a metric space, $g:
X\rightarrow X$ a mapping and $\ast\in\mathfrak{B_{n}}$. Then, for
any U=$(x_{1},x_{2},...,x_{n})$,V=$(y_{1},y_{2},...,y_{n})\in X^n$
and for each $i\in I_n$,
\begin{enumerate}
\item [{(i)}] $\frac{1}{n}\sum\limits_{k=1}^{n}d(gx_{i_k},gy_{i_k})=\frac{1}{n}\sum\limits_{j=1}^{n}d(gx_j,gy_j)=\Delta_n(G{\rm U},G{\rm V})$\; provided $\ast$
is permuted,\\
\item [{(ii)}] $\max\limits_{k\in
I_n}d(gx_{i_k},gy_{i_k})=\max\limits_{j\in
I_n}d(gx_j,gy_j)=\nabla_n(G{\rm U},G{\rm V})$\;provided $\ast$
is permuted,\\
\item [{(iii)}] $\max\limits_{k\in
I_n}d(gx_{i_k},gy_{i_k})\leq\max\limits_{j\in
I_n}d(gx_j,gy_j)=\nabla_n(G{\rm U},G{\rm V})$.
\end{enumerate}
\noindent{\bf Proof.} The result is followed by using Remark 5 (item
(i)) and Proposition 4.
\begin{prop} Let $(X,d)$ be a metric space. Then for any sequence ${\rm U}^{(m)}\subset X^n$
and any ${\rm U}\in X^n$, where
${\rm U}^{(m)}=(x^{(m)}_1,x^{(m)}_2,...,x^{(m)}_n)$ and
${\rm U}=(x_1,x_2,...,x_n)$
\begin{enumerate}
\item [{(i)}] ${\rm U}^{(m)}\stackrel{\Delta_n}{\longrightarrow} {\rm U}\Leftrightarrow x_i^{(m)}\stackrel{d}{\longrightarrow} x_i\;{\rm for~each}~i\in I_n.$
\item [{(ii)}] ${\rm U}^{(m)}\stackrel{\nabla_n}{\longrightarrow} {\rm U}\Leftrightarrow x_i^{(m)}\stackrel{d}{\longrightarrow} x_i\;{\rm for~each}~i\in I_n.$
\end{enumerate}
\end{prop}

\noindent{\bf Lemma 7.}  Let $(X,d)$ be a metric space, $F:X^{n}
\rightarrow X$ and $g:X\rightarrow X$ two mappings and
$\ast\in\mathfrak{B_{n}}$.
\begin{enumerate}
\item [{(i)}] If $g$ is continuous then $G$ is continuous in both
metric spaces $(X^n,\Delta_n)$ and $(X^n,\nabla_n)$,
\item [{(ii)}] If $F$ is continuous then $F_{\ast}$ is continuous in
both metric spaces $(X^n,\Delta_n)$ and $(X^n,\nabla_n)$.
\end{enumerate}
\noindent{\bf Proof.} (i) Take a sequence ${\rm U}^{(m)}\subset X^n$
and a ${\rm U}\in X^n$, where
U$^{(m)}=(x^{(m)}_1,x^{(m)}_2,...,x^{(m)}_n)$ and
U$=(x_1,x_2,...,x_n)$ such that
$${\rm U}^{(m)}\stackrel{\Delta_n}{\longrightarrow} {\rm U}\;\;\;({\rm resp.~U}^{(m)}\stackrel{\nabla_n}{\longrightarrow} {\rm U}),$$
which, on using Proposition 6 implies that
$$x_i^{(m)}\stackrel{d}{\longrightarrow} x_i\;{\rm for~each}~i\in
I_n.\eqno(8)$$ Using (8) and continuity of $g$, we get
$$g(x_i^{(m)})\stackrel{d}{\longrightarrow} g(x_i)\;{\rm for~each}~i\in
I_n,$$ which, again by using Proposition 6 gives rise
$$G({\rm U}^{(m)})\stackrel{\Delta_n}{\longrightarrow} G({\rm U})\;\;\;({\rm resp.}~G({\rm U}^{(m)})\stackrel{\nabla_n}{\longrightarrow} G({\rm U})).$$
Hence, $G$ is continuous in metric space $(X^n,\Delta_n)$
(resp. $(X^n,\nabla_n)$)\\
(ii) Take a sequence ${\rm U}^{(m)}\subset X^n$ and a ${\rm U}\in
X^n$, where U$^{(m)}=(x^{(m)}_1,x^{(m)}_2,...,x^{(m)}_n)$ and
U$=(x_1,x_2,...,x_n)$ such that
$${\rm U}^{(m)}\stackrel{\Delta_n}{\longrightarrow} {\rm U}\;\;\;({\rm resp.~U}^{(m)}\stackrel{\nabla_n}{\longrightarrow} {\rm U}),$$
which, on using Proposition 6 implies that
$$x_i^{(m)}\stackrel{d}{\longrightarrow} x_i\;{\rm for~each}~i\in
I_n.$$ It follows for each $i\in I_n$ that
$$x_{i_1}^{(m)}\stackrel{d}{\longrightarrow} x_{i_1},x_{i_2}^{(m)}\stackrel{d}{\longrightarrow} x_{i_2},..., x_{i_n}^{(m)}\stackrel{d}{\longrightarrow} x_{i_n}.
\eqno(9)$$ Using (9) and continuity of $F$, we get
$$F(x_{i_1}^{(m)},x_{i_2}^{(m)},...,x_{i_n}^{(m)})\stackrel{d}{\longrightarrow}F(x_{i_1},x_{i_2},...,x_{i_n})$$
so that $$F({\rm U}^{(m)\ast}_i)\stackrel{d}{\longrightarrow}F(\rm
U)\;{\rm for~each}~i\in I_n.$$
 which, again by using Proposition 6 gives rise
$$F_\ast({\rm U}^{(m)})\stackrel{\Delta_n}{\longrightarrow}F_\ast({\rm U})\;\;\;({\rm resp.}~F_\ast({\rm U}^{(m)})\stackrel{\nabla_n}{\longrightarrow} F_\ast({\rm U})).$$
Hence, $F_\ast$ is continuous in metric space $(X^n,\Delta_n)$
(resp. $(X^n,\nabla_n)$)\\
\begin{prop} Let $(X,d,\preceq)$ be an ordered metric space and $\{{\rm
U}^{(m)}\}$ a sequence in $X^n,$ where ${\rm
U}^{(m)}=(x_1^{(m)},x_2^{(m)},...,x_n^{(m)})$.
\begin{enumerate}
\item [{(i)}] If $\{{\rm
U}^{(m)}\}$ is monotone in $(X^n,\sqsubseteq_{\iota_n})$ then each
of $\{x_1^{(m)}\}$,$\{x_2^{(m)}\}$,...,$\{x_n^{(m)}\}$ is monotone
in $(X,\preceq)$.
\item [{(ii)}] If $\{{\rm
U}^{(m)}\}$ is Cauchy in $(X^n,\Delta_n)$ (similarly in
$(X^n,\nabla_n)$) then each of $\{x_1^{(m)}\}$,$\{x_2^{(m)}\}$,
...,$\{x_n^{(m)}\}$ is Cauchy in $(X,d)$.
\end{enumerate}
\end{prop}
{\noindent\bf{Lemma 8.}} Let $(X,d,\preceq)$ be an ordered metric
space, $E\subseteq X$ and $\ast\in\mathfrak{B_{n}}$. Let $F:X^{n}
\rightarrow X$ and $g:X\rightarrow X$ be two mappings.
\begin{enumerate}
\item [{(i)}] If $(E,d,\preceq)$ is O-complete then
$(E^n,\Delta_n,\sqsubseteq_{\iota_n})$ and
$(E^n,\nabla_n,\sqsubseteq_{\iota_n})$ both are O-complete.
\item [{(ii)}] If $F$ and $g$ are $(\ast,{\rm
O})$-compatible then $F_{\ast}$ and $G$ are O-compatible in both
ordered metric spaces $(X^n,\Delta_n,\sqsubseteq_{\iota_n})$ and
$(X^n,\nabla_n,\sqsubseteq_{\iota_n})$,
\item [{(iii)}] If $g$ is O-continuous then $G$ is O-continuous in both ordered metric spaces $(X^n,\Delta_n,\sqsubseteq_{\iota_n})$
and $(X^n,\nabla_n,\sqsubseteq_{\iota_n})$,
\item [{(iv)}] If $F$ is O-continuous then
$F_{\ast}$ is O-continuous in both ordered metric spaces
$(X^n,\Delta_n,\sqsubseteq_{\iota_n})$ and
$(X^n,\nabla_n,\sqsubseteq_{\iota_n})$,
\item [{(v)}] If $F$ is $(g,{\rm O})$-continuous then
$F_{\ast}$ is $(G,{\rm O})$-continuous in both ordered metric spaces
$(X^n,\Delta_n,\sqsubseteq_{\iota_n})$ and
$(X^n,\nabla_n,\sqsubseteq_{\iota_n})$,
\item [{(vi)}] If $(E,d,\preceq)$ has {\it g-MCB} property then both $(E^n,\Delta_n,\sqsubseteq_{\iota_n})$ and $(E^n,\nabla_n,\sqsubseteq_{\iota_n})$ have
{\it G-MCB} property,
\item [{(vii)}] If $(E,d,\preceq)$ has {\it MCB} property then both $(E^n,\Delta_n,\sqsubseteq_{\iota_n})$ and $(E^n,\nabla_n,\sqsubseteq_{\iota_n})$ have {\it MCB}
property.
\end{enumerate}
\noindent{\bf Proof.} (i) Let $\{{\rm U}^{(m)}\}$ be a monotone
Cauchy sequence in $(E^n,\Delta_n,\sqsubseteq_{\iota_n})$ (resp. in
$(E^n,\nabla_n,\sqsubseteq_{\iota_n})$). Denote
U$^{(m)}=(x^{(m)}_1,x^{(m)}_2,...,x^{(m)}_n)$, then by Proposition
7, each of $\{x_1^{(m)}\}$,$\{x_2^{(m)}\}$,...,$\{x_n^{(m)}\}$ is a
monotone Cauchy sequence in $(E,d,\preceq)$. By O-completeness of
$(E,d,\preceq)$, $\exists~x_1,x_2,...,x_n\in E$ such that
$$x_i^{(m)}\stackrel{d}{\longrightarrow} x_i\;{\rm for~each}~i\in
I_n,$$ which using Proposition 6, implies that $${\rm
U}^{(m)}\stackrel{\Delta_n}{\longrightarrow} {\rm U}\;\;({\rm
resp.}~ {\rm U}^{(m)}\stackrel{\nabla_n}{\longrightarrow} {\rm
U}),$$ where ${\rm U}=(x_1,x_2,...,x_n)$. It follows that
$(E^n,\Delta_n,\sqsubseteq_{\iota_n})$ (resp.
$(E^n,\nabla_n,\sqsubseteq_{\iota_n})$) is O-complete.\\

(ii) Take a sequence $\{{\rm U}^{(m)}\}\subset X^n$ such that
$\{G{\rm U}^{(m)}\}$ and $\{F_\ast{\rm U}^{(m)}\}$ are monotone
(w.r.t. partial order $\sqsubseteq_{\iota_n}$) and
$$G({\rm U}^{(m)})\stackrel{\Delta_n}{\longrightarrow}{\rm W}\;{\rm and}\;F_\ast({\rm U}^{(m)})\stackrel{\Delta_n}{\longrightarrow}{\rm W}$$ for some W$\in X^n$. Write
U$^{(m)}=(x^{(m)}_1,x^{(m)}_2,...,x^{(m)}_n)$ and
W$=(z_1,z_2,...,z_n)$. Then, by using Propositions 6 and 7, we
obtain
$$g(x_i^{(m)})\uparrow\downarrow z_i\;{\rm and}\;F(x_{i_1}^{(m)},x_{i_2}^{(m)},...,x_{i_n}^{(m)})\uparrow\downarrow z_i\;\;{\rm for~each}~i\in I_n.\eqno(10)$$
On using (10) and $(\ast,{\rm O})$-compatibility of $F$ and $g$, we
have
$$\lim\limits_{m\to \infty}d(gF(x_{i_1}^{(m)},x_{i_2}^{(m)},...,x_{i_n}^{(m)}),F(gx_{i_1}^{(m)},gx_{i_2}^{(m)},...,gx_{i_n}^{(m)}))=0\;\;{\rm for~each}~i\in I_n$$
$i.e.$
$$\lim\limits_{m\to \infty}d(g(F{\rm U}_i^{(m)\ast}),F(G{\rm U}_i^{(m)\ast}))=0\;\;{\rm for~each}~i\in I_n.\eqno(11)$$
Now, owing to (11), we have
\begin{eqnarray*}
\Delta_n(GF_\ast{\rm U}^{(m)},F_\ast G{\rm U}^{(m)})&=&\frac{1}{n}\sum\limits_{i=1}^{n}d(g(F{\rm U}_i^{(m)\ast}),F(G{\rm U}_i^{(m)\ast}))\\
&\rightarrow& 0\;{\rm as}\;n\rightarrow \infty.
\end{eqnarray*}
It follows that $F_{\ast}$ and $G$ are O-compatible in ordered
metric space $(X^n,\Delta_n,\sqsubseteq_{\iota_n})$. In the similar
manner, one can prove the same for ordered metric space
$(X^n,\nabla_n,\sqsubseteq_{\iota_n})$.\\

The procedures of the proofs of parts (iii) and (iv) are similar to
Lemma 7 and the part (v) and hence the proof is left for readers.\\

(v) Take a sequence $\{{\rm U}^{(m)}\}\subset X^n$ and a ${\rm U}\in
X^n$ such that $\{G{\rm U}^{(m)}\}$ is monotone (w.r.t. partial
order $\sqsubseteq_{\iota_n}$) and
$$G({\rm U}^{(m)})\stackrel{\Delta_n}{\longrightarrow} G({\rm U})\;\;\;({\rm resp.}~G({\rm U}^{(m)})\stackrel{\nabla_n}{\longrightarrow} G({\rm U})).$$
Write U$^{(m)}=(x^{(m)}_1,x^{(m)}_2,...,x^{(m)}_n)$ and
U$=(x_1,x_2,...,x_n)$. Then, by using Propositions 6 and 7, we
obtain
$$g(x_i^{(m)})\uparrow\downarrow  g(x_i)\;{\rm for~each}~i\in
I_n.$$ It follows for each $i\in I_n$ that
$$g(x_{i_1}^{(m)})\uparrow\downarrow g(x_{i_1}),g(x_{i_2}^{(m)})\uparrow\downarrow g(x_{i_2}),..., g(x_{i_n}^{(m)})\uparrow\downarrow g(x_{i_n}). \eqno(12)$$ Using (12) and $(g,{\rm
O})$-continuity of $F$, we get
$$F(x_{i_1}^{(m)},x_{i_2}^{(m)},...,x_{i_n}^{(m)})\stackrel{d}{\longrightarrow}F(x_{i_1},x_{i_2},...,x_{i_n})$$
so that $$F({\rm U}^{(m)\ast}_i)\stackrel{d}{\longrightarrow}F({\rm U}^{\ast}_i)\;{\rm for~each}~i\in I_n,$$
 which, by using Proposition 6 gives rise
$$F_\ast({\rm U}^{(m)})\stackrel{\Delta_n}{\longrightarrow}F_\ast({\rm U})\;\;\;({\rm resp.}~F_\ast({\rm U}^{(m)})\stackrel{\nabla_n}{\longrightarrow} F_\ast({\rm U})).$$
Hence, $F_{\ast}$ is $(G,{\rm O})$-continuous in both ordered
metric spaces $(X^n,\Delta_n,\sqsubseteq_{\iota_n})$ and
$(X^n,\nabla_n,\sqsubseteq_{\iota_n})$.\\

(vi) Take a sequence $\{{\rm U}^{(m)}\}\subset E^n$ and a ${\rm
U}\in E^n$ such that $\{{\rm U}^{(m)}\}$ is monotone (w.r.t. partial
order $\sqsubseteq_{\iota_n}$) and
$${\rm U}^{(m)}\stackrel{\Delta_n}{\longrightarrow} {\rm U}\;\;\;({\rm resp.}~{\rm U}^{(m)}\stackrel{\nabla_n}{\longrightarrow} {\rm U}).$$
Write U$^{(m)}=(x^{(m)}_1,x^{(m)}_2,...,x^{(m)}_n)$ and
U$=(x_1,x_2,...,x_n)$. Then, by Proposition 6, we obtain
$$x_i^{(m)}\stackrel{d}{\longrightarrow} x_i\;{\rm for~each}~i\in
I_n.\eqno(13)$$ Now, there are two possibilities:\\
Case $(a):$ If $\{{\rm U}^{(m)}\}$ is increasing, then for all $m,l
\in \mathbb{N}_0$ with $m<l$, we have
$${\rm U}^{(m)}\sqsubseteq_{\iota_n}{\rm U}^{(l)},$$
or equivalently,
$$x^{(m)}_i)\preceq x^{(l)}_i\;{\rm {for\; each}
}\; i\in A\; {\rm {and} }\;x^{(m)}_i\succeq x^{(l)}_i\; {\rm {for\;
each} }\; i\in B.\eqno(14)$$ On combining (13) and (14), we obtain
$$x^{(m)}_i\uparrow x_i\;{\rm {for\; each}
}\; i\in A\; {\rm {and} }\;x^{(m)}_i\downarrow x_i\; {\rm {for\;
each} }\; i\in B,$$ which on using {\it g-MCB} property of
$(E,d,\preceq)$, gives rise
$$g(x^{(m)}_i)\preceq g(x_i)\;{\rm {for\; each}
}\; i\in A\; {\rm {and} }\;g(x^{(m)}_i)\succeq g(x_i)\; {\rm {for\;
each} }\; i\in B,$$ or equivalently,
$$G({\rm U}^{(m)})\sqsubseteq_{\iota_n}G({\rm U}).$$
It follows that $(E^n,\Delta_n,\sqsubseteq_{\iota_n})$ (resp.
$(E^n,\nabla_n,\sqsubseteq_{\iota_n})$) has {\it G-ICU} property.\\

Case $(b):$ If $\{{\rm U}^{(m)}\}$ is decreasing, then for all $m,l
\in \mathbb{N}_0$ with $m<l$, we have
$${\rm U}^{(m)}\sqsupseteq_{\iota_n}{\rm U}^{(l)},$$
or equivalently,
$$x^{(m)}_i)\succeq x^{(l)}_i\;{\rm {for\; each}
}\; i\in A\; {\rm {and} }\;x^{(m)}_i\preceq x^{(l)}_i\; {\rm {for\;
each} }\; i\in B.\eqno(15)$$ On combining (13) and (15), we obtain
$$x^{(m)}_i\downarrow x_i\;{\rm {for\; each}
}\; i\in A\; {\rm {and} }\;x^{(m)}_i\uparrow x_i\; {\rm {for\; each}
}\; i\in B,$$ which on using {\it g-MCB} property of
$(E,d,\preceq)$, gives rise
$$g(x^{(m)}_i)\succeq g(x_i)\;{\rm {for\; each}
}\; i\in A\; {\rm {and} }\;g(x^{(m)}_i)\preceq g(x_i)\; {\rm {for\;
each} }\; i\in B,$$ or equivalently,
$$G({\rm U}^{(m)})\sqsupseteq_{\iota_n}G({\rm U}).$$
It follows that $(E^n,\Delta_n,\sqsubseteq_{\iota_n})$ (resp.
$(E^n,\nabla_n,\sqsubseteq_{\iota_n})$) has {\it G-DCL} property.
Hence, in both the cases, $(E^n,\Delta_n,\sqsubseteq_{\iota_n})$
(resp. $(E^n,\nabla_n,\sqsubseteq_{\iota_n})$) has {\it G-MCB}
property.\\

(vii) This result is directly followed from (vi) by setting $g=I,$
the identity mapping.\\

\section{Multi-tupled Coincidence Theorems for Compatible Mappings}
\label{SC:Multi-tupled Coincidence Theorems for Compatible Mappings}
In this section, we prove the results regarding the existence and
uniqueness of $\ast$-coincidence points in
ordered metric spaces for O-compatible mappings.\\

\noindent{\bf Theorem 1.} Let $(X,d,\preceq)$ be an ordered metric
space, $E$ an O-complete subspace of $X$ and
$\ast\in\mathcal{U}_{\iota_n}$. Let $F:X^{n}\rightarrow X$ and
$g:X\to X $ be two mappings. Suppose that the following conditions
hold:
\begin{enumerate}
\item [{(i)}] $F(X^n)\subseteq g(X)\cap E$,
\item [{(ii)}] $F$ has $\iota_n$-mixed $g$-monotone property,
\item [{(iii)}] $F$ and $g$ are $(\ast,{\rm
O})$-compatible,
\item [{(iv)}] $g$ is O-continuous,
\item [{(v)}] either $F$ is O-continuous or $(E,d,\preceq)$ has {\it g-MCB}
property,
\item [{(vi)}] there exist $x^{(0)}_1,x^{(0)}_2,...,x^{(0)}_n
\in X$ such that
$${\begin{cases} g(x^{(0)}_{i}) \preceq
F(x^{(0)}_{i_1},x^{(0)}_{i_2},...,x^{(0)}_{i_n})\;{\rm for~ each}~
i\in A\cr g(x^{(0)}_{i}) \succeq
F(x^{(0)}_{i_1},x^{(0)}_{i_2},...,x^{(0)}_{i_n})\;{\rm for~ each}~
i\in B\cr\end{cases}}$$ or
$${\begin{cases} g(x^{(0)}_{i}) \succeq
F(x^{(0)}_{i_1},x^{(0)}_{i_2},...,x^{(0)}_{i_n})\;{\rm for~ each}~
i\in A\cr g(x^{(0)}_{i}) \preceq
F(x^{(0)}_{i_1},x^{(0)}_{i_2},...,x^{(0)}_{i_n})\;{\rm for~ each}~
i\in B,\cr\end{cases}}$$
\item [{(vii)}] there
exists $\varphi\in \Omega$ such that
$$\frac{1}{n}\sum\limits_{i=1}^{n}d(F(x_{i_1},x_{i_2},...,x_{i_n}),F(y_{i_1},y_{i_2},...,y_{i_n}))\leq\varphi\Big(\frac{1}{n}\sum\limits_{i=1}^{n}d(gx_i,gy_i)\Big)$$
for all $x_{1},x_{2},...,x_{n},y_{1},y_{2},...,y_{n}\in X$ with
[$g(x_i)\preceq g(y_i)$ for each $i\in A$ and $g(x_i)\succeq g(y_i)$
for each $i\in B$] or [$g(x_i)\succeq g(y_i)$ for each $i\in A$ and
$g(x_i)\preceq g(y_i)$ for each $i\in B$],
\end{enumerate}
or alternately
\begin{enumerate}
\item [{(vii$^\prime$)}] there
exists $\varphi\in \Omega$ such that
$$\max\limits_{i\in I_n}d(F(x_{i_1},x_{i_2},...,x_{i_n}),F(y_{i_1},y_{i_2},...,y_{i_n}))\leq\varphi\Big(\max\limits_{i\in I_n}d(gx_i,gy_i)\Big)$$
for all $x_{1},x_{2},...,x_{n},y_{1},y_{2},...,y_{n}\in X$ with
[$g(x_i)\preceq g(y_i)$ for each $i\in A$ and $g(x_i)\succeq g(y_i)$
for each $i\in B$] or [$g(x_i)\succeq g(y_i)$ for each $i\in A$ and
$g(x_i)\preceq g(y_i)$ for each $i\in B$].
\end{enumerate}
Then $F$ and $g$ have an $\ast$-coincidence point.\\

\noindent{\bf Proof.} We can induce two metrics $\Delta_n$ and
$\nabla_n$, patrial order $\sqsubseteq_{\iota_n}$ and two
self-mappings $F_\ast$ and $G$ on $X^n$ defined as in Section 4. By
item (i) of Lemma 8, both ordered metric subspaces
$(E^n,\Delta_n,\sqsubseteq_{\iota_n})$ and
$(E^n,\nabla_n,\sqsubseteq_{\iota_n})$ are O-complete. Further,
\begin{enumerate}
\item [{(i)}] implies that $F_\ast(X^n)\subseteq G(X^n)\cap E^n$ by item (i) of Lemma 3,
\item[{(ii)}] implies that $F_\ast$ is $G$-increasing in ordered set $(X^n,\sqsubseteq_{\iota_n})$ by Lemma 5,
\item [{(iii)}] implies that $F_\ast$ and $G$ are O-compatible in both $(X^n,\Delta_n,\sqsubseteq_{\iota_n})$ and $(X^n,\nabla_n,\sqsubseteq_{\iota_n})$ by
item (ii) of Lemma 8,
\item [{(iv)}] implies that $G$ is O-continuous in both $(X^n,\Delta_n,\sqsubseteq_{\iota_n})$ and
$(X^n,\nabla_n,\sqsubseteq_{\iota_n})$ by item (iii) of Lemma 8,
\item [{(v)}] implies that either $F_\ast$ is O-continuous in both $(X^n,\Delta_n,\sqsubseteq_{\iota_n})$ and $(X^n,\nabla_n,\sqsubseteq_{\iota_n})$ or both
$(E^n,\Delta_n,\sqsubseteq_{\iota_n})$ and
$(E^n,\nabla_n,\sqsubseteq_{\iota_n})$ have $G$-{\it MCB} property
by items (iv) and (vi) of Lemma 8
\item [{(vi)}] is equivalent to $G({\rm
U}^{(0)})\sqsubseteq_{\iota_n}F_\ast({\rm U}^{(0)})$ or $G({\rm
U}^{(0)})\sqsupseteq_{\iota_n}F_\ast({\rm U}^{(0)})$ where
U$^{(0)}=(x^{(0)}_1,x^{(0)}_2,...,x^{(0)}_n) \in X^n$,
\item [{(vii)}] means that $\Delta_n(F_\ast{\rm U},F_\ast{\rm V})\leq \varphi(\Delta_n(G{\rm U},G{\rm
V}))$ for all U=$(x_{1},x_{2},...,x_{n})$,
V=$(y_{1},y_{2},...,y_{n})\in X^n$ with $G({\rm U})\sqsubseteq_{\iota_n}G({\rm V})$ or $G({\rm U})\sqsupseteq_{\iota_n}G({\rm V})$,
\item [{(vii$^\prime$)}] means that $\nabla_n(F_\ast{\rm U},F_\ast{\rm V})\leq \varphi(\nabla_n(G{\rm U},G{\rm
V}))$ for all U=$(x_{1},x_{2},...,x_{n})$,
V=$(y_{1},y_{2},...,y_{n})\in X^n$ with $G({\rm U})\sqsubseteq_{\iota_n}G({\rm V})$ or $G({\rm U})\sqsupseteq_{\iota_n}G({\rm V})$.
\end{enumerate}
Therefore, the conditions (i)-(vii) of Lemma 1 are satisfied in the
context of ordered metric space
$(X^n,\Delta_n,\sqsubseteq_{\iota_n})$ or
$(X^n,\nabla_n,\sqsubseteq_{\iota_n})$ and two self-mappings
$F_\ast$ and $G$ on $X^n$. Thus, by Lemma 1, $F_\ast$ and $G$ have a
coincidence point, which is a $\ast$-coincidence point of $F$ and
$g$ by item (iii) of Lemma 3.\\

\noindent{\bf Corollary 1.} Let $(X,d,\preceq)$ be an O-complete
ordered metric space, $F:X^{n}\rightarrow X$ and $g:X\to X $ two
mappings and $\ast\in\mathcal{U}_{\iota_n}$. Suppose that the
following conditions hold:
\begin{enumerate}
\item [{(i)}] $F(X^n)\subseteq g(X)$,
\item [{(ii)}] $F$ has $\iota_n$-mixed $g$-monotone property,
\item [{(iii)}] $F$ and $g$ are $(\ast,{\rm
O})$-compatible,
\item [{(iv)}] $g$ is O-continuous,
\item [{(v)}] either $F$ is O-continuous or $(X,d,\preceq)$ has {\it g-MCB}
property,
\item [{(vi)}] there exist $x^{(0)}_1,x^{(0)}_2,...,x^{(0)}_n
\in X$ such that
$${\begin{cases} g(x^{(0)}_{i}) \preceq
F(x^{(0)}_{i_1},x^{(0)}_{i_2},...,x^{(0)}_{i_n})\;{\rm for~ each}~
i\in A\cr g(x^{(0)}_{i}) \succeq
F(x^{(0)}_{i_1},x^{(0)}_{i_2},...,x^{(0)}_{i_n})\;{\rm for~ each}~
i\in B\cr\end{cases}}$$ or
$${\begin{cases} g(x^{(0)}_{i}) \succeq
F(x^{(0)}_{i_1},x^{(0)}_{i_2},...,x^{(0)}_{i_n})\;{\rm for~ each}~
i\in A\cr g(x^{(0)}_{i}) \preceq
F(x^{(0)}_{i_1},x^{(0)}_{i_2},...,x^{(0)}_{i_n})\;{\rm for~ each}~
i\in B,\cr\end{cases}}$$
\item [{(vii)}] there
exists $\varphi\in \Omega$ such that
$$\frac{1}{n}\sum\limits_{i=1}^{n}d(F(x_{i_1},x_{i_2},...,x_{i_n}),F(y_{i_1},y_{i_2},...,y_{i_n}))\leq\varphi\Big(\frac{1}{n}\sum\limits_{i=1}^{n}d(gx_i,gy_i)\Big)$$
for all $x_{1},x_{2},...,x_{n},y_{1},y_{2},...,y_{n}\in X$ with
[$g(x_i)\preceq g(y_i)$ for each $i\in A$ and $g(x_i)\succeq g(y_i)$
for each $i\in B$] or [$g(x_i)\succeq g(y_i)$ for each $i\in A$ and
$g(x_i)\preceq g(y_i)$ for each $i\in B$],
\end{enumerate}
or alternately
\begin{enumerate}
\item [{(vii$^\prime$)}] there
exists $\varphi\in \Omega$ such that
$$\max\limits_{i\in I_n}d(F(x_{i_1},x_{i_2},...,x_{i_n}),F(y_{i_1},y_{i_2},...,y_{i_n}))\leq\varphi\Big(\max\limits_{i\in I_n}d(gx_i,gy_i)\Big)$$
for all $x_{1},x_{2},...,x_{n},y_{1},y_{2},...,y_{n}\in X$ with
[$g(x_i)\preceq g(y_i)$ for each $i\in A$ and $g(x_i)\succeq g(y_i)$
for each $i\in B$] or [$g(x_i)\succeq g(y_i)$ for each $i\in A$ and
$g(x_i)\preceq g(y_i)$ for each $i\in B$].
\end{enumerate}
Then $F$ and $g$ have an $\ast$-coincidence point.\\

On using Remarks 2, 4, 8 and 9, we obtain a natural version of
Theorem 1 as a
consequence, which runs below:\\

\noindent{\bf Corollary 2.} Theorem 1 remains true if the usual
metrical terms namely: completeness,
$\ast$-compatibility/commutativity and continuity are used
instead of their respective O-analogous.\\

As increasing requirement on $g$ together with {\it MCB} property
implies {\it g-MCB} property, therefore the following consequence of
Theorem 1 is immediately.\\

\noindent{\bf Corollary 3.} Theorem 1 remains true if we replace the
condition (v) by the following condition:
\begin{enumerate}
\item [{(v)}$^\prime$] $g$ is increasing and $(E,d,\preceq)$ has {\it MCB}
property.\\
\end{enumerate}

\noindent{\bf Corollary 4.} Theorem 1 remains true if we replace the
condition (vii) by the following condition:
\begin{enumerate}
\item [{(vii)}$^\prime$] there exists $\varphi\in \Omega$ such that
$$d(F(x_1,x_2,...,x_n),F(y_1,y_2,...,y_n))\leq\varphi\Big(\frac{1}{n}\sum\limits_{i=1}^{n}d(gx_i,gy_i)\Big)$$
for all $x_{1},x_{2},...,x_{n},y_{1},y_{2},...,y_{n}\in X$ with
[$g(x_i)\preceq g(y_i)$ for each $i\in A$ and $g(x_i)\succeq g(y_i)$
for each $i\in B$] or [$g(x_i)\succeq g(y_i)$ for each $i\in A$ and
$g(x_i)\preceq g(y_i)$ for each $i\in B$] provided that $\ast$ is
permuted.
\end{enumerate}

\noindent{\bf Proof.} Set U=$(x_{1},x_{2},...,x_{n})$,
V=$(y_{1},y_{2},...,y_{n})$ then we have $G({\rm
U})\sqsubseteq_{\iota_n}G({\rm V})$ or $G({\rm
U})\sqsupseteq_{\iota_n}G({\rm V})$. As $G({\rm U})$ and $G({\rm
V})$ are comparable, for each $i\in I_n$, $G({\rm U}^\ast_i)$ and
$G({\rm V}^\ast_i)$ are comparable w.r.t. partial order
$\sqsubseteq_{\iota_n}$ (owing to Lemma 4). Applying the contractivity condition
(vii)$^\prime$ on these points and using Lemma 6, for each $i\in
I_n$, we obtain
\begin{eqnarray*}
d(F(x_{i_1},x_{i_2},...,x_{i_n}),F(y_{i_1},y_{i_2},...,y_{i_n}))&\leq& \varphi\Big(\frac{1}{n}\sum\limits_{k=1}^{n}d(gx_{i_k},gy_{i_k})\Big)\\
&=&
\varphi\Big(\frac{1}{n}\sum\limits_{j=1}^{n}d(gx_j,gy_j)\Big)\;{\rm
as\;} \ast\;{\rm is\; permuted}
\end{eqnarray*}
so that
$$d(F(x_{i_1},x_{i_2},...,x_{i_n}),F(y_{i_1},y_{i_2},...,y_{i_n}))\leq\varphi\Big(\frac{1}{n}\sum\limits_{j=1}^{n}d(gx_j,gy_j)\Big)\;{\rm for~each~}i\in I_n.$$
Taking summation over $i\in I_n$ on both the sides of above
inequality, we obtain
$$\sum\limits_{i=1}^{n}d(F(x_{i_1},x_{i_2},...,x_{i_n}),F(y_{i_1},y_{i_2},...,y_{i_n}))\leq n\varphi\Big(\frac{1}{n}\sum\limits_{j=1}^{n}d(gx_j,gy_j)\Big)$$
so that
$$\frac{1}{n}\sum\limits_{i=1}^{n}d(F(x_{i_1},x_{i_2},...,x_{i_n}),F(y_{i_1},y_{i_2},...,y_{i_n}))\leq\varphi\Big(\frac{1}{n}\sum\limits_{j=1}^{n}d(gx_j,gy_j)\Big)$$
for all $x_{1},x_{2},...,x_{n},y_{1},y_{2},...,y_{n}\in X$ with
[$g(x_i)\preceq g(y_i)$ for each $i\in A$ and $g(x_i)\succeq g(y_i)$
for each $i\in B$] or [$g(x_i)\succeq g(y_i)$ for each $i\in A$ and
$g(x_i)\preceq g(y_i)$ for each $i\in B$].\\

Therefore, the contractivity condition (vii) of Theorem 1 holds and
hence Theorem 1 is
applicable.\\

\noindent{\bf Corollary 5.} Theorem 1 remains true if we replace the
condition (vii$^\prime$) by the following condition:
\begin{enumerate}
\item [{(vii$^\prime$)}$^\prime$] there exists $\varphi\in \Omega$ such that
$$d(F(x_1,x_2,...,x_n),F(y_1,y_2,...,y_n))\leq\varphi\Big(\max\limits_{i\in I_n}d(gx_i,gy_i)\Big)$$
for all $x_{1},x_{2},...,x_{n},y_{1},y_{2},...,y_{n}\in X$ with
[$g(x_i)\preceq g(y_i)$ for each $i\in A$ and $g(x_i)\succeq g(y_i)$
for each $i\in B$] or [$g(x_i)\succeq g(y_i)$ for each $i\in A$ and
$g(x_i)\preceq g(y_i)$ for each $i\in B$] provided that either
$\ast$ is permuted or $\varphi$ is increasing on $[0,\infty)$.
\end{enumerate}

\noindent{\bf Proof.} Set U=$(x_{1},x_{2},...,x_{n})$,
V=$(y_{1},y_{2},...,y_{n})$, then similar to previous corollary, for
each $i\in I_n$, $G({\rm U}^\ast_i)$ and $G({\rm V}^\ast_i)$ are
comparable w.r.t. partial order $\sqsubseteq_{\iota_n}$. Applying
the contractivity condition (vii$^\prime$)}$^\prime$ on these points
and using Lemma 6, for each $i\in I_n$, we obtain
$$d(F(x_{i_1},x_{i_2},...,x_{i_n}),F(y_{i_1},y_{i_2},...,y_{i_n}))\leq \varphi\Big(\max\limits_{k\in
I_n}d(gx_{i_k},gy_{i_k})\Big)$$
\indent\hspace{8.8cm}${\begin{cases}=\varphi\Big(\max\limits_{j\in
I_n}d(gx_j,gy_j)\Big)\;{\rm if\;} \ast\;{\rm is\; permuted,}\cr
\hspace{0.0in}\leq\varphi\Big(\max\limits_{j\in
I_n}d(gx_j,gy_j)\Big)\;{\rm if\;}\varphi\;{\rm is\;
inceasing.}\cr\end{cases}}$\\
so that
$$d(F(x_{i_1},x_{i_2},...,x_{i_n}),F(y_{i_1},y_{i_2},...,y_{i_n}))\leq\varphi\Big(\max\limits_{j\in
I_n}d(gx_i,gy_i)\Big)\;{\rm for~each~}i\in I_n.$$ Taking maximum
over $i\in I_n$ on both the sides of above inequality, we obtain
$$\max\limits_{i\in
I_n}d(F(x_{i_1},x_{i_2},...,x_{i_n}),F(y_{i_1},y_{i_2},...,y_{i_n}))\leq\varphi\Big(\max\limits_{j\in
I_n}d(gx_j,gy_j)\Big)$$ for all
$x_{1},x_{2},...,x_{n},y_{1},y_{2},...,y_{n}\in X$ with
[$g(x_i)\preceq g(y_i)$ for each $i\in A$ and $g(x_i)\succeq g(y_i)$
for each $i\in B$] or [$g(x_i)\succeq g(y_i)$ for each $i\in A$ and
$g(x_i)\preceq g(y_i)$ for each $i\in B$].\\

Therefore, the contractivity condition (vii$^\prime$) of Theorem 1
holds and hence Theorem 1 is
applicable.\\

Now, we present multi-tupled coincidence theorems for linear and
generalized
linear contractions.\\

\noindent{\bf Corollary 6.} In addition to the hypotheses (i)-(vi)
of Theorem 1, suppose that one of the following conditions holds:
\begin{enumerate}
\item [{(viii)}] there exists $\alpha\in [0,1)$ such that
$$\frac{1}{n}\sum\limits_{i=1}^{n}d(F(x_{i_1},x_{i_2},...,x_{i_n}),F(y_{i_1},y_{i_2},...,y_{i_n}))\leq\frac{\alpha}{n} \sum\limits_{i=1}^{n}d(gx_i,gy_i)$$
for all $x_{1},x_{2},...,x_{n},y_{1},y_{2},...,y_{n}\in X$ with
[$g(x_i)\preceq g(y_i)$ for each $i\in A$ and $g(x_i)\succeq g(y_i)$
for each $i\in B$] or [$g(x_i)\succeq g(y_i)$ for each $i\in A$ and
$g(x_i)\preceq g(y_i)$ for each $i\in B$],
\item [{(ix)}] there exists $\alpha\in [0,1)$ such that
$$\max\limits_{i\in I_n}d(F(x_{i_1},x_{i_2},...,x_{i_n}),F(y_{i_1},y_{i_2},...,y_{i_n}))\leq\alpha\max\limits_{i\in I_n}d(gx_i,gy_i)$$
for all $x_{1},x_{2},...,x_{n},y_{1},y_{2},...,y_{n}\in X$ with
[$g(x_i)\preceq g(y_i)$ for each $i\in A$ and $g(x_i)\succeq g(y_i)$
for each $i\in B$] or [$g(x_i)\succeq g(y_i)$ for each $i\in A$ and
$g(x_i)\preceq g(y_i)$ for each $i\in B$].
\end{enumerate}
Then $F$ and $g$ have an $\ast$-coincidence point.\\

\noindent{\bf Proof.} On setting $\varphi(t)=\alpha t$ with
$\alpha\in [0,1)$ in Theorem 1, we get our
result.\\

\noindent{\bf Corollary 7.} In addition to the hypotheses (i)-(vi)
of Theorem 1, suppose that one of the following conditions holds:
\begin{enumerate}
\item [{(x)}] there exists $\alpha\in [0,1)$ such that
$$d(F(x_1,x_2,...,x_n),F(y_1,y_2,...,y_n))\leq\alpha\max\limits_{i\in I_n}d(gx_i,gy_i)$$
for all $x_{1},x_{2},...,x_{n},y_{1},y_{2},...,y_{n}\in X$ with
[$g(x_i)\preceq g(y_i)$ for each $i\in A$ and $g(x_i)\succeq g(y_i)$
for each $i\in B$] or [$g(x_i)\succeq g(y_i)$ for each $i\in A$ and
$g(x_i)\preceq g(y_i)$ for each $i\in B$],
\item [{(xi)}] there exists $\alpha_1,\alpha_2,...,\alpha_n\in [0,1)$ with $\sum\limits_{i=1}^n\alpha_i<1$ such that
$$d(F(x_1,x_2,...,x_n),F(y_1,y_2,...,y_n))\leq\sum\limits_{i=1}^{n}\alpha_i d(gx_i,gy_i)$$
for all $x_{1},x_{2},...,x_{n},y_{1},y_{2},...,y_{n}\in X$ with
[$g(x_i)\preceq g(y_i)$ for each $i\in A$ and $g(x_i)\succeq g(y_i)$
for each $i\in B$] or [$g(x_i)\succeq g(y_i)$ for each $i\in A$ and
$g(x_i)\preceq g(y_i)$ for each $i\in B$],
\item [{(xii)}] there exists $\alpha\in [0,1)$ such that
$$d(F(x_1,x_2,...,x_n),F(y_1,y_2,...,y_n))\leq\frac{\alpha}{n} \sum\limits_{i=1}^{n}d(gx_i,gy_i)$$
for all $x_{1},x_{2},...,x_{n},y_{1},y_{2},...,y_{n}\in X$ with
[$g(x_i)\preceq g(y_i)$ for each $i\in A$ and $g(x_i)\succeq g(y_i)$
for each $i\in B$] or [$g(x_i)\succeq g(y_i)$ for each $i\in A$ and
$g(x_i)\preceq g(y_i)$ for each $i\in B$].
\end{enumerate}
Then $F$ and $g$ have an $\ast$-coincidence point.\\

\noindent{\bf Proof.} Setting $\varphi(t)=\alpha t$ with $\alpha\in
[0,1)$ in Corollary 5, we get the result
corresponding to the contractivity condition (x). Notice that here $\varphi$ is increasing on $[0,\infty)$.\\

To prove the result corresponding to (xi), let
$\beta=\sum\limits_{i=1}^n\alpha_i<1$, then we have
\begin{align*}
d(F(x_1,x_2,...,x_n),F(y_1,y_2,...,y_n))&\leq \sum\limits_{i=1}^n\alpha_id(gx_i,gy_i)\nonumber\\
& \leq \Big(\sum\limits_{i=1}^n\alpha_i\Big) \max\limits_{j\in I_n}d(gx_j,gy_j)\nonumber\\
& =\beta \max\limits_{j\in I_n}d(gx_j,gy_j)\nonumber\\
\end{align*}

so that our result follows from the result corresponding to (x).\\

Finally, setting $\alpha_i=\frac{\alpha}{n}$ for all $i\in I_n,$
where $\alpha\in [0,1)$ in (xi), we get the result corresponding to
(xii). Notice that here $\sum\limits_{i=1}^n\alpha_i=\alpha<1$.\\

Now, we present uniqueness results corresponding to Theorem 1, which
run as follows:\\

\noindent{\bf Theorem 2.} In addition to the hypotheses of Theorem
1, suppose that for every pair $(x_1,x_2,...,x_n)$,
$(y_1,y_2,...,y_n)\in X^n$, there exists $(z_1,z_2,...,z_n)\in X^n$
such that $(gz_1,gz_2,...,gz_n)$ is comparable to
$(gx_1,gx_2,...,gx_n)$ and $(gy_1,gy_2,...,gy_n)$ w.r.t. partial
order $\sqsubseteq_{\iota_n}$, then $F$ and $g$ have a unique point
of $\ast$-coincidence, which remains also a unique common
$\ast$-fixed point.\\

\noindent{\bf Proof.} Set U=$(x_{1},x_{2},...,x_{n})$,
V=$(y_{1},y_{2},...,y_{n})$ and W=$(z_{1},z_{2},...,z_{n})$, then by
one of our assumptions $G({\rm W})$ is comparable to $G({\rm U})$
and $G({\rm V})$. Therefore, all the conditions of Lemma 1 are
satisfied. Hence, by Lemma 1, $F_\ast$ and $G$ have a unique point
of coincidence as well as a unique common fixed point, which is
indeed a unique point of $\ast$-coincidence as well as a unique
common $\ast$-fixed point of
$F$ and $g$ by items (iv) and (v) of Lemma 3.\\

\noindent{\bf Theorem 3.} In addition to the hypotheses of Theorem
2, suppose that $g$ is one-one, then $F$ and $g$ have
a unique $\ast$-coincidence point.\\

\noindent{\bf Proof.} Let U=$(x_1,x_2,...,x_n)$ and
V=$(y_1,y_2,...,y_n)$ be two $\ast$-coincidence point of $F$ and $g$
then then using Theorem 2, we obtain
$$(gx_1,gx_2,...,gx_n)=(gy_1,gy_2,...,gy_n)$$
or equivalently
$$g(x_i)=g(y_i)\;{\rm for~each~}i\in I_n.$$
As $g$ is one-one, we have
$$x_i=y_i\;{\rm for~each~}i\in I_n.$$
It follows that U=V, $i.e.$, $F$ and $g$ have a unique $\ast$-coincidence
point.\\

\section{Multi-tupled Coincidence Theorems without Compatibility of mappings}
\label{SC:Multi-tupled Coincidence Theorems without Compatibility of
mappings}

In this section, we prove the results regarding the existence and
uniqueness of $\ast$-coincidence points in an ordered metric space
$X$ for a pair of mappings $F:X^{n}\rightarrow X$ and $g:X\to X$,
which are not necessarily O-compatible.\\

\noindent{\bf Theorem 4.} Let $(X,d,\preceq)$ be an ordered metric
space, $E$ an O-complete subspace of $X$ and
$\ast\in\mathcal{U}_{\iota_n}$. Let $F:X^{n}\rightarrow X$ and
$g:X\to X $ be two mappings. Suppose that the following conditions
hold:
\begin{enumerate}
\item [{(i)}] $F(X^n)\subseteq E\subseteq g(X)$,
\item [{(ii)}]  $F$ has $\iota_n$-mixed $g$-monotone property,
\item [{(iii)}] either $F$ is $(g,{\rm O})$-continuous or $F$ and $g$ are continuous or $(E,d,\preceq)$ has {\it MCB} property,
\item [{(iv)}] there exist $x^{(0)}_1,x^{(0)}_2,...,x^{(0)}_n
\in X$ such that
$${\begin{cases} g(x^{(0)}_{i}) \preceq
F(x^{(0)}_{i_1},x^{(0)}_{i_2},...,x^{(0)}_{i_n})\;{\rm for~ each}~
i\in A\cr g(x^{(0)}_{i}) \succeq
F(x^{(0)}_{i_1},x^{(0)}_{i_2},...,x^{(0)}_{i_n})\;{\rm for~ each}~
i\in B\cr\end{cases}}$$ or
$${\begin{cases} g(x^{(0)}_{i}) \succeq
F(x^{(0)}_{i_1},x^{(0)}_{i_2},...,x^{(0)}_{i_n})\;{\rm for~ each}~
i\in A\cr g(x^{(0)}_{i}) \preceq
F(x^{(0)}_{i_1},x^{(0)}_{i_2},...,x^{(0)}_{i_n})\;{\rm for~ each}~
i\in B,\cr\end{cases}}$$
\item [{(v)}] there
exists $\varphi\in \Omega$ such that
$$\frac{1}{n}\sum\limits_{i=1}^{n}d(F(x_{i_1},x_{i_2},...,x_{i_n}),F(y_{i_1},y_{i_2},...,y_{i_n}))=\varphi\Big(\frac{1}{n}\sum\limits_{i=1}^{n}d(gx_i,gy_i)\Big)$$
for all $x_{1},x_{2},...,x_{n},y_{1},y_{2},...,y_{n}\in X$ with
[$g(x_i)\preceq g(y_i)$ for each $i\in A$ and $g(x_i)\succeq g(y_i)$
for each $i\in B$] or [$g(x_i)\succeq g(y_i)$ for each $i\in A$ and
$g(x_i)\preceq g(y_i)$ for each $i\in B$],
\end{enumerate}
or alternately
\begin{enumerate}
\item [{(v$^\prime$)}] there
exists $\varphi\in \Omega$ such that
$$\max\limits_{i\in I_n}d(F(x_{i_1},x_{i_2},...,x_{i_n}),F(y_{i_1},y_{i_2},...,y_{i_n}))=\varphi\Big(\max\limits_{i\in I_n}d(gx_i,gy_i)\Big)$$
for all $x_{1},x_{2},...,x_{n},y_{1},y_{2},...,y_{n}\in X$ with
[$g(x_i)\preceq g(y_i)$ for each $i\in A$ and $g(x_i)\succeq g(y_i)$
for each $i\in B$] or [$g(x_i)\succeq g(y_i)$ for each $i\in A$ and
$g(x_i)\preceq g(y_i)$ for each $i\in B$].
\end{enumerate}
Then $F$ and $g$ have an $\ast$-coincidence point.\\

\noindent{\bf Proof.} We can induce two metrics $\Delta_n$ and
$\nabla_n$, patrial order $\sqsubseteq_{\iota_n}$ and two
self-mappings $F_\ast$ and $G$ on $X^n$ defined as in Section 4. By
item (i) of Lemma 8, both ordered metric subspaces
$(E^n,\Delta_n,\sqsubseteq_{\iota_n})$ and
$(E^n,\nabla_n,\sqsubseteq_{\iota_n})$ are O-complete. Further,
\begin{enumerate}
\item [{(i)}] implies that $F_\ast(X^n)\subseteq E^n\subseteq G(X^n)$ by item (ii) of Lemma 3,
\item [{(ii)}] implies that $F_\ast$ is $G$-increasing in ordered set $(X^n,\sqsubseteq_{\iota_n})$ by Lemma 5,
\item [{(iii)}] implies that either $F_\ast$ is $(G,{\rm O})$-continuous in both $(X^n,\Delta_n,\sqsubseteq_{\iota_n})$ and $(X^n,\nabla_n,\sqsubseteq_{\iota_n})$ or
$F_\ast$ and $G$ are continuous in both $(X^n,\Delta_n)$ and
$(X^n,\nabla_n)$ or both $(E^n,\Delta_n,\sqsubseteq_{\iota_n})$ and
$(E^n,\nabla_n,\sqsubseteq_{\iota_n})$ have {\it MCB} property by
Lemma 7 and items (v) and (vii) of Lemma 8,
\item [{(iv)}] is equivalent to $G({\rm
U}^{(0)})\sqsubseteq_{\iota_n}F_\ast({\rm U}^{(0)})$ or $G({\rm
U}^{(0)})\sqsupseteq_{\iota_n}F_\ast({\rm U}^{(0)})$ where
U$^{(0)}=(x^{(0)}_1,x^{(0)}_2,...,x^{(0)}_n) \in X^n$,
\item [{(v)}] means that $\Delta_n(F_\ast{\rm U},F_\ast{\rm V})\leq \varphi(\Delta_n(G{\rm U},G{\rm
V}))$ for all U=$(x_{1},x_{2},...,x_{n})$,
V=$(y_{1},y_{2},...,y_{n})\in X^n$ with $G({\rm U})\sqsubseteq_{\iota_n}G({\rm V})$ or $G({\rm U})\sqsupseteq_{\iota_n}G({\rm V})$,
\item [{(v$^\prime$)}] means that $\nabla_n(F_\ast{\rm U},F_\ast{\rm V})\leq \varphi(\nabla_n(G{\rm U},G{\rm
V}))$ for all U=$(x_{1},x_{2},...,x_{n})$,
V=$(y_{1},y_{2},...,y_{n})\in X^n$ with $G({\rm U})\sqsubseteq_{\iota_n}G({\rm V})$ or $G({\rm U})\sqsupseteq_{\iota_n}G({\rm V})$.
\end{enumerate}
Therefore, the conditions (i)-(v) of Lemma 2 are satisfied in the
context of ordered metric space
$(X^n,\Delta_n,\sqsubseteq_{\iota_n})$ or
$(X^n,\nabla_n,\sqsubseteq_{\iota_n})$ and two self-mappings
$F_\ast$ and $G$ on $X^n$. Thus, by Lemma 2, $F_\ast$ and $G$ have a
coincidence point, which is a $\ast$-coincidence point of $F$ and
$g$ by item (iii) of Lemma 3.\\

\noindent{\bf Corollary 8.} Let $(X,d,\preceq)$ be an O-complete
ordered metric space, $F:X^{n}\rightarrow X$ and $g:X\to X $ two
mappings and $\ast\in\mathcal{U}_{\iota_n}$. Suppose that the
following conditions hold:
\begin{enumerate}
\item [{(i)}] either $g$ is onto or there exists an ${\rm O}$-closed subspace $E$ of $X$ such that $f(X)\subseteq E\subseteq g(X)$,
\item [{(ii)}]  $F$ has $\iota_n$-mixed $g$-monotone property,
\item [{(iii)}] either $F$ is $(g,{\rm O})$-continuous or $F$ and $g$ are continuous or $(X,d,\preceq)$ has {\it MCB} property,
\item [{(iv)}] there exist $x^{(0)}_1,x^{(0)}_2,...,x^{(0)}_n
\in X$ such that
$${\begin{cases} g(x^{(0)}_{i}) \preceq
F(x^{(0)}_{i_1},x^{(0)}_{i_2},...,x^{(0)}_{i_n})\;{\rm for~ each}~
i\in A\cr g(x^{(0)}_{i}) \succeq
F(x^{(0)}_{i_1},x^{(0)}_{i_2},...,x^{(0)}_{i_n})\;{\rm for~ each}~
i\in B\cr\end{cases}}$$ or
$${\begin{cases} g(x^{(0)}_{i}) \succeq
F(x^{(0)}_{i_1},x^{(0)}_{i_2},...,x^{(0)}_{i_n})\;{\rm for~ each}~
i\in A\cr g(x^{(0)}_{i}) \preceq
F(x^{(0)}_{i_1},x^{(0)}_{i_2},...,x^{(0)}_{i_n})\;{\rm for~ each}~
i\in B,\cr\end{cases}}$$
\item [{(v)}] there
exists $\varphi\in \Omega$ such that
$$\frac{1}{n}\sum\limits_{i=1}^{n}d(F(x_{i_1},x_{i_2},...,x_{i_n}),F(y_{i_1},y_{i_2},...,y_{i_n}))=\varphi\Big(\frac{1}{n}\sum\limits_{i=1}^{n}d(gx_i,gy_i)\Big)$$
for all $x_{1},x_{2},...,x_{n},y_{1},y_{2},...,y_{n}\in X$ with
[$g(x_i)\preceq g(y_i)$ for each $i\in A$ and $g(x_i)\succeq g(y_i)$
for each $i\in B$] or [$g(x_i)\succeq g(y_i)$ for each $i\in A$ and
$g(x_i)\preceq g(y_i)$ for each $i\in B$],
\end{enumerate}
or alternately
\begin{enumerate}
\item [{(v$^\prime$)}] there
exists $\varphi\in \Omega$ such that
$$\max\limits_{i\in I_n}d(F(x_{i_1},x_{i_2},...,x_{i_n}),F(y_{i_1},y_{i_2},...,y_{i_n}))=\varphi\Big(\max\limits_{i\in I_n}d(gx_i,gy_i)\Big)$$
for all $x_{1},x_{2},...,x_{n},y_{1},y_{2},...,y_{n}\in X$ with
[$g(x_i)\preceq g(y_i)$ for each $i\in A$ and $g(x_i)\succeq g(y_i)$
for each $i\in B$] or [$g(x_i)\succeq g(y_i)$ for each $i\in A$ and
$g(x_i)\preceq g(y_i)$ for each $i\in B$].
\end{enumerate}
Then $F$ and $g$ have an $\ast$-coincidence point.\\

\noindent{\bf Proof.} The result corresponding to first part of
(i) ($i.e$ in case that $g$ is onto) is followed by taking
$E=X=g(X)$ in Theorem 4. While the result corresponding to second alternating part of (i)
($i.e$ in case that $E$ is O-closed)  is
followed by using Proposition 1.\\

On using Remarks 2,3 and 8, we obtain a natural version of Theorem 4
as a consequence, which runs below:\\

\noindent{\bf Corollary 9.} Theorem 4 (also Corollary 8) remains
true if the usual metrical terms namely: completeness, closedness
and $g$-continuity are used
instead of their respective O-analogous.\\

Similar to Corollaries 4-6, the following consequences of
Theorem 4 hold.\\

\noindent{\bf Corollary 10.} Theorem 4 remains true if we replace
the condition (v) by the following condition:
\begin{enumerate}
\item [{(v)}$^\prime$] there exists $\varphi\in \Omega$ such that
$$d(F(x_1,x_2,...,x_n),F(y_1,y_2,...,y_n))\leq\varphi\Big(\frac{1}{n}\sum\limits_{i=1}^{n}d(gx_i,gy_i)\Big)$$
for all $x_{1},x_{2},...,x_{n},y_{1},y_{2},...,y_{n}\in X$ with
[$g(x_i)\preceq g(y_i)$ for each $i\in A$ and $g(x_i)\succeq g(y_i)$
for each $i\in B$] or [$g(x_i)\succeq g(y_i)$ for each $i\in A$ and
$g(x_i)\preceq g(y_i)$ for each $i\in B$] provided that $\ast$ is
permuted.\\
\end{enumerate}

\noindent{\bf Corollary 11.} Theorem 4 remains true if we replace
the condition (v$^\prime$) by the following condition:
\begin{enumerate}
\item [{(v$^\prime$)}$^\prime$] there exists $\varphi\in \Omega$ such that
$$d(F(x_1,x_2,...,x_n),F(y_1,y_2,...,y_n))\leq\varphi\Big(\max\limits_{i\in I_n}d(gx_i,gy_i)\Big)$$
for all $x_{1},x_{2},...,x_{n},y_{1},y_{2},...,y_{n}\in X$ with
[$g(x_i)\preceq g(y_i)$ for each $i\in A$ and $g(x_i)\succeq g(y_i)$
for each $i\in B$] or [$g(x_i)\succeq g(y_i)$ for each $i\in A$ and
$g(x_i)\preceq g(y_i)$ for each $i\in B$] provided that either
$\ast$ is permuted or $\varphi$ is increasing on $[0,\infty)$.\\
\end{enumerate}

\noindent{\bf Corollary 12.} In addition to the hypotheses (i)-(iv)
of Theorem 4, suppose that one of the following conditions holds:
\begin{enumerate}
\item [{(vi)}] there exists $\alpha\in [0,1)$ such that
$$\frac{1}{n}\sum\limits_{i=1}^{n}d(F(x_{i_1},x_{i_2},...,x_{i_n}),F(y_{i_1},y_{i_2},...,y_{i_n}))\leq\frac{\alpha}{n} \sum\limits_{i=1}^{n}d(gx_i,gy_i)$$
for all $x_{1},x_{2},...,x_{n},y_{1},y_{2},...,y_{n}\in X$ with
[$g(x_i)\preceq g(y_i)$ for each $i\in A$ and $g(x_i)\succeq g(y_i)$
for each $i\in B$] or [$g(x_i)\succeq g(y_i)$ for each $i\in A$ and
$g(x_i)\preceq g(y_i)$ for each $i\in B$],
\item [{(vii)}] there exists $\alpha\in [0,1)$ such that
$$\max\limits_{i\in I_n}d(F(x_{i_1},x_{i_2},...,x_{i_n}),F(y_{i_1},y_{i_2},...,y_{i_n}))\leq\alpha\max\limits_{i\in I_n}d(gx_i,gy_i)$$
for all $x_{1},x_{2},...,x_{n},y_{1},y_{2},...,y_{n}\in X$ with
[$g(x_i)\preceq g(y_i)$ for each $i\in A$ and $g(x_i)\succeq g(y_i)$
for each $i\in B$] or [$g(x_i)\succeq g(y_i)$ for each $i\in A$ and
$g(x_i)\preceq g(y_i)$ for each $i\in B$].
\end{enumerate}
Then $F$ and $g$ have an $\ast$-coincidence point.\\

\noindent{\bf Corollary 13.} In addition to the hypotheses (i)-(iv)
of Theorem 4, suppose that one of the following conditions holds:
\begin{enumerate}
\item [{(viii)}] there exists $\alpha\in [0,1)$ such that
$$d(F(x_1,x_2,...,x_n),F(y_1,y_2,...,y_n))\leq\alpha\max\limits_{i\in I_n}d(gx_i,gy_i)$$
for all $x_{1},x_{2},...,x_{n},y_{1},y_{2},...,y_{n}\in X$ with
[$g(x_i)\preceq g(y_i)$ for each $i\in A$ and $g(x_i)\succeq g(y_i)$
for each $i\in B$] or [$g(x_i)\succeq g(y_i)$ for each $i\in A$ and
$g(x_i)\preceq g(y_i)$ for each $i\in B$],
\item [{(ix)}] there exists $\alpha_1,\alpha_2,...,\alpha_n\in [0,1)$ with $\sum\limits_{i=1}^n\alpha_i<1$ such that
$$d(F(x_1,x_2,...,x_n),F(y_1,y_2,...,y_n))\leq\sum\limits_{i=1}^{n}\alpha_i d(gx_i,gy_i)$$
for all $x_{1},x_{2},...,x_{n},y_{1},y_{2},...,y_{n}\in X$ with
[$g(x_i)\preceq g(y_i)$ for each $i\in A$ and $g(x_i)\succeq g(y_i)$
for each $i\in B$] or [$g(x_i)\succeq g(y_i)$ for each $i\in A$ and
$g(x_i)\preceq g(y_i)$ for each $i\in B$],
\item [{(x)}] there exists $\alpha\in [0,1)$ such that
$$d(F(x_1,x_2,...,x_n),F(y_1,y_2,...,y_n))\leq\frac{\alpha}{n} \sum\limits_{i=1}^{n}d(gx_i,gy_i)$$
for all $x_{1},x_{2},...,x_{n},y_{1},y_{2},...,y_{n}\in X$ with
[$g(x_i)\preceq g(y_i)$ for each $i\in A$ and $g(x_i)\succeq g(y_i)$
for each $i\in B$] or [$g(x_i)\succeq g(y_i)$ for each $i\in A$ and
$g(x_i)\preceq g(y_i)$ for each $i\in B$].
\end{enumerate}
Then $F$ and $g$ have an $\ast$-coincidence point.\\

Now, we present uniqueness results corresponding to Theorem 4, which
run as follows:\\

\noindent{\bf Theorem 5.} In addition to the hypotheses of Theorem
4, suppose that for every pair $(x_1,x_2,...,x_n)$,
$(y_1,y_2,...,y_n)\in X^n$, there exists $(z_1,z_2,...,z_n)\in X^n$
such that $(gz_1,gz_2,...,gz_n)$ is comparable to
$(gx_1,gx_2,...,gx_n)$ and $(gy_1,gy_2,...,gy_n)$ w.r.t. partial
order $\sqsubseteq_{\iota_n}$, then $F$ and $g$ have a unique point
of $\ast$-coincidence.\\

\noindent{\bf Proof.} Set U=$(x_{1},x_{2},...,x_{n})$,
V=$(y_{1},y_{2},...,y_{n})$ and W=$(z_{1},z_{2},...,z_{n})$, then by
one of our assumptions $G({\rm W})$ is comparable to $G({\rm U})$
and $G({\rm V})$. Therefore, all the conditions of Lemma 2 are
satisfied. Hence, by Lemma 2, $F_\ast$ and $G$ have a unique point
of coincidence, which is indeed a unique point of $\ast$-coincidence
of $F$ and $g$ by item (iv) of Lemma 3.\\

\noindent{\bf Theorem 6.} In addition to the hypotheses of Theorem
5, suppose that $g$ is one-one, then $F$ and $g$ have
a unique $\ast$-coincidence point.\\

\noindent{\bf Proof.} The proof of Theorem 6 is similar to that of Theorem 3.\\

\noindent{\bf Theorem 7.} In addition to the hypotheses of Theorem
5, suppose that $F$ and $g$ are weakly $\ast$-compatible, then $F$
and $g$ have a unique common $\ast$-fixed point.\\

{\noindent\bf{Proof.}} Let $(x_{1},x_{2},...,x_{n})$ be a
$\ast$-coincidence point of $F$ and $g$. Write
$F(x_{i_1},x_{i_2},...,x_{i_n})=g(x_i)=\overline{x}_i$ for each
$i\in I_n$. Then, by Proposition 5,
$(\overline{x}_1,\overline{x}_2,...,\overline{x}_n)$ being a point
of $\ast$-coincidence of $F$ and $g$ is also a $\ast$-coincidence
point of $F$ and $g.$ It follows from Theorem 5 that
$$(gx_1,gx_2,...,gx_n)=(g\overline{x}_1,g\overline{x}_2,...,g\overline{x}_n)$$ $i.e.$,
$\overline{x_i}=g(\overline{x_i})$ for each $i\in I_n$, which for
each $i\in I_n$ yields that
$$F(x_{i_1},x_{i_2},...,x_{i_n})=g(\overline{x_i})=\overline{x_i}.$$
Hence, $(\overline{x}_1,\overline{x}_2,...,\overline{x}_n)$ is a
common $\ast$-fixed point of $F$ and $g$. To prove uniqueness,
assume that $(x^*_1,x^*_2,...,x^*_n)$ is another common $\ast$-fixed
point of $F$ and $g$. Then again from Theorem 5,
$$(gx^*_1,gx^*_2,...,gx^*_n)=(g\overline{x}_1,g\overline{x}_2,...,g\overline{x}_n)$$
$i.e.$
$$(x^*_1,x^*_2,...,x^*_n)=(\overline{x}_1,\overline{x}_2,...,\overline{x}_n).$$
This completes the proof.\\

\section{Multi-tupled Fixed Point Theorems}
\label{SC:Multi-tupled Fixed Point Theorems}

On particularizing $g=I$, the identity mapping on $X$, in the
foregoing results contained in Sections 5 and 6, we obtain the
corresponding $\ast$-fixed point results, which run as follows:\\

\noindent{\bf Theorem 8.} Let $(X,d,\preceq)$ be an ordered metric
space, $F:X^{n}\rightarrow X$ a mapping and
$\ast\in\mathcal{U}_{\iota_n}$. Let $E$ be an O-complete subspace of
$X$ such that $F(X^n)\subseteq E$. Suppose that the following
conditions hold:
\begin{enumerate}
\item [{(i)}] $F$ has $\iota_n$-mixed monotone property,
\item [{(ii)}] either $F$ is O-continuous or $(E,d,\preceq)$ has {\it MCB} property,
\item [{(iii)}] there exist $x^{(0)}_1,x^{(0)}_2,...,x^{(0)}_n
\in X$ such that
$${\begin{cases} x^{(0)}_{i} \preceq
F(x^{(0)}_{i_1},x^{(0)}_{i_2},...,x^{(0)}_{i_n})\;{\rm for~ each}~
i\in A\cr x^{(0)}_{i} \succeq
F(x^{(0)}_{i_1},x^{(0)}_{i_2},...,x^{(0)}_{i_n})\;{\rm for~ each}~
i\in B\cr\end{cases}}$$ or
$${\begin{cases} x^{(0)}_{i} \succeq
F(x^{(0)}_{i_1},x^{(0)}_{i_2},...,x^{(0)}_{i_n})\;{\rm for~ each}~
i\in A\cr x^{(0)}_{i} \preceq
F(x^{(0)}_{i_1},x^{(0)}_{i_2},...,x^{(0)}_{i_n})\;{\rm for~ each}~
i\in B,\cr\end{cases}}$$
\item [{(iv)}] there
exists $\varphi\in \Omega$ such that
$$\frac{1}{n}\sum\limits_{i=1}^{n}d(F(x_{i_1},x_{i_2},...,x_{i_n}),F(y_{i_1},y_{i_2},...,y_{i_n}))\leq\varphi\Big(\frac{1}{n}\sum\limits_{i=1}^{n}d(x_i,y_i)\Big)$$
for all $x_{1},x_{2},...,x_{n},y_{1},y_{2},...,y_{n}\in X$ with
[$x_i\preceq y_i$ for each $i\in A$ and $x_i\succeq y_i$ for each
$i\in B$] or [$x_i\succeq y_i$ for each $i\in A$ and $x_i\preceq
y_i$ for each $i\in B$],
\end{enumerate}
or alternately
\begin{enumerate}
\item [{(iv$^\prime$)}] there
exists $\varphi\in \Omega$ such that
$$\max\limits_{i\in I_n}d(F(x_{i_1},x_{i_2},...,x_{i_n}),F(y_{i_1},y_{i_2},...,y_{i_n}))\leq\varphi\Big(\max\limits_{i\in I_n}d(x_i,y_i)\Big)$$
for all $x_{1},x_{2},...,x_{n},y_{1},y_{2},...,y_{n}\in X$ with
[$x_i\preceq y_i$ for each $i\in A$ and $x_i\succeq y_i$ for each
$i\in B$] or [$x_i\succeq y_i$ for each $i\in A$ and $x_i\preceq
y_i$ for each $i\in B$].
\end{enumerate}
Then $F$ has an $\ast$-fixed point.\\

\noindent{\bf Corollary 14.} Let $(X,d,\preceq)$ be an O-complete
ordered metric space, $F:X^{n}\rightarrow X$ a mapping and
$\ast\in\mathcal{U}_{\iota_n}$. Suppose that the following
conditions hold:
\begin{enumerate}
\item [{(i)}] $F$ has $\iota_n$-mixed monotone property,
\item [{(ii)}] either $F$ is O-continuous or $(X,d,\preceq)$ has {\it MCB} property,
\item [{(iii)}] there exist $x^{(0)}_1,x^{(0)}_2,...,x^{(0)}_n
\in X$ such that
$${\begin{cases} x^{(0)}_{i} \preceq
F(x^{(0)}_{i_1},x^{(0)}_{i_2},...,x^{(0)}_{i_n})\;{\rm for~ each}~
i\in A\cr x^{(0)}_{i} \succeq
F(x^{(0)}_{i_1},x^{(0)}_{i_2},...,x^{(0)}_{i_n})\;{\rm for~ each}~
i\in B\cr\end{cases}}$$ or
$${\begin{cases} x^{(0)}_{i} \succeq
F(x^{(0)}_{i_1},x^{(0)}_{i_2},...,x^{(0)}_{i_n})\;{\rm for~ each}~
i\in A\cr x^{(0)}_{i} \preceq
F(x^{(0)}_{i_1},x^{(0)}_{i_2},...,x^{(0)}_{i_n})\;{\rm for~ each}~
i\in B,\cr\end{cases}}$$
\item [{(iv)}] there
exists $\varphi\in \Omega$ such that
$$\frac{1}{n}\sum\limits_{i=1}^{n}d(F(x_{i_1},x_{i_2},...,x_{i_n}),F(y_{i_1},y_{i_2},...,y_{i_n}))\leq\varphi\Big(\frac{1}{n}\sum\limits_{i=1}^{n}d(x_i,y_i)\Big)$$
for all $x_{1},x_{2},...,x_{n},y_{1},y_{2},...,y_{n}\in X$ with
[$x_i\preceq y_i$ for each $i\in A$ and $x_i\succeq y_i$ for each
$i\in B$] or [$x_i\succeq y_i$ for each $i\in A$ and $x_i\preceq
y_i$ for each $i\in B$],
\end{enumerate}
or alternately
\begin{enumerate}
\item [{(iv$^\prime$)}] there
exists $\varphi\in \Omega$ such that
$$\max\limits_{i\in I_n}d(F(x_{i_1},x_{i_2},...,x_{i_n}),F(y_{i_1},y_{i_2},...,y_{i_n}))\leq\varphi\Big(\max\limits_{i\in I_n}d(x_i,y_i)\Big)$$
for all $x_{1},x_{2},...,x_{n},y_{1},y_{2},...,y_{n}\in X$ with
[$x_i\preceq y_i$ for each $i\in A$ and $x_i\succeq y_i$ for each
$i\in B$] or [$x_i\succeq y_i$ for each $i\in A$ and $x_i\preceq
y_i$ for each $i\in B$].
\end{enumerate}
Then $F$ has an $\ast$-fixed point.\\

\noindent{\bf Corollary 15.} Theorem 8 remains true if the usual
metrical terms namely: completeness and continuity are used
instead of their respective O-analogous.\\

\noindent{\bf Corollary 16.} Theorem 8 remains true if we replace
the condition (iv) by the following condition:
\begin{enumerate}
\item [{(iv)}$^\prime$] there exists $\varphi\in \Omega$ such that
$$d(F(x_1,x_2,...,x_n),F(y_1,y_2,...,y_n))\leq\varphi\Big(\frac{1}{n}\sum\limits_{i=1}^{n}d(x_i,y_i)\Big)$$
for all $x_{1},x_{2},...,x_{n},y_{1},y_{2},...,y_{n}\in X$ with
[$x_i\preceq y_i$ for each $i\in A$ and $x_i\succeq y_i$ for each
$i\in B$] or [$x_i\succeq y_i$ for each $i\in A$ and $x_i\preceq
y_i$ for each $i\in B$] provided that $\ast$ is permuted.\\
\end{enumerate}

\noindent{\bf Corollary 17.} Theorem 8 remains true if we replace
the condition (iv$^\prime$) by the following condition:
\begin{enumerate}
\item [{(iv$^\prime$)}$^\prime$] there exists $\varphi\in \Omega$ such that
$$d(F(x_1,x_2,...,x_n),F(y_1,y_2,...,y_n))\leq\varphi\Big(\max\limits_{i\in I_n}d(x_i,y_i)\Big)$$
for all $x_{1},x_{2},...,x_{n},y_{1},y_{2},...,y_{n}\in X$ with
[$x_i\preceq y_i$ for each $i\in A$ and $x_i\succeq y_i$ for each
$i\in B$] or [$x_i\succeq y_i$ for each $i\in A$ and $x_i\preceq
y_i$ for each $i\in B$] provided that either $\ast$ is permuted or
$\varphi$ is increasing on $[0,\infty)$.\\
\end{enumerate}

\noindent{\bf Corollary 18.} Theorem 8 remains true if we replace
the condition (iv) by the following condition:
\begin{enumerate}
\item [{(v)}] there exists $\alpha\in [0,1)$ such that
$$\frac{1}{n}\sum\limits_{i=1}^{n}d(F(x_{i_1},x_{i_2},...,x_{i_n}),F(y_{i_1},y_{i_2},...,y_{i_n}))\leq\frac{\alpha}{n} \sum\limits_{i=1}^{n}d(x_i,y_i)$$
for all $x_{1},x_{2},...,x_{n},y_{1},y_{2},...,y_{n}\in X$ with
[$x_i\preceq y_i$ for each $i\in A$ and $x_i\succeq y_i)$ for each
$i\in B$] or [$x_i\succeq y_i$ for each $i\in A$ and $x_i\preceq
y_i$ for each $i\in B$],
\item [{(vi)}] there exists $\alpha\in [0,1)$ such that
$$\max\limits_{i\in I_n}d(F(x_{i_1},x_{i_2},...,x_{i_n}),F(y_{i_1},y_{i_2},...,y_{i_n}))\leq\alpha\max\limits_{i\in I_n}d(x_i,y_i)$$
for all $x_{1},x_{2},...,x_{n},y_{1},y_{2},...,y_{n}\in X$ with
[$x_i\preceq y_i$ for each $i\in A$ and $x_i\succeq y_i$ for each
$i\in B$] or [$x_i\succeq y_i$ for each $i\in A$ and
$x_i\preceq y_i$ for each $i\in B$].\\
\end{enumerate}

\noindent{\bf Corollary 19.} Theorem 8 remains true if we replace
the condition (iv) by the following condition:
\begin{enumerate}
\item [{(vii)}] there exists $\alpha\in [0,1)$ such that
$$d(F(x_1,x_2,...,x_n),F(y_1,y_2,...,y_n))\leq\alpha\max\limits_{i\in I_n}d(x_i,y_i)$$
for all $x_{1},x_{2},...,x_{n},y_{1},y_{2},...,y_{n}\in X$ with
[$x_i\preceq y_i$ for each $i\in A$ and $x_i\succeq y_i$ for each
$i\in B$] or [$x_i\succeq y_i$ for each $i\in A$ and $x_i\preceq
y_i$ for each $i\in B$].
\item [{(viii)}] there exist $\alpha_1,\alpha_2,...,\alpha_n\in [0,1)$ with $\sum\limits_{i=1}^n\alpha_i<1$ such that
$$d(F(x_1,x_2,...,x_n),F(y_1,y_2,...,y_n))\leq\sum\limits_{i=1}^{n}\alpha_i d(x_i,y_i)$$
for all $x_{1},x_{2},...,x_{n},y_{1},y_{2},...,y_{n}\in X$ with
[$x_i\preceq y_i$ for each $i\in A$ and $x_i\succeq y_i$ for each
$i\in B$] or [$x_i\succeq y_i$ for each $i\in A$ and $x_i\preceq
y_i$ for each $i\in B$].
\item [{(ix)}] there exists $\alpha\in [0,1)$ such that
$$d(F(x_1,x_2,...,x_n),F(y_1,y_2,...,y_n))\leq\frac{\alpha}{n} \sum\limits_{i=1}^{n}d(x_i,y_i)$$
for all $x_{1},x_{2},...,x_{n},y_{1},y_{2},...,y_{n}\in X$ with
[$x_i\preceq y_i$ for each $i\in A$ and $x_i\succeq y_i$ for each
$i\in B$] or [$x_i\succeq y_i$ for each $i\in A$ and $x_i\preceq
y_i$ for each $i\in B$].\\
\end{enumerate}

\noindent{\bf Theorem 9.} In addition to the hypotheses of Theorem
8, suppose that for every pair $(x_1,x_2,...,x_n)$,
$(y_1,y_2,...,y_n)\in X^n$, there exists $(z_1,z_2,...,z_n)\in X^n$
such that $(z_1,z_2,...,z_n)$ is comparable to $(x_1,x_2,...,x_n)$
and $(y_1,y_2,...,y_n)$ w.r.t. partial order
$\sqsubseteq_{\iota_n}$, then $F$ has a unique $\ast$-fixed point.\\

\section{Particular Cases}
\label{SC:Particular Cases}

\subsection{Coupled Fixed/Coincidence Theorems} \hspace{1cm}\\

 On setting $n=2$, $\iota_2=\{\{1\},\{2\}\}$ and $\ast=\left[\begin{matrix}
1 &2\\
2 &1\\
\end{matrix}\right]$ in Corollaries 2,3,4,10,16,18,19, we obtain
the following results ($i.e.$ Corollaries 20-26).\\

\noindent{\bf Corollary 20} (Bhaskar and Lakshmikantham \cite{C1}).
Let $(X,d,\preceq)$ be an ordered complete metric space and
$F:X^{2}\rightarrow X$ a mapping. Suppose that the following
conditions hold:
\begin{enumerate}
\item [{(i)}] $F$ has mixed monotone property,
\item [{(ii)}] either $F$ is continuous or $(X,d,\preceq)$ has {\it MCB} property,
\item [{(iii)}] there exist $x^{(0)},y^{(0)}
\in X$ such that
$$ x^{(0)} \preceq
F(x^{(0)},y^{(0)})\;{\rm and }\; y^{(0)} \succeq
F(y^{(0)},x^{(0)})$$
\item [{(iv)}] there exists $\alpha\in [0,1)$ such that
$$d(F(x,y),F(u,v))\leq\frac{\alpha}{2} [d(x,u)+d(y,v)]$$
for all $x,y,u,v\in X$ with $x\preceq u$ and $y\succeq v$.
\end{enumerate}

Then $F$ has a coupled fixed point.\\

\noindent{\bf Corollary 21} (Berinde \cite{C6}). Corollary 20
remains true if we replace conditions (iii) and (iv) by the
following respective conditions:
\begin{enumerate}
\item [{(iii)$^\prime$}] there exist $x^{(0)},y^{(0)}
\in X$ such that
$$ x^{(0)} \preceq
F(x^{(0)},y^{(0)})\;{\rm and }\; y^{(0)} \succeq
F(y^{(0)},x^{(0)})$$ or
$$ x^{(0)} \succeq
F(x^{(0)},y^{(0)})\;{\rm and }\; y^{(0)} \preceq
F(y^{(0)},x^{(0)})$$
\item [{(iv)$^\prime$}] there exists $\alpha\in [0,1)$ such that
$$d(F(x,y),F(u,v))+d(F(y,x),F(v,u))\leq \alpha [d(x,u)+d(y,v)]$$
for all $x,y,u,v\in X$ with $x\preceq u$ and $y\succeq v$.\\
\end{enumerate}

\noindent{\bf Corollary 22} (Wu and Liu \cite{PGF12}, Sintunavarat
and Kumam \cite{C5}). Corollary 20 remains true if we replace
condition (iv) by the following condition:
\begin{enumerate}
\item [{(iv)$^{\prime\prime}$}] there exists $\varphi\in \Phi$ such that
$$d(F(x,y),F(u,v))\leq\varphi\bigg(\frac{d(x,u)+d(y,v)}{2}\bigg)$$
for all $x,y,u,v\in X$ with $x\preceq u$ and $y\succeq v$.\\
\end{enumerate}

\noindent{\bf Corollary 23} (Lakshmikantham and \'{C}iri\'{c}
\cite{C2}). Let $(X,d,\preceq)$ be an ordered complete metric space
and $F:X^{2}\rightarrow X$ and $g:X\to X $ two mappings. Suppose
that the following conditions hold:
\begin{enumerate}
\item [{(i)}] $F(X^2)\subseteq g(X)$,
\item [{(ii)}] $F$ has mixed $g$-monotone property,
\item [{(iii)}] $F$ and $g$ are commuting,
\item [{(iv)}] $g$ is continuous,
\item [{(v)}] either $F$ is continuous or $(X,d,\preceq)$ has $g$-{\it MCB}
property,
\item [{(vi)}] there exist $x^{(0)},y^{(0)}
\in X$ such that
$$ x^{(0)} \preceq
F(x^{(0)},y^{(0)})\;{\rm and }\; y^{(0)} \succeq
F(y^{(0)},x^{(0)})$$
\item [{(vii)}] there exists $\varphi\in \Phi$ such that
$$d(F(x,y),F(u,v))\leq\varphi\bigg(\frac{d(gx,gu)+d(yg,gv)}{2}\bigg)$$
for all $x,y,u,v\in X$ with $g(x)\preceq g(u)$ and $g(y)\succeq
g(v)$.
\end{enumerate}

Then $F$ and $g$ have a coupled coincidence point.\\

\noindent{\bf Corollary 24} (Choudhury and Kundu \cite{C3}).
Corollary 23 remains true if we replace conditions (iii), (iv) and
(v) by the following respective conditions:
\begin{enumerate}
\item [{(iii)$^\prime$}] $F$ and $g$ are compatible,
\item [{(iv)$^\prime$}] $g$ is continuous and increasing,
\item [{(v)$^\prime$}] either $F$ is continuous or $(X,d,\preceq)$ has {\it MCB}
property.\\
\end{enumerate}

\noindent{\bf Corollary 25} (Berinde \cite{C7}). Corollary 23
remains true if we replace conditions (vi) and (vii) by the
following respective conditions:
\begin{enumerate}
\item [{(vi)$^\prime$}] there exist $x^{(0)},y^{(0)}
\in X$ such that
$$ g(x^{(0)}) \preceq
F(x^{(0)},y^{(0)})\;{\rm and }\; g(y^{(0)}) \succeq
F(y^{(0)},x^{(0)})$$ or
$$ g(x^{(0)}) \succeq
F(x^{(0)},y^{(0)})\;{\rm and }\; g(y^{(0)}) \preceq
F(y^{(0)},x^{(0)})$$
\item [{(vii)$^\prime$}] there exists $\varphi\in \Phi$ such that
$$d(F(x,y),F(u,v))+d(F(y,x),F(v,u))\leq2\varphi\bigg(\frac{d(gx,gu)+d(gy,gv)}{2}\bigg)$$
for all $x,y,u,v\in X$ with $g(x)\preceq g(u)$ and $g(y)\succeq
g(v)$.\\
\end{enumerate}

{\noindent\bf{Corollary 26}} (Husain $et$ $al.$ \cite{C4},
Sintunavarat and Kumam \cite{C5}). Let $(X,d,\preceq)$ be an ordered
metric space and $F:X^2\rightarrow X$ and $g: X\rightarrow X$ two
mappings. Let $(gX,d)$ be complete subspace. Suppose that the
following conditions hold:
\begin{enumerate}
\item [{(i)}] $F(X^2)\subseteq g(X)$,
\item [{(ii)}] $F$ has mixed $g$-monotone property,
\item [{(iii)}] $g$ is continuous,
\item [{(iv)}] either $F$ is continuous or $(X,d,\preceq)$ has {\it MCB} property,
\item [{(v)}]there exist $x^{(0)},y^{(0)}
\in X$ such that
$$ x^{(0)} \preceq
F(x^{(0)},y^{(0)})\;{\rm and }\; y^{(0)} \succeq
F(y^{(0)},x^{(0)})$$
\item [{(vi)}] there exists $\varphi\in \Phi$ such that
$$d(F(x,y),F(u,v))\leq\varphi\bigg(\frac{d(gx,gu)+d(gy,gv)}{2}\bigg)$$
for all $x,y,u,v\in X$ with $g(x)\preceq g(u)$ and $g(y)\succeq
g(v)$.
\end{enumerate}
Then $F$ and $g$ have a coupled coincidence point.\\

\subsection{Tripled Fixed/Coincidence Theorems}\hspace{1cm}\\

On setting $n=3$, $\iota_3=\{\{1,3\},\{2\}\}$ and
$\ast=\left[\begin{matrix}
1 &2 &3\\
2 &1 &2\\
3 &2 &1\\
\end{matrix}\right]$ in Corollaries 2,3,5,7,9,13,19, we obtain
the following results ($i.e.$ Corollaries 27-32).\\

\noindent{\bf Corollary 27} (Berinde and Borcut \cite{T1}). Let
$(X,d,\preceq)$ be an ordered complete metric space and
$F:X^{3}\rightarrow X$ a mapping. Suppose that the following
conditions hold:
\begin{enumerate}
\item [{(i)}] $F$ has alternating mixed monotone property,
\item [{(ii)}] either $F$ is continuous or $(X,d,\preceq)$ has {\it MCB} property,
\item [{(iii)}] there exist $x^{(0)},y^{(0)},z^{(0)}
\in X$ such that
$$ x^{(0)} \preceq
F(x^{(0)},y^{(0)},z^{(0)}), y^{(0)}\succeq
F(y^{(0)},x^{(0)},y^{(0)})\;{\rm and }\; z^{(0)} \preceq
F(z^{(0)},y^{(0)},x^{(0)}),$$
\item [{(iv)}] there exist $\alpha,\beta,\gamma\in [0,1)$ with $\alpha+\beta+\gamma<1$ such that
$$d(F(x,y,z),F(u,v,w))\leq\alpha d(x,u)+\beta d(y,v)+ \gamma d(z,w)$$
for all $x,y,z,u,v,w\in X$ with $x\preceq u$, $y\succeq v$ and
$z\preceq w$.
\end{enumerate}
Then $F$ has a tripled fixed point (in the sense of Berinde and
Borcut \cite{T1}), $i.e.$, there exist $x,y,z\in X$ such that
$F(x,y,z)=x$, $F(y,x,y)=y$ and $F(z,y,x)=z$.\\

\noindent{\bf Corollary 28} (Borcut and Berinde \cite{T2}). Let
$(X,d,\preceq)$ be an ordered complete metric space and
$F:X^{3}\rightarrow X$ and $g:X\to X $ two mappings. Suppose that
the following conditions hold:
\begin{enumerate}
\item [{(i)}] $F(X^3)\subseteq g(X)$,
\item [{(ii)}] $F$ has alternating mixed $g$-monotone property,
\item [{(iii)}] $F$ and $g$ are commuting,
\item [{(iv)}] $g$ is continuous,
\item [{(v)}] either $F$ is continuous or $(X,d,\preceq)$ has $g$-{\it MCB}
property,
\item [{(vi)}] there exist $x^{(0)},y^{(0)},z^{(0)}
\in X$ such that
$$ g(x^{(0)}) \preceq
F(x^{(0)},y^{(0)},z^{(0)}), g(y^{(0)})\succeq
F(y^{(0)},x^{(0)},y^{(0)})\;{\rm and }\; g(z^{(0)}) \preceq
F(z^{(0)},y^{(0)},x^{(0)}),$$
\item [{(vii)}] there exist $\alpha,\beta,\gamma\in [0,1)$ with $\alpha+\beta+\gamma<1$ such that
$$d(F(x,y,z),F(u,v,w))\leq\alpha d(gx,gu)+\beta d(gy,gv)+ \gamma d(gz,gw)$$
for all $x,y,z,u,v,w\in X$ with $g(x)\preceq g(u)$, $g(y)\succeq
g(v)$ and $g(z)\preceq g(w)$.
\end{enumerate}

Then $F$ and $g$ have a tripled coincidence point (in the sense of
Berinde and Borcut \cite{T1}), $i.e.$, there exist $x,y,z\in X$ such
that
$F(x,y,z)=g(x)$, $F(y,x,y)=g(y)$ and $F(z,y,x)=g(z)$.\\

\noindent{\bf Corollary 29} (Borcut \cite{T3}). Corollary 28
remains true if we replace condition (vii) by the following
condition:
\begin{enumerate}
\item [{(vii)$^{\prime}$}] there exists $\varphi\in \Phi$ provided $\varphi$ is increasing
such that
$$d(F(x,y,z),F(u,v,w))\leq\varphi\big(\max \{d(gx,gu),d(gy,gv),d(gz,gw)\}\big)$$
for all $x,y,z,u,v,w\in X$ with $g(x)\preceq g(u)$, $g(y)\succeq
g(v)$ and $g(z)\preceq g(w)$.\\
\end{enumerate}

\noindent{\bf Corollary 30} (Choudhury $et\;al.$ \cite{T4}).
Corollary 29 remains true if we replace conditions (iii) and (v) by
the following conditions respectively:
\begin{enumerate}
\item [{(iii)$^{\prime}$}] $F$ and $g$ are compatible,
\item [{(v)$^{\prime}$}] either $F$ is continuous or $(X,d,\preceq)$ has {\it MCB}
property provided $g$ is increasing.\\
\end{enumerate}

{\noindent\bf{Corollary 31}} (Husain $et$ $al.$ \cite{C4}). Let
$(X,d,\preceq)$ be an ordered metric space and $F:X^3\rightarrow X$
and $g: X\rightarrow X$ two mappings. Let $(gX,d)$ be complete
subspace. Suppose that the following conditions hold:
\begin{enumerate}
\item [{(i)}] $F(X^3)\subseteq g(X)$,
\item [{(ii)}] $F$ has alternating mixed $g$-monotone property,
\item [{(iii)}] $g$ is continuous,
\item [{(iv)}] either $F$ is continuous or $(X,d,\preceq)$ has {\it MCB} property,
\item [{(v)}] there exist $x^{(0)},y^{(0)},z^{(0)}
\in X$ such that
$$ g(x^{(0)}) \preceq
F(x^{(0)},y^{(0)},z^{(0)}), g(y^{(0)})\succeq
F(y^{(0)},x^{(0)},y^{(0)})\;{\rm and }\; g(z^{(0)}) \preceq
F(z^{(0)},y^{(0)},x^{(0)}),$$
\item [{(vi)}] there exist $\alpha,\beta,\gamma\in [0,1)$ with $\alpha+\beta+\gamma<1$ such that
$$d(F(x,y,z),F(u,v,w))\leq\alpha d(gx,gu)+\beta d(gy,gv)+ \gamma d(gz,gw)$$
for all $x,y,z,u,v,w\in X$ with $g(x)\preceq g(u)$, $g(y)\succeq
g(v)$ and $g(z)\preceq g(w)$.
\end{enumerate}
Then $F$ and $g$ have a tripled coincidence point (in the sense of
Berinde and Borcut \cite{T1}), $i.e.$, there exist $x,y,z\in X$ such
that
$F(x,y,z)=g(x)$, $F(y,x,y)=g(y)$ and $F(z,y,x)=g(z)$.\\

\noindent{\bf Corollary 32} (Radenovi$\acute{\rm c}$ \cite{R3}). Let
$(X,d,\preceq)$ be an ordered metric space and $F:X^{3}\rightarrow
X$ and $g:X\to X $ two mappings. Suppose that the following
conditions hold:
\begin{enumerate}
\item [{(i)}] $F(X^3)\subseteq g(X)$,
\item [{(ii)}] $F$ has alternating mixed $g$-monotone property,
\item [{(iii)}] there exist $x^{(0)},y^{(0)},z^{(0)}
\in X$ such that
$$g(x^{(0)}) \preceq
F(x^{(0)},y^{(0)},z^{(0)}), g(y^{(0)})\succeq
F(y^{(0)},x^{(0)},y^{(0)})\;{\rm and }\; g(z^{(0)}) \preceq
F(z^{(0)},y^{(0)},x^{(0)})~{\rm or}$$
$$g(x^{(0)}) \succeq
F(x^{(0)},y^{(0)},z^{(0)}), g(y^{(0)})\preceq
F(y^{(0)},x^{(0)},y^{(0)})\;{\rm and }\; g(z^{(0)}) \succeq
F(z^{(0)},y^{(0)},x^{(0)}),$$
\item [{(iv)}] there exists $\varphi\in \Phi$ provided $\varphi$ is increasing
such that
$$\max \{d(F(x,y,z),F(u,v,w)),d(F(y,x,y),F(v,u,v)),d(F(z,y,x),F(w,v,u))\}$$
$$\leq\varphi\big(\max \{d(gx,gu),d(gy,gv),d(gz,gw)\}\big)$$ for all
$x,y,z,u,v,w\in X$ with [$g(x)\preceq g(u)$, $g(y)\succeq g(v)$ and
$g(z)\preceq g(w)$] or [$g(x)\succeq g(u)$, $g(y)\preceq g(v)$ and
$g(z)\succeq g(w)$],
\item [{(v)}] $F$ and $g$ are continuous and compatible and $(X,d)$ is complete, or
\item [{(v$^\prime$)}] $(X,d,\preceq)$ has {\it MCB} property and one of $F(X^3)$ or $g(X)$ is complete.
\end{enumerate}

Then $F$ and $g$ have a tripled coincidence point (in the sense of
Berinde and Borcut \cite{T1}), $i.e.$, there exist $x,y,z\in X$ such
that
$F(x,y,z)=g(x)$, $F(y,x,y)=g(y)$ and $F(z,y,x)=g(z)$.\\

On setting $n=3$, $\iota_3=\{\{1,3\},\{2\}\}$ and
$\ast=\left[\begin{matrix}
1 &2 &3\\
2 &3 &2\\
3 &2 &1\\
\end{matrix}\right]$ in Corollary 19, we obtain
the following result:\\

\noindent{\bf Corollary 33} (Wu and Liu \cite{TQ}). Let
$(X,d,\preceq)$ be an ordered complete metric space and
$F:X^{3}\rightarrow X$ a mapping. Suppose that the following
conditions hold:
\begin{enumerate}
\item [{(i)}] $F$ has alternating mixed monotone property,
\item [{(ii)}] either $F$ is continuous or $(X,d,\preceq)$ has {\it MCB} property,
\item [{(iii)}] there exist $x^{(0)},y^{(0)},z^{(0)}
\in X$ such that
$$ x^{(0)} \preceq
F(x^{(0)},y^{(0)},z^{(0)}), y^{(0)}\succeq
F(y^{(0)},z^{(0)},y^{(0)})\;{\rm and }\; z^{(0)} \preceq
F(z^{(0)},y^{(0)},x^{(0)}),$$
\item [{(iv)}] there exist $\alpha,\beta,\gamma\in [0,1)$ with $\alpha+\beta+\gamma<1$ such that
$$d(F(x,y,z),F(u,v,w))\leq\alpha d(x,u)+\beta d(y,v)+ \gamma d(z,w)$$
for all $x,y,z,u,v,w\in X$ with $x\preceq u$, $y\succeq v$ and
$z\preceq w$.
\end{enumerate}
Then $F$ has a tripled fixed point (in the sense of Wu and Liu
\cite{TQ}), $i.e.$, there exist $x,y,z\in X$ such that
$F(x,y,z)=x$, $F(y,z,y)=y$ and $F(z,y,x)=z$.\\

On setting $n=3$, $\iota_3=\{\{1,2\},\{3\}\}$ and
$\ast=\left[\begin{matrix}
1 &2 &3\\
2 &1 &3\\
3 &3 &2\\
\end{matrix}\right]$ in Corollary 19, we obtain
the following result:\\

\noindent{\bf Corollary 34} (Berzig and Samet \cite{HD1}). Let
$(X,d,\preceq)$ be an ordered complete metric space and
$F:X^{3}\rightarrow X$ a mapping. Suppose that the following
conditions hold:
\begin{enumerate}
\item [{(i)}] $F$ has 2-mixed monotone property,
\item [{(ii)}] either $F$ is continuous or $(X,d,\preceq)$ has {\it MCB} property,
\item [{(iii)}] there exist $x^{(0)},y^{(0)},z^{(0)}
\in X$ such that
$$ x^{(0)} \preceq
F(x^{(0)},y^{(0)},z^{(0)}), y^{(0)}\preceq
F(y^{(0)},x^{(0)},z^{(0)})\;{\rm and }\; z^{(0)} \succeq
F(z^{(0)},z^{(0)},y^{(0)}),$$
\item [{(iv)}] there exist $\alpha,\beta,\gamma\in [0,1)$ with $\alpha+\beta+\gamma<1$ such that
$$d(F(x,y,z),F(u,v,w))\leq\alpha d(x,u)+\beta d(y,v)+ \gamma d(z,w)$$
for all $x,y,z,u,v,w\in X$ with $x\preceq u$, $y\preceq v$ and
$z\succeq w$.
\end{enumerate}
Then $F$ has a tripled fixed point (in the sense of Berzig and Samet
\cite{HD1}), $i.e.$, there exist $x,y,z\in X$ such that $F(x,y,z)=x$, $F(y,x,z)=y$,
$F(z,z,y)=z.$

\subsection{Quadrupled Fixed/Coincidence Theorems} \hspace{1cm}\\

On setting $n=4$, $\iota_4=\{\{1,3\},\{2,4\}\}$ and
$\ast=\left[\begin{matrix}
1 &2 &3 &4\\
2 &3 &4 &1\\
3 &4 &1 &2\\
4 &1 &2 &3\\
\end{matrix}\right]$ in Corollaries 4,7,19, we obtain
the following results ($i.e.$ Corollaries 35-37).\\

\noindent{\bf Corollary 35} (Karapinar and Luong \cite{Q1}). Let
$(X,d,\preceq)$ be an ordered complete metric space and
$F:X^{4}\rightarrow X$ a mapping. Suppose that the following
conditions hold:
\begin{enumerate}
\item [{(i)}] $F$ has alternating mixed monotone property,
\item [{(ii)}] either $F$ is continuous or $(X,d,\preceq)$ has {\it MCB} property,
\item [{(iii)}] there exist $x^{(0)},y^{(0)},z^{(0)},w^{(0)}
\in X$ such that
$$ x^{(0)}\preceq F(x^{(0)},y^{(0)},z^{(0)},w^{(0)}),$$
$$y^{(0)}\succeq F(y^{(0)},z^{(0)},w^{(0)},x^{(0)}),$$
$$z^{(0)}\preceq F(z^{(0)},w^{(0)},x^{(0)},y^{(0)}),$$
$$w^{(0)}\succeq F(w^{(0)},x^{(0)},y^{(0)},z^{(0)}),$$
\item [{(iv)}] there exist $\alpha,\beta,\gamma,\delta\in [0,1)$ with $\alpha+\beta+\gamma+\delta<1$ such that
$$d(F(x,y,z,w),F(u,v,r,t))\leq\alpha d(gx,gu)+\beta d(y,v)+ \gamma d(z,r)+\delta d(w,t)$$
for all $x,y,z,w,u,v,r,t\in X$ with $x\preceq u$, $y\succeq v$,
$z\preceq r$ and $w\succeq t$.
\end{enumerate}

Then $F$ has a quadrupled fixed point (in the sense of Karapinar and
Luong \cite{Q1}), $i.e.$, there exist $x,y,z,w\in X$ such that
$$F(x,y,z,w)=x, F(y,z,w,x)=y,$$
$$F(z,w,x,y)=z, F(w,x,y,z)=w.$$

\noindent{\bf Corollary 36} (Liu \cite{Q3}). Let $(X,d,\preceq)$ be
an ordered complete metric space and $F:X^{4}\rightarrow X$ and
$g:X\to X $ two mappings. Suppose that the following conditions
hold:
\begin{enumerate}
\item [{(i)}] $F(X^4)\subseteq g(X)$,
\item [{(ii)}] $F$ has alternating mixed $g$-monotone property,
\item [{(iii)}] $F$ and $g$ are commuting,
\item [{(iv)}] $g$ is continuous,
\item [{(v)}] either $F$ is continuous or $(X,d,\preceq)$ has $g$-{\it MCB}
property,
\item [{(vi)}] there exist $x^{(0)},y^{(0)},z^{(0)},w^{(0)} \in X$ such that
$$g(x^{(0)}) \preceq
F(x^{(0)},y^{(0)},z^{(0)},w^{(0)}),$$ $$g(y^{(0)}) \succeq
F(y^{(0)},z^{(0)},w^{(0)},x^{(0)}),$$
$$g(z^{(0)}) \preceq
F(z^{(0)},w^{(0)},x^{(0)},y^{(0)}),$$ $$g(w^{(0)}) \succeq
F(w^{(0)},x^{(0)},y^{(0)},z^{(0)}),$$
\item [{(vii)}] there exist $\alpha,\beta,\gamma,\delta\in [0,1)$ with $\alpha+\beta+\gamma+\delta<1$ such that
$$d(F(x,y,z,w),F(u,v,r,t))\leq\alpha d(gx,gu)+\beta d(gy,gv)+ \gamma d(gz,gr)+\delta d(gw,gt)$$
for all $x,y,z,w,u,v,r,t\in X$ with $g(x)\preceq g(u)$, $g(y)\succeq
g(v)$, $g(z)\preceq g(r)$ and $g(w)\succeq g(t)$.
\end{enumerate}

Then $F$ and $g$ have a quadrupled coincidence point (in the sense
of Karapinar and Luong \cite{Q1}), $i.e.$, there exist $x,y,z,w\in
X$ such that
$$F(x,y,z,w)=g(x), F(y,z,w,x)=g(y),$$
$$F(z,w,x,y)=g(z), F(w,x,y,z)=g(w).$$

\noindent{\bf Corollary 37} (Karapinar and Berinde \cite{Q2}).
Corollary 36 remains true if we replace condition (vii) by the
following condition:
\begin{enumerate}
\item [{(vii)$^{\prime}$}] there exists $\varphi\in \Phi$ such that
$$d(F(x,y,z,w),F(u,v,r,t))\leq\varphi\bigg(\frac{d(gx,gu)+d(gy,gv)+d(gz,gr)+d(gw,gt)}
{4}\bigg)$$ for all $x,y,z,w,u,v,r,t\in X$ with $g(x)\preceq g(u)$,
$g(y)\succeq g(v)$, $g(z)\preceq g(r)$ and $g(w)\succeq g(t)$.\\
\end{enumerate}

On setting $n=4$, $\iota_4=\{\{1,3\},\{2,4\}\}$ and
$\ast=\left[\begin{matrix}
1 &4 &3 &2\\
2 &1 &4 &3\\
3 &2 &1 &4\\
4 &3 &2 &1\\
\end{matrix}\right]$ in Corollary 19, we obtain
the following result:\\

\noindent{\bf Corollary 38} (Wu and Liu \cite{TQ}). Let
$(X,d,\preceq)$ be an ordered complete metric space and
$F:X^{4}\rightarrow X$ a mapping. Suppose that the following
conditions hold:
\begin{enumerate}
\item [{(i)}] $F$ has alternating mixed monotone property,
\item [{(ii)}] either $F$ is continuous or $(X,d,\preceq)$ has {\it MCB} property,
\item [{(iii)}] there exist $x^{(0)},y^{(0)},z^{(0)},w^{(0)}
\in X$ such that
$$ x^{(0)}\preceq F(x^{(0)},w^{(0)},z^{(0)},y^{(0)}),$$
$$y^{(0)}\succeq F(y^{(0)},x^{(0)},w^{(0)},z^{(0)}),$$
$$z^{(0)}\preceq F(z^{(0)},y^{(0)},x^{(0)},w^{(0)}),$$
$$w^{(0)}\succeq F(w^{(0)},z^{(0)},y^{(0)},x^{(0)}),$$
\item [{(iv)}] there exist $\alpha,\beta,\gamma,\delta\in [0,1)$ with $\alpha+\beta+\gamma+\delta<1$ such that
$$d(F(x,y,z,w),F(u,v,r,t))\leq\alpha d(x,u)+\beta d(y,v)+ \gamma d(z,r)+\delta d(w,t)$$
for all $x,y,z,w,u,v,r,t\in X$ with $x\preceq u$, $y\succeq v$,
$z\preceq r$ and $w\succeq t$.
\end{enumerate}

Then $F$ has a quadrupled fixed point (in the sense of Wu and Liu
\cite{TQ} ), $i.e.$, there exist $x,y,z,w\in X$ such that
$$F(x,w,z,y)=x, F(y,x,w,z)=y,$$
$$F(z,y,x,w)=z, F(w,z,y,x)=w.$$

On setting $n=4$, $\iota_4=\{\{1,2\},\{3,4\}\}$ and
$\ast=\left[\begin{matrix}
1 &2 &3 &4\\
1 &2 &4 &3\\
3 &4 &2 &1\\
3 &4 &1 &2\\
\end{matrix}\right]$ in Corollary 19, we obtain respectively
the following result:\\

\noindent{\bf Corollary 39} (Berzig and Samet \cite{HD1}). Let
$(X,d,\preceq)$ be an ordered complete metric space and
$F:X^{4}\rightarrow X$ a mapping. Suppose that the following
conditions hold:
\begin{enumerate}
\item [{(i)}] $F$ has 2-mixed monotone property,
\item [{(ii)}] either $F$ is continuous or $(X,d,\preceq)$ has {\it MCB} property,
\item [{(iii)}] there exist $x^{(0)},y^{(0)},z^{(0)},w^{(0)}
\in X$ such that
$$ x^{(0)}\preceq F(x^{(0)},y^{(0)},z^{(0)},w^{(0)}),$$
$$y^{(0)}\preceq F(x^{(0)},y^{(0)},w^{(0)},z^{(0)}),$$
$$z^{(0)}\succeq F(z^{(0)},w^{(0)},y^{(0)},x^{(0)}),$$
$$w^{(0)}\succeq F(z^{(0)},w^{(0)},x^{(0)},y^{(0)}),$$
\item [{(iv)}] there exist $\alpha,\beta,\gamma,\delta\in [0,1)$ with $\alpha+\beta+\gamma+\delta<1$ such that
$$d(F(x,y,z,w),F(u,v,r,t))\leq\alpha d(x,u)+\beta d(y,v)+ \gamma d(z,r)+\delta d(w,t)$$
for all $x,y,z,w,u,v,r,t\in X$ with $x\preceq u$, $y\preceq v$,
$z\succeq r$ and $w\succeq t$.
\end{enumerate}

Then $F$ has a quadrupled fixed point (in the sense of Berzig and
Samet \cite{HD1}), $i.e.$, there exist $x,y,z,w\in X$ such that
$$F(x,y,z,w)=x, F(x,y,w,z)=y,$$
$$F(z,w,y,x)=z, F(z,w,x,y)=w.$$

\subsection{Four fundamental $n$-tupled Coincidence Theorems} \hspace{1cm}\\

In this subsection, we assume $\iota_n=\{O_n, E_n\}$, where
$$O_n=\bigg\{2p-1~:~p\in \bigg\{1,2,...,\displaystyle\bigg[\frac{n+1}{2}\bigg]\bigg\}\bigg\},$$
$i.e.$, the set of all odd natural numbers in $I_n$ and
$$E_n=\bigg\{2p~:~p\in \bigg\{1,2,...,\displaystyle\bigg[\frac{n}{2}\bigg]\bigg\}\bigg\},$$
$i.e.$, the set of all even natural numbers in $I_n.$\\

On setting

\indent\hspace{3cm}$\ast(i,k)=i_k={\begin{cases}i+k-1\;\;\;\;\;\;\;\;\;\;1\leq
k\leq{n-i+1}\cr
\hspace{0.0in}i+k-n-1\;\;\;\;{n-i+2}\leq k\leq{n}\cr\end{cases}}$\\

for even $n$ in Corollaries 4 and 10, we obtain the following
result, which extends the main results of Imdad $et\;al.$ \cite{n1},
Imdad $et\;al.$ \cite{n2}, Husain $et\;al.$ \cite{n3} and Dalal
$et\;al.$ \cite{n4}.\\

\noindent{\bf Corollary 40.} Let $(X,d,\preceq)$ be an ordered
metric space, $E$ an O-complete subspace of $X$ and $n$ an even
natural number. Let $F:X^n\rightarrow X$ and $g: X\rightarrow X$ be
two mappings. Suppose that the following conditions hold:\\
\indent\hspace{0.5mm}$(a)$ $F(X^n)\subseteq g(X)\cap E$,\\
\indent\hspace{0.5mm}$(b)$ $F$ has alternating mixed $g$-monotone property,\\
\indent\hspace{0.5mm}$(c)$ there exist $x^{(0)}_1,x^{(0)}_2,...,x^{(0)}_n \in X$
such that
$${\begin{cases} g(x^{(0)}_{i}) \preceq
F(x^{(0)}_{i},x^{(0)}_{i+1},...,x^{(0)}_n,x^{(0)}_1,...,x^{(0)}_{i-1})\;{\rm
if}~i~{\rm is~odd}\cr g(x^{(0)}_{i}) \succeq
F(x^{(0)}_{i},x^{(0)}_{i+1},...,x^{(0)}_n,x^{(0)}_1,...,x^{(0)}_{i-1})\;{\rm
if}~i~{\rm is~even}\cr\end{cases}}$$ or
$${\begin{cases} g(x^{(0)}_{i}) \succeq
F(x^{(0)}_{i},x^{(0)}_{i+1},...,x^{(0)}_n,x^{(0)}_1,...,x^{(0)}_{i-1})\;{\rm
if}~i~{\rm is~odd}\cr g(x^{(0)}_{i}) \preceq
F(x^{(0)}_{i},x^{(0)}_{i+1},...,x^{(0)}_n,x^{(0)}_1,...,x^{(0)}_{i-1})\;{\rm
if}~i~{\rm is~even,}\cr\end{cases}}$$
\indent\hspace{0.5mm}$(d)$
there exists $\varphi\in \Omega$ such that
$$d(F(x_1,x_2,...,x_n),F(y_1,y_2,...,y_n))\leq\varphi\Big(\frac{1}{n}\sum\limits_{i=1}^{n}d(gx_i,gy_i)\Big)$$
for all $x_{1},x_{2},...,x_{n},y_{1},y_{2},...,y_{n}\in X$ with
[$g(x_i)\preceq g(y_i)$ if $i$ is odd and $g(x_i)\succeq g(y_i)$ if
$i$ is even] or [$g(x_i)\succeq g(y_i)$ if $i$ is odd and
$g(x_i)\preceq g(y_i)$ if $i$ is even],\\
\indent\hspace{0.5mm} $(e)$ $(e1)$ $F$ and $g$ are O-compatible,\\
\indent\hspace{8mm}$(e2)$ $g$ is O-continuous,\\
\indent\hspace{8mm}$(e3)$ either $F$ is O-continuous or $(E,d,\preceq)$ has {\it g-MCB} property\\
\indent\hspace{2mm}or alternately\\
\indent\hspace{0.5mm} $(e^\prime)$ $(e^\prime1)$ $E\subseteq g(X)$,\\
\indent\hspace{9mm}$(e^\prime2)$ either $F$ is $(g,{\rm O})$-continuous or $F$ and $g$ are continuous or $(E,d,\preceq)$\\
\indent\hspace{1.7cm} has {\it MCB} property.\\

Then $F$ and $g$ have a forward cyclic $n$-tupled coincidence
point, $i.e.$, there exist $x_1,x_2,...,x_n\in X$ such that
$$F(x_{i},x_{i+1},...,x_n,x_1,...,x_{i-1})=g(x_i)\;{\rm for\;each}\;i\in I_n.$$

On setting

\indent\hspace{3cm}$\ast(i,k)=i_k={\begin{cases}i+k-1\;\;\;\;\;\;\;\;\;\;1\leq
k\leq{n-i+1}\cr
\hspace{0.0in}i+k-n-1\;\;\;\;{n-i+2}\leq k\leq{n}\cr\end{cases}}$\\

for even $n$ in Corollaries 5 and 11, we obtain the following
result, which extends the main results of Dalal \cite{n5}.\\

\noindent{\bf Corollary 41}. Corollary 40 remains true if we replace
condition $(d)$ by the following condition:
\begin{enumerate}
\item [{$(d)^{\prime}$}] there exists $\varphi\in \Omega$ such that
$$d(F(x_1,x_2,...,x_n),F(y_1,y_2,...,y_n))\leq\varphi\Big(\max\limits_{i\in I_n}d(gx_i,gy_i)\Big)$$
for all $x_{1},x_{2},...,x_{n},y_{1},y_{2},...,y_{n}\in X$ with
[$g(x_i)\preceq g(y_i)$ if $i$ is odd and $g(x_i)\succeq g(y_i)$ if
$i$ is even] or [$g(x_i)\succeq g(y_i)$ if $i$ is odd and
$g(x_i)\preceq g(y_i)$ if $i$ is even].\\
\end{enumerate}

On setting

\indent\hspace{3cm}$\ast(i,k)=i_k={\begin{cases}i-k+1\;\;\;\;\;\;\;\;\;\;1\leq
k\leq{i}\cr
\hspace{0.0in}n+i-k+1\;\;\;\;{i+1}\leq k\leq{n-1}\cr\end{cases}}$\\

for even $n$ in Corollaries 4 and 10 (similarly Corollaries 5 and
11), we obtain the following result:\\

\noindent{\bf Corollary 42}. If in the hypotheses of Corollary 40
(similarly Corollary 41), the condition $(c)$ is replaced by the
following condition:
\begin{enumerate}
\item [{$(c)^{\prime}$}] there exist $x^{(0)}_1,x^{(0)}_2,...,x^{(0)}_n \in X$ such that
$${\begin{cases} g(x^{(0)}_{i}) \preceq
F(x^{(0)}_{i},x^{(0)}_{i-1},...,x^{(0)}_{1},x^{(0)}_n,x^{(0)}_{n-1},...,x^{(0)}_{i+1})\;{\rm
if}~i~{\rm is~odd}\cr g(x^{(0)}_{i}) \succeq
F(x^{(0)}_{i},x^{(0)}_{i-1},...,x^{(0)}_{1},x^{(0)}_n,x^{(0)}_{n-1},...,x^{(0)}_{i+1})\;{\rm
if}~i~{\rm is~even}\cr\end{cases}}$$ or
$${\begin{cases} g(x^{(0)}_{i}) \succeq
F(x^{(0)}_{i},x^{(0)}_{i-1},...,x^{(0)}_{1},x^{(0)}_n,x^{(0)}_{n-1},...,x^{(0)}_{i+1})\;{\rm
if}~i~{\rm is~odd}\cr g(x^{(0)}_{i}) \preceq
F(x^{(0)}_{i},x^{(0)}_{i-1},...,x^{(0)}_{1},x^{(0)}_n,x^{(0)}_{n-1},...,x^{(0)}_{i+1})\;{\rm
if}~i~{\rm is~even,}\cr\end{cases}}$$
\end{enumerate}
then $F$ and $g$ have a backward cyclic $n$-tupled coincidence
point, $i.e.$, there exist $x_1,x_2,\cdots,x_n\in X$ such that
$$F(x_{i},x_{i-1},...,x_{1},x_n,x_{n-1},...,x_{i+1})=g(x_i)\;{\rm for\;each}\;i\in I_n.$$

The following result improves Theorem 2.1 of Karapinar and
Rold$\acute{\rm a}$n \cite{NX4+}.\\

\noindent{\bf Corollary 43}. Corollary 40 (resp. Corollary 41 or Corollary 42) is not valid for any odd natural number $n$.\\
\noindent{\bf Proof.} In view Remark 6, to ensure the existence of $\ast$-fixed
point (for a mapping satisfying $\iota_n$-mixed monotone property),
$\ast\in\mathcal{U}_{\iota_n}$ but in these cases
$\ast\not\in\mathcal{U}_{\iota_n}$. To substantiate this, take
particularly, $n=3$ and $\ast=\left[\begin{matrix}
1 &2 &3\\
2 &3 &1\\
3 &1 &2\\
\end{matrix}\right]$ (in case
of forward cyclic $n$-tupled fixed points). Then
$\ast(2,3)=1\not\in B$, $\ast(3,2)=1\not\in B$, $\ast(3,3)=2\not\in
A$. Similar arguments can be produced in case
of backward cyclic $n$-tupled fixed points.\\

On setting

\indent\hspace{3cm}$\ast(i,k)=i_k={\begin{cases}i-k+1\;\;\;\;\;\;\;\;\;\;1\leq
k\leq{i}\cr
\hspace{0.0in}k-i+1\;\;\;\;{i+1}\leq k\leq{n}\cr\end{cases}}$\\

 in Corollaries 5 and 11, we obtain the following result, which extends
the main results of Gordji and Ramezani \cite{NX1} and Imdad $et\;al.$ \cite{NX2}.\\

\noindent{\bf Corollary 44}. Let $(X,d,\preceq)$ be an ordered
metric space and $E$ an O-complete subspace of $X$. Let
$F:X^n\rightarrow X$ and $g: X\rightarrow X$ be
two mappings. Suppose that the following conditions hold:\\
\indent\hspace{0.5mm}$(a)$ $F(X^n)\subseteq g(X)\cap E$,\\
\indent\hspace{0.5mm}$(b)$ $F$ has alternating mixed $g$-monotone property,\\
\indent\hspace{0.5mm}$(c)$ there exist
$x^{(0)}_1,x^{(0)}_2,...,x^{(0)}_n \in X$ such that
$${\begin{cases} g(x^{(0)}_{i}) \preceq
F(x^{(0)}_{i},x^{(0)}_{i-1},...,x^{(0)}_2,x^{(0)}_1,x^{(0)}_2,...,x^{(0)}_{n-i+1})\;{\rm
if}~i~{\rm is~odd}\cr g(x^{(0)}_{i}) \succeq
F(x^{(0)}_{i},x^{(0)}_{i-1},...,x^{(0)}_2,x^{(0)}_1,x^{(0)}_2,...,x^{(0)}_{n-i+1})\;{\rm
if}~i~{\rm is~even}\cr\end{cases}}$$ or
$${\begin{cases} g(x^{(0)}_{i}) \succeq
F(x^{(0)}_{i},x^{(0)}_{i-1},...,x^{(0)}_2,x^{(0)}_1,x^{(0)}_2,...,x^{(0)}_{n-i+1})\;{\rm
if}~i~{\rm is~odd}\cr g(x^{(0)}_{i}) \preceq
F(x^{(0)}_{i},x^{(0)}_{i-1},...,x^{(0)}_2,x^{(0)}_1,x^{(0)}_2,...,x^{(0)}_{n-i+1})\;{\rm
if}~i~{\rm is~even,}\cr\end{cases}}$$ \indent\hspace{0.5mm}$(d)$
there exists $\varphi\in \Omega$ provided $\varphi$ is increasing
such that
$$d(F(x_1,x_2,...,x_n),F(y_1,y_2,...,y_n))\leq\varphi\Big(\max\limits_{i\in I_n}d(gx_i,gy_i)\Big)$$
for all $x_{1},x_{2},...,x_{n},y_{1},y_{2},...,y_{n}\in X$ with
[$g(x_i)\preceq g(y_i)$ if $i$ is odd and $g(x_i)\succeq g(y_i)$ if
$i$ is even] or [$g(x_i)\succeq g(y_i)$ if $i$ is odd and
$g(x_i)\preceq g(y_i)$ if $i$ is even],\\
\indent\hspace{0.5mm} $(e)$ $(e1)$ $F$ and $g$ are O-compatible,\\
\indent\hspace{8mm}$(e2)$ $g$ is O-continuous,\\
\indent\hspace{8mm}$(e3)$ either $F$ is O-continuous or $(E,d,\preceq)$ has {\it g-MCB} property\\
\indent\hspace{2mm}or alternately\\
\indent\hspace{0.5mm} $(e^\prime)$ $(e^\prime1)$ $E\subseteq g(X)$,\\
\indent\hspace{9mm}$(e^\prime2)$ either $F$ is $(g,{\rm O})$-continuous or $F$ and $g$ are continuous or $(E,d,\preceq)$\\
\indent\hspace{1.7cm} has {\it MCB} property.\\

Then $F$ and $g$ have a 1-skew cyclic $n$-tupled coincidence
point, $i.e.$, there exist $x_1,x_2,\cdots,x_n\in X$ such that
$$F(x_{i},x_{i-1},...,x_2,x_1,x_2,...,x_{n-i+1})=g(x_i)\;{\rm for\;each}\;i\in I_n.$$

On setting

\indent\hspace{3cm}$\ast(i,k)=i_k={\begin{cases}i+k-1\;\;\;\;\;\;\;\;\;\;1\leq
k\leq{n-i+1}\cr \hspace{0.0in}2n-i-k+1\;\;\;\;{n-i+2}\leq
k\leq{n}\cr\end{cases}}$\\

 in Corollaries 5 and 11, we obtain the following result:\\

\noindent{\bf Corollary 45}. If in the hypotheses of Corollary 44,
the condition $(c)$ is replaced by the following condition
\begin{enumerate}
\item [{$(c)^{\prime}$}] there exist $x^{(0)}_1,x^{(0)}_2,...,x^{(0)}_n \in X$ such that
$${\begin{cases} g(x^{(0)}_{i}) \preceq
F(x^{(0)}_{i},x^{(0)}_{i+1},...,x^{(0)}_{n-1},x^{(0)}_n,x^{(0)}_{n-1},...,x^{(0)}_{n-i+1})\;{\rm
if}~i~{\rm is~odd}\cr g(x^{(0)}_{i}) \succeq
F(x^{(0)}_{i},x^{(0)}_{i+1},...,x^{(0)}_{n-1},x^{(0)}_n,x^{(0)}_{n-1},...,x^{(0)}_{n-i+1})\;{\rm
if}~i~{\rm is~even}\cr\end{cases}}$$ or
$${\begin{cases} g(x^{(0)}_{i}) \succeq
F(x^{(0)}_{i},x^{(0)}_{i+1},...,x^{(0)}_{n-1},x^{(0)}_n,x^{(0)}_{n-1},...,x^{(0)}_{n-i+1})\;{\rm
if}~i~{\rm is~odd}\cr g(x^{(0)}_{i}) \preceq
F(x^{(0)}_{i},x^{(0)}_{i+1},...,x^{(0)}_{n-1},x^{(0)}_n,x^{(0)}_{n-1},...,x^{(0)}_{n-i+1})\;{\rm
if}~i~{\rm is~even,}\cr\end{cases}}$$
\end{enumerate}

then $F$ and $g$ have an $n$-skew cyclic $n$-tupled coincidence
point, $i.e.$, there exist $x_1,x_2,\cdots,x_n\in X$ such that
$$F(x_{i},x_{i+1},...,x_{n-1},x_n,x_{n-1},...,x_{n-i+1})=g(x_i)\;{\rm for\;each}\;i\in I_n.$$

\subsection{Berzig-Samet higher dimensional fixed/coincidence point Theorems} \hspace{1cm}\\
On setting $\iota_n=\{\{1,2,....,p\},\{p+1,....,n\}\}$ and
$\ast(i,k)=i_k={\begin{cases}\varphi_i(k)\;\;\;\;\;\;1\leq k\leq
p\cr \hspace{0.0in}\psi_i(k)\;\;\;\;\;\;p< k\leq n\cr\end{cases}}$
(where
$\varphi_1$,...,$\varphi_n$,$\psi_1$,...,$\psi_n$ are arbitrary) in Corollary 19 and Corollary 5, we obtain respectively the following results:\\

\noindent{\bf Corollary 46} (Berzig and Samet \cite{HD1}). Let
$(X,d,\preceq)$ be an ordered complete metric space,
$F:X^{n}\rightarrow X$ a mapping and $p$ a natural number such that
$1\leq p< n$. Let $\varphi_1,...,\varphi_p:\{1,...,p\}\rightarrow
\{1,...,p\}$, $\psi_1,...,\psi_p:\{p+1,...,n\}\rightarrow
\{p+1,...,n\}$, $\varphi_{p+1},...,\varphi_n:\{1,...,p\}\rightarrow
\{p+1,...,n\}$ and $\psi_{p+1},...,\psi_n:\{p+1,...,n\}\rightarrow
\{1,...,p\}$ be $2n$ mappings. Also denote
$x[\varphi(i:i+j)]:=(x_{\varphi(i)},x_{\varphi(i+1)},...,x_{\varphi(i+j)})$.
Suppose that the following conditions hold:
\begin{enumerate}
\item [{(i)}] $F$ has $p$-mixed monotone property,
\item [{(ii)}] either $F$ is continuous or $(X,d,\preceq)$ has {\it MCB} property,
\item [{(iii)}] there exists U$^{(0)}=(x^{(0)}_1,x^{(0)}_2,...,x^{(0)}_n)
\in X^n$ such that
$$ x^{(0)}_1 \preceq F\big(x^{(0)}[\varphi_1(1:p)],x^{(0)}[\psi_1(p+1:n)]\big)$$
$$\vdots$$
$$ x^{(0)}_p \preceq F\big(x^{(0)}[\varphi_p(1:p)],x^{(0)}[\psi_p(p+1:n)]\big)$$
$$\;\;\;\;\;\; x^{(0)}_{p+1} \succeq F\big(x^{(0)}[\varphi_{p+1}(1:p)],x^{(0)}[\psi_{p+1}(p+1:n)]\big)$$
$$\vdots$$
$$\;\;x^{(0)}_{n} \succeq F\big(x^{(0)}[\varphi_{n}(1:p)],x^{(0)}[\psi_{n}(p+1:n)]\big)$$
\item [{(iv)}] there exist $\alpha_i\in [0,1)$ $(1\leq i\leq n)$ with $\sum\limits_{i=1}^{n}\alpha_i<1$ such that
$$d(F({\rm U}),F({\rm V}))\leq\sum\limits_{i=1}^{n}\alpha_id(x_i,y_i)$$
for all U$=(x_1,x_2,...,x_n)$, V$=(y_1,y_2,...,y_n)\in X^n$ with
$$x_1\preceq y_1,...,x_p\preceq y_p,$$
$$x_{p+1}\succeq y_{p+1},...,x_n\succeq y_n.$$
\end{enumerate}

Then there exist $x_1,x_2,...,x_n\in X$ such that
$$F\big(x[\varphi_i(1:p)],x[\psi_i(p+1:n)]\big)=x_i\;{\rm for~each~}i\in I_n.$$

\noindent{\bf Corollary 47} (Aydi and Berzig \cite{HD2}).  Let
$(X,d,\preceq)$ be an ordered complete metric space,
$F:X^{n}\rightarrow X$ and $g:X\to X $ two mappings and $p$ a
natural number such that $1\leq p< n$. Let
$\varphi_1,...,\varphi_p:\{1,...,p\}\rightarrow \{1,...,p\}$,
$\psi_1,...,\psi_p:\{p+1,...,n\}\rightarrow \{p+1,...,n\}$,
$\varphi_{p+1},...,\varphi_n:\{1,...,p\}\rightarrow \{p+1,...,n\}$
and $\psi_{p+1},...,\psi_n:\{p+1,...,n\}\rightarrow \{1,...,p\}$ be
$2n$ mappings. Also denote
$x[\varphi(i:i+j)]:=(x_{\varphi(i)},x_{\varphi(i+1)},...,x_{\varphi(i+j)})$.
Suppose that the following conditions hold:
\begin{enumerate}
\item [{(i)}] $F(X^n)\subseteq g(X)$,
\item [{(ii)}] $F$ has $p$-mixed $g$-monotone property,
\item [{(iii)}] $F$ and $g$ are commuting,
\item [{(iv)}] $g$ is continuous,
\item [{(v)}] either $F$ is continuous or $(X,d,\preceq)$ has $g$-{\it MCB}
property,
\item [{(vi)}] there exists U$^{(0)}=(x^{(0)}_1,x^{(0)}_2,...,x^{(0)}_n)
\in X^n$ such that
$$ g(x^{(0)}_1) \preceq F\big(x^{(0)}[\varphi_1(1:p)],x^{(0)}[\psi_1(p+1:n)]\big)$$
$$\vdots$$
$$ g(x^{(0)}_p) \preceq F\big(x^{(0)}[\varphi_p(1:p)],x^{(0)}[\psi_p(p+1:n)]\big)$$
$$\;\;\;\;\;\; g(x^{(0)}_{p+1}) \succeq F\big(x^{(0)}[\varphi_{p+1}(1:p)],x^{(0)}[\psi_{p+1}(p+1:n)]\big)$$
$$\vdots$$
$$\;\;g(x^{(0)}_{n}) \succeq F\big(x^{(0)}[\varphi_{n}(1:p)],x^{(0)}[\psi_{n}(p+1:n)]\big)$$
\item [{(vii)}] there exists $\varphi\in \Phi$ provided $\varphi$ is increasing
such that
$$d(F({\rm U}),F({\rm V}))\leq\varphi(\max\limits_{i\in I_n} d(gx_i,gy_i))$$
for all U$=(x_1,x_2,...,x_n)$, V$=(y_1,y_2,...,y_n)\in X^n$ with
$$g(x_1)\preceq g(y_1),...,g(x_p)\preceq g(y_p),$$
$$g(x_{p+1})\succeq g(y_{p+1}),...,g(x_n)\succeq g(y_n).$$
\end{enumerate}

Then there exist $x_1,x_2,...,x_n\in X$ such that
$$F\big(x[\varphi_i(1:p)],x[\psi_i(p+1:n)]\big)=g(x_i)\;{\rm for~each~}i\in I_n.$$

\subsection{Rold$\acute{\rm a}$n-Martinez-Moreno-Rold$\acute{\rm a}$n multidimensional Coincidence
Theorems} \hspace{1cm}\\
On setting $\ast(i,k)=i_k=\sigma_i(k)$ (where
$\sigma_1$,$\sigma_2$,...,$\sigma_n$ are arbitrary) in Corollary 7,
we
obtain the following result:\\

\noindent{\bf Corollary 48} (Rold$\acute{\rm a}$n $et\;al.$
\cite{MD1}). Let $(X,d,\preceq)$ be an ordered complete metric space
and $F:X^{n}\rightarrow X$ and $g:X\to X $ two mappings. Let
$\Upsilon=(\sigma_1,\sigma_2,...,\sigma_n)$ be a $n$-tuple of
mappings from $I_n$ into itself verifying $\sigma_i\in\Omega_{A,B}$
if $i\in A$ and $\sigma_i\in\Omega^\prime_{A,B}$ if $i\in B$. Suppose that the following conditions hold:
\begin{enumerate}
\item [{(i)}] $F(X^n)\subseteq g(X)$,
\item [{(ii)}] $F$ has $\iota_n$-mixed $g$-monotone property,
\item [{(iii)}] $F$ and $g$ are commuting,
\item [{(iv)}] $g$ is continuous,
\item [{(v)}] either $F$ is continuous or $(X,d,\preceq)$ has $g$-{\it MCB}
property,
\item [{(vi)}] there exist $x^{(0)}_1,x^{(0)}_2,...,x^{(0)}_n
\in X$ such that
$${\begin{cases} g(x^{(0)}_{i}) \preceq
F(x^{(0)}_{\sigma_i(1)},x^{(0)}_{\sigma_i(2)},...,x^{(0)}_{\sigma_i(n)})\;{\rm
for~ each}~ i\in A\cr g(x^{(0)}_{i}) \succeq
F(x^{(0)}_{\sigma_i(1)},x^{(0)}_{\sigma_i(2)},...,x^{(0)}_{\sigma_i(n)})\;{\rm
for~ each}~ i\in B\cr\end{cases}}$$
\item [{(vii)}] there exists $\alpha\in [0,1)$ such that
$$d(F(x_1,x_2,...,x_n),F(y_1,y_2,...,y_n))\leq\alpha\max\limits_{i\in I_n}d(gx_i,gy_i)$$
for all $x_{1},x_{2},...,x_{n},y_{1},y_{2},...,y_{n}\in X$ with
$g(x_i)\preceq g(y_i)$ for each $i\in A$ and $g(x_i)\succeq g(y_i)$
for each $i\in B$.
\end{enumerate}

Then $F$ and $g$ have, at least, one $\Upsilon$-coincidence point.\\

On setting $\ast(i,k)=i_k=\sigma_i(k)$ (where
$\sigma_1$,$\sigma_2$,...,$\sigma_n$ are arbitrary) in Corollaries
1,2,3,9, we
obtain the following result:\\

\noindent{\bf Corollary 49} (Al-Mezel $et.\;al$ \cite{MD3}). Let
$(X,d,\preceq)$ be an ordered metric space and $F:X^{n}\rightarrow
X$ and $g:X\to X $ two mappings. Let
$\Upsilon=(\sigma_1,\sigma_2,...,\sigma_n)$ be an $n$-tuple of
mappings from $I_n$ into itself verifying $\sigma_i\in\Omega_{A,B}$
if $i\in A$ and $\sigma_i\in\Omega^\prime_{A,B}$ if $i\in B$. Suppose that the following properties are fulfilled:

\begin{enumerate}
\item [{(i)}] $F(X^n)\subseteq g(X)$,
\item [{(ii)}] $F$ has $\iota_n$-mixed $g$-monotone property,
\item [{(iii)}] there exist $x^{(0)}_1,x^{(0)}_2,...,x^{(0)}_n
\in X$ such that
$${\begin{cases} g(x^{(0)}_{i}) \preceq
F(x^{(0)}_{\sigma_i(1)},x^{(0)}_{\sigma_i(2)},...,x^{(0)}_{\sigma_i(n)})\;{\rm
for~ each}~ i\in A\cr g(x^{(0)}_{i}) \succeq
F(x^{(0)}_{\sigma_i(1)},x^{(0)}_{\sigma_i(2)},...,x^{(0)}_{\sigma_i(n)})\;{\rm
for~ each}~ i\in B\cr\end{cases}}$$
\item [{(iv)}] there
exists $\varphi\in \Phi$ such that
$$\frac{1}{n}\sum\limits_{i=1}^{n}d\big(F(x_{\sigma_i(1)},x_{\sigma_i(2)},...,x_{\sigma_i(n)}),F(y_{\sigma_i(1)},y_{\sigma_i(2)},...,y_{\sigma_i(n)})\big)\leq\varphi\Big(\frac{1}{n}\sum\limits_{i=1}^{n}d(gx_i,gy_i)\Big)$$
for all $x_{1},x_{2},...,x_{n},y_{1},y_{2},...,y_{n}\in X$ with
$g(x_i)\preceq g(y_i)$ for each $i\in A$ and $g(x_i)\succeq g(y_i)$
for each $i\in B$.
\end{enumerate}

Also assume that at least one of the following conditions holds:
\begin{enumerate}
\item [{$(a)$}] $(X,d)$ is complete, $F$ and $g$ are continuous and $F$ and $g$ are $({\rm
O},\Upsilon)$-compatible,
\item [{$(b)$}] $(X,d)$ is complete and $F$ and $g$ are continuous and
commuting,
\item [{$(c)$}] $(gX,d)$ is complete and $(X,d,\preceq)$ has {\it MCB}
property,
\item [{$(d)$}] $(X,d)$ is complete, $g(X)$ is closed and $(X,d,\preceq)$ has {\it MCB}
property,
\item [{$(e)$}] $(X,d)$ is complete, $g$ is continuous and increasing, $F$ and $g$ are $({\rm
O},\Upsilon)$-compatible and $(X,d,\preceq)$ has {\it MCB} property.
\end{enumerate}

Then $F$ and $g$ have, at least, one $\Upsilon$-coincidence point.\\

\end{document}